\documentclass[a4paper,12pt]{article}
\usepackage{amsthm}
\usepackage{amsmath}
\allowdisplaybreaks[4]%
\usepackage{bm}%
\usepackage{amssymb}
\usepackage{tabularx}
\usepackage{indentfirst}
\usepackage{enumerate}
\usepackage{dsfont}
\usepackage{graphicx}
\usepackage{subfigure}
\usepackage{epsfig}
\usepackage{graphics}
\usepackage{cases}
\usepackage[compress]{cite}
\usepackage{txfonts}
\usepackage{geometry}
 \usepackage{epstopdf}
 \usepackage{lineno}
 \usepackage[justification=centering]{caption}
 \usepackage{color}
 \usepackage{hyperref}
\hypersetup{
colorlinks=true,
linkcolor=blue,
anchorcolor=blue,
citecolor=blue}
\newtheorem{theorem}{Theorem}[section]
\newtheorem{lemma}{Lemma}[section]

\newtheorem{exm}{Example}[section]

\newtheorem{proposition}{Proposition}[section]
\newtheorem{remark}{Remark}[section]
\newtheorem{definition}{Definition}[section]

\numberwithin{equation}{section}

\begin{document}

\thispagestyle{empty}

\title{Nonlocal diffusion and pulse intervention in a faecal-oral model with moving infected fronts \thanks{The first author is supported by Postgraduate Research \& Practice Innovation Program of Jiangsu Province (KYCX24\_3711) and the third author acknowledges the support of the National Natural Science Foundation of China (No. 12271470).}}

\date{\empty}

\author{Qi Zhou$^1$, Michael Pedersen$^2$, and Zhigui Lin$^1 \thanks{Corresponding author. Email: zglin@yzu.edu.cn (Z. Lin).}$
\\
{\small 1 School of Mathematical Science, Yangzhou University, Yangzhou 225002, China}\\
{\small 2 Department of Applied Mathematics and Computer Science,}\\
{\small Technical University of Denmark, DK 2800, Lyngby, Denmark}
}
 \maketitle
\begin{quote}
\noindent
{\bf Abstract.}{\footnotesize\small~How individual dispersal patterns and human intervention behaviours affect the spread of infectious diseases constitutes a central problem in epidemiological research. This paper develops an impulsive nonlocal faecal-oral model with free boundaries, where pulses are introduced to capture a periodic spraying of disinfectant, and nonlocal diffusion describes the long-range dispersal of individuals, and free boundaries represent moving infected fronts. We first check that the model has a unique nonnegative global classical solution. Then, the principal eigenvalue, which depends on the infected region, the impulse intensity, and the kernel functions for nonlocal diffusion, is examined by using the theory of resolvent positive operators and  their perturbations. Based on this value, this paper obtains that the diseases are either vanishing or spreading, and provides criteria for determining when vanishing  and spreading occur. At the end, a numerical example is presented in order to corroborate the theoretical findings and to gain further understanding of the effect of the pulse intervention. This work shows that the pulsed intervention is beneficial in combating the diseases, but the effect of the nonlocal diffusion depends on the choice of the kernel functions.
}


\end{quote}

\section{Introduction}
Human beings have always been more seriously threatened by diseases than by wars, poverty, disasters, etc. It is estimated that there are about 1.3 to 4 million cases of cholera globally each year, causing about 21,000 to 143,000 deaths ~\cite{Ail-Nelson-Lopze}. Typhoid disease is estimated to cause between 9.9 and 14.7 million infection cases and between 75,000 and 208,000 deaths annually~\cite{Mogasale-Maskery}. Global adult mortality rates significantly increased during the COVID-19 pandemic in 2020 and 2021, reversing past decreasing trends~\cite{Hay-collaborators}. A common feature of the diseases mentioned above is that they can be spread by the faecal-oral
route~\cite{Chen-Tang-Teng, Capasso-Maddalena-JMB, Heller-Mota}. Additionally, poliomyelitis, infectious hepatitis, hand-foot-mouth diseases, and enteric viruses like norovirus, rotavirus and astrovirus can also be transmitted by the faecal-oral route~\cite{Chen-Tang-Teng, Capasso-Maddalena-JMB, Ghosh-Kumar-Santiana}.

Mathematical models have been widely used to predict the trends, design the control measures, and understand the underlying mechanisms of infectious diseases~\cite{Yang-Gong-Sun}. The first faecal-oral model has the form
\begin{eqnarray}\label{Capasso-Paveri-Fontana}
\left\{
\begin{array}{ll}
u'(t)=-a_{11}u+a_{12}v,~~~t>0,\\
v'(t)=-a_{22}v+G(u),~~~t>0\\
\end{array}
\right.
\end{eqnarray}
with suitable initial conditions~\cite{Capasso-Paveri-Fontana}. Model \eqref{Capasso-Paveri-Fontana}, formulated by Capasso and Paveri-Fontana, is used to understand how cholera was transmitted in the Mediterranean regions of Europe in 1973. In \eqref{Capasso-Paveri-Fontana}, $u(t)$ denotes the average density of the infectious agents (bacteria, virus, etc.) at time $t$, while $v(t)$ denotes the average density of the infected individuals at time $t$. The constants $a_{11}$, $a_{12}$, and $a_{22}$ stand for the natural mortality of the infectious agents, the growing rate of the infectious agents stemming from the infected individuals, and the death rate of the infected individuals, respectively. All these constants are positive. $G(u)$ stands for the infective rate of humans, and satisfies
\begin{eqnarray*}
\bf{(G):}
\left\{
\begin{array}{l}
G\in \mathbb{C}^{1}([0,\infty)), ~G(0)=0~\text{and}~G'(u)>0~\text{for}~u\in[0, \infty), \\[2mm]
\frac{G(u)}{u}~\text{is~strongly~decreasing}~\text{for}~u>0~\text{and}~\lim\limits_{u\rightarrow+\infty}\frac{G(u)}{u}<\frac{a_{11}a_{22}}{a_{12}}.
\end{array} \right.
\end{eqnarray*}
Taking into account the spatial movement of individuals, Capasso and Maddalena further developed the corresponding reaction-diffusion model in 1981~\cite{Capasso-Maddalena-JMB}. Subsequently, there have been a series of reaction-diffusion faecal-oral models based on the work~\cite{Capasso-Maddalena-JMB}. We refer to~\cite{Capasso-JMAA} for an integro-differential model,~\cite{Thieme-Zhao-JDE} for an integro-differential model with a time delay,~\cite{Wu-Hsu-TAMS} for an model with two time delays, and~\cite{Hsu-Yang-Nonlinearity} for a model with general growing function.

However, the infected region of infectious disease always moves gradually outwards, but the models mentioned above do not capture this. In 2010, Du and Lin introduced the free boundary condition to study the invasion of a single species~\cite{Du-Lin}. Motivated by this work, Ahn, Baek and Lin developed a diffusive faecal-oral model with free boundaries, where the free boundaries successfully capture the movement of the infected region~\cite{Ahn-Baek-Lin}. This model does not take into account the movement of the infected humans. Considering this point, Wang and Du developed the following model
\begin{eqnarray}\label{Wang-Du-1}
\left\{
\begin{array}{ll}
u_{t}=d_{1}\Delta u-a_{11}u+a_{12}v,\; &\,  x\in(g(t),h(t)), t>0, \\[2mm]
v_{t}=d_{2}\Delta v-a_{22}v+G(u),\; &\,  x\in(g(t),h(t)), t>0,\\[2mm]
u=0, v=0,\; &\,  x\in\{g(t), h(t)\},t>0,\\[2mm]
g'(t)=-\mu_{1}u_{x}(g(t),t)-\mu_{2}v_{x}(g(t),t),\; &\, t>0,\\[2mm]
h'(t)=-\mu_{1}u_{x}(h(t),t)-\mu_{2}v_{x}(h(t),t),\; &\, t>0,\\[2mm]
g(0)=-h_{0}, u(x,0)=u_{0}(x),\; &\,x\in[-h_{0},h_{0}],\\[2mm]
h(0)=h_{0},~~ v(x,0)=v_{0}(x),\; &\,x\in[-h_{0},h_{0}]
\end{array} \right.
\end{eqnarray}
with suitable initial functions $u_{0}(x)$ and $v_{0}(x)$. All constants in \eqref{Wang-Du-1} are positive. The biological meanings of some of the symbols in \eqref{Wang-Du-1}
are listed in \autoref{Tab:01}, and the rest have the same meanings as in \eqref{Capasso-Paveri-Fontana}.
\begin{table}[!htb]
\centering
\caption{List of notations and their meanings in model \eqref{Wang-Du-1}}
\vspace{0.1cm}
\label{Tab:01}
\small
\begin{tabular}{ll}
\hline
 Notation& Biological meaning\\
\hline
$u(x,t)$ & Spatial density of the infectious agents at location $x$ and time $t$\\
$v(x,t)$ & Spatial density of the infectious persons at location $x$ and time $t$\\
$g(t)$ & Left boundary of the infected region at time $t$\\
$h(t)$ & Right boundary of the infected region at time $t$\\
$d_{1}$  & Diffusion coefficient of the infectious agents\\
$d_{2}$  & Diffusion coefficient of the infectious persons\\
$\mu_{1}$   &Expansion capacity of the infectious agents \\
$\mu_{2}$   &Expansion capacity of the infectious persons \\
$h_{0}$   &Length of the right boundary of the initial infection region\\
\hline
\end{tabular}
\end{table}

Note that Fickian diffusion applies to a diffusion process that corresponds to the random walk only when the step size and time size are small compared with the spatial variable and time, respectively~\cite{Zhao-Ruan}. Therefore, the Laplace operator used in model \eqref{Wang-Du-1} only describes the spatial properties of the infectious persons and the infectious agents locally. In 2003, Murray emphasized the importance and intuitively necessity of the long range effects in the biological areas~\cite[Chapter 17]{Murray}. An extensively used nonlocal diffusion operator to replace the local diffusion term is given by
\begin{equation}\label{operator}
(J*u-u)(x,t):=\int_{\mathds{R}}J(x-y)[u(y,t)-u(x,t)]dy,
\end{equation}
where the kernel function $J(x)$ for nonlocal diffusion satisfies
\begin{eqnarray*}
\bf{(J):}
\left\{
\begin{array}{l}
J\in \mathbb{C}(\mathbb{R})\cap\mathbb{L}^{\infty}(\mathbb{R}), ~\text{and}~\int_{-\infty}^{+\infty}J(x)dx=1, \\[2mm]
J(0)>0, ~\text{and}~J(x)=J(-x)\geq 0~\text{for}~x\in\mathbb{R}.
\end{array} \right.
\end{eqnarray*}
Through appropriately selecting $J(x)$, the nonlocal diffusion operator \eqref{operator} may be used to describe both local and nonlocal diffusions. In 2019, Cao et al. introduced and studied a class of free boundary models with nonlocal diffusion in~\cite{Cao-Du-Li-Li}. Based on the work of~\cite{Cao-Du-Li-Li}, there has been extensive research on faecal-oral epidemic models with nonlocal diffusion and free boundaries, see, for instance,~\cite{Wang-Du-1,Zhao-Zhang-Li, Wang-Du-2, Du-Ni,Zhao-Li-Du,Du-Li-Ni,Chang-Du} and the references therein. In particular, Chang and Du considered a nonlocal version of the free boundary model \eqref{Wang-Du-1} in~\cite{Chang-Du}, which has the form
\begin{eqnarray}\label{Wang-Du}
\left\{
\begin{array}{ll}
u_{t}=d_{1}\int\limits^{h(t)}_{g(t)}J_{1}(x-y)u(y,t)dy-d_{1}u-a_{11}u+a_{12}v,\; &\, x\in(g(t),h(t)), t\in\mathbb{R}^{+}, \\[2mm]
v_{t}=d_{2}\int\limits^{h(t)}_{g(t)}J_{2}(x-y)v(y,t)dy-d_{2}v-a_{22}v+G(u),\; &\,  x\in(g(t),h(t)), t\in\mathbb{R}^{+}, \\[2mm]
g'(t)=-\mu_{1}\int\limits^{h(t)}_{g(t)}\int\limits^{g(t)}_{-\infty}J_{1}(x-y)u(x,t)dydx\\
~~~~~~~~~~~~~~~~~~~~~~~~~~~-\mu_{2}\int\limits^{h(t)}_{g(t)}\int\limits^{g(t)}_{-\infty}J_{2}(x-y)v(x,t)dydx,\; &\, t\in\mathbb{R}^{+},\\[2mm]
h'(t)=\mu_{1}\int\limits^{h(t)}_{g(t)}\int\limits^{+\infty}_{h(t)}J_{1}(x-y)u(x,t)dydx\\
~~~~~~~~~~~~~~~~~~~~~~~~~~~+\mu_{2}\int\limits^{h(t)}_{g(t)}\int\limits^{+\infty}_{h(t)}J_{2}(x-y)v(x,t)dydx,\; &\, t\in\mathbb{R}^{+},\\[2mm]
u(g(t),t)=u(h(t),t)=0, v(g(t),t)=v(h(t),t)=0,\; &\, t\in\mathbb{R}^{+},\\[2mm]
u(x,0)=u_{0}(x), v(x,0)=v_{0}(x), -g(0)=h(0)=h_{0}, \; &\,x\in[-h_{0},h_{0}],
\end{array} \right.
\end{eqnarray}
where the kernel functions $J_{i}(i=1,2)$ satisfies the assumption \textbf{(J)}. The biological meanings of the symbols in model \eqref{Wang-Du} are the same as in the local version.

Human intervention can prevent and control the large-scale outbreaks of infectious diseases. For the diseases transmitted by the faecal-oral route, the infectious agents in the environment have a significant impact on the evolution of the diseases. The spraying of disinfectant liquids is effective in order to dramatically reduce the number of  infectious agents rapidly. Classical reaction diffusion models do not capture this instantaneous and abrupt phenomenon. However, this phenomenon can be described well  by impulse differential equation models.

The impulse differential models have the characteristics of both continuous and discrete models. In the last three decades, the theoretical research on impulsive differential equations has yielded numerous results. In the field of biomathematics, Lewis and Li proposed simple impulsive reaction-diffusion equation models to study persistence and spread of species with a reproductive stage and a dispersal stage in bounded and unbounded domains~\cite{Lewis-Li}. Subsequently, some impulsive partial differential equation models rooted in biomathematical background were developed and analyzed. We refer to~\cite{Fazly-Lewis-Wang} for higher dimensional models extended from~\cite{Lewis-Li},~\cite{Zhou-Lin-Santos} for an impulse epidemic model in a periodically evolving environment,~\cite{Liang-Yan-Xiang} for a pest growth model with multiple pulse perturbations,~\cite{Fazly-Lewis-Wang-2} for a hybrid impulsive reaction-advection-diffusion model,~\cite{Wu-Zhao-1} for a birth pulse population model with nonlocal dispersal,~\cite{Zhang-Yi-Chen} for a birth pulse population model with shifting environments, and~\cite{Wu-Zhao-2} for an impulsive hybrid population model in a heterogeneous landscape.

In this paper, we develop a pulsed nonlocal faecal-oral model with free boundaries by introducing the impulsive intervention into model \eqref{Wang-Du}, in order to examine how impulsive intervention and nonlocal diffusion together affect the spread of the diseases. The developed model can be reduced to model \eqref{Wang-Du} by choosing the impulse intervention function properly. The introduction of the impulsive intervention results in the appearance of a periodic solution and reduces the regularity of the solution. As a result, it may make analysis hard and the results become complex and different, and will naturally raise some new problems.

With the introduction of the impulse intervention, it is natural to ask what kind of regularity does the solution possess? whether the spreading-vanishing dichotomy is still true? what are the new criteria for determining spreading and vanishing? whether the impulsive intervention and the nonlocal diffusion influence the movement speed of the infected area? whether the combination of the impulsive intervention and the nonlocal diffusion affects or even alters the dynamical behaviour? These are the research motivations. The main contributions of this research are listed as follows.
\begin{itemize}
\item{ Introduce pulsed intervention into a faecal-oral model with free boundaries and nonlocal diffusion.}
\item{ Investigate the principal eigenvalue, which depends on the impulse intensity, the infected region, and the kernel functions.}
\item{ Prove a vanishing-spreading dichotomy, and provide the sufficient conditions for determining vanishing and spreading.}
\end{itemize}

The paper is arranged as follows. \autoref{Section-2} develops an impulsive faecal-oral model with  free boundaries and  nonlocal diffusion, and present some preliminaries including necessary notations, a frequently used comparison principle, a prior estimate, and the well-posedness of the solution. The related impulse problem with fixed boundaries and the corresponding eigenvalue problem with impulse are treated in \autoref{Section-3}. Here, the existence of the principal eigenvalue is obtained by applying the theory of resolvent positive operators with their perturbations. \autoref{Section-4} proves a vanishing-spreading dichotomy and obtains the sufficient conditions for determining vanishing and spreading based on the results obtained in \autoref{Section-3}. In \autoref{Section-5}, numerical simulation is used to present the effect of the impulsive intervention and to corroborate the validity of the theoretical findings. Finally, a brief discussion and conclusion is given in \autoref{Section-6}.
\section{Model formulation and preliminaries}\label{Section-2}
In this section, we formulate a pulsed nonlocal faecal-oral model in a moving infected region, the global existence of its unique nonnegative classical solution is then given, and finally a frequently used comparison principle is provided.
\subsection{Model formulation}
As stated in the introduction, the transmission of  faecal-oral diseases is usually suppressed by human control measures. As an example,  periodic disinfection of the infected area plays an important role in cutting off the transmission route and protecting the susceptible population. This kind of measure generally results in a dramatic decrease in the spatial density of the infectious agents within a very small time span, and there is no direct impact on the infected individuals.

Motivated by this biological phenomenon, we consider here a situation with two stages: a human intervention stage and a natural developmental stage. In a human intervention stage, the spatial density of the infectious agents decreases impulsively through a discrete-time map. We use $H$ to describe the spatial density at the end of a human intervention stage, which is a function of the density at the beginning of the stage. The spatial density of the infected individuals remains unchanged during this stage. In a natural developmental stage, the evolution of the diseases is described by model \eqref{Wang-Du}.

Let $\tau$ be the time between two adjacent impulses. We always take $k=0,1,2, \cdots$ unless otherwise specified. Denote the start and end times of the pulse by $k\tau$ and $(k\tau)^{+}$, respectively. Then, a mathematical model capturing the spatial dynamics of the diseases transmitted by faecal-oral route under impulsive intervention is given here:
\begin{eqnarray}\label{Zhou-Lin}
\left\{
\begin{array}{ll}
u_{t}=d_{1}\int\limits^{h(t)}_{g(t)}J_{1}(x-y)u(y,t)dy-d_{1}u-a_{11}u+a_{12}v,\; &\,x\in(g(t),h(t)), t\in((k\tau)^{+}, (k+1)\tau],  \\[2mm]
v_{t}=d_{2}\int\limits^{h(t)}_{g(t)}J_{2}(x-y)v(y,t)dy-d_{2}v-a_{22}v+G(u),\; &\, x\in(g(t),h(t)), t\in((k\tau)^{+}, (k+1)\tau],  \\[2mm]
g'(t)=-\mu_{1}\int\limits^{h(t)}_{g(t)}\int\limits^{g(t)}_{-\infty}J_{1}(x-y)u(x,t)dydx\\
~~~~~~~~~~~~~~~~~~~~~~~~~~~-\mu_{2}\int\limits^{h(t)}_{g(t)}\int\limits^{g(t)}_{-\infty}J_{2}(x-y)v(x,t)dydx,\; &\, t\in(k\tau, (k+1)\tau],\\[2mm]
h'(t)=\mu_{1}\int\limits^{h(t)}_{g(t)}\int\limits^{+\infty}_{h(t)}J_{1}(x-y)u(x,t)dydx\\
~~~~~~~~~~~~~~~~~~~~~~~~~~~+\mu_{2}\int\limits^{h(t)}_{g(t)}\int\limits^{+\infty}_{h(t)}J_{2}(x-y)v(x,t)dydx,\; &\, t\in(k\tau, (k+1)\tau],\\[2mm]
u(g(t),t)=u(h(t),t)=0, v(g(t),t)=v(h(t),t)=0,\; &\, t\in(k\tau, (k+1)\tau], \\[2mm]
u(x,(k\tau)^{+})=H(u(x,k\tau)), ~~~~~v(x,(k\tau)^{+})=v(x, k\tau),\; &\, x\in(g(k\tau), h(k\tau)), \\[2mm]
u(x,0)=u_{0}(x), v(x,0)=v_{0}(x), -g(0)=h(0)=h_{0}, \; &\,x\in[-h_{0},h_{0}].
\end{array} \right.
\end{eqnarray}
In \eqref{Zhou-Lin}, it is always assumed that $u(x,t)=v(x,t)=0$ for $x\in\mathbb{R}\setminus [g(t), h(t)]$ and $t\geq 0$, and $u_{0}(x)$ and $v_{0}(x)$ are assumed to satisfy
\begin{eqnarray}\label{initial value}
\left\{
\begin{array}{l}
u_{0}\in\mathbb{C}[-h_{0},h_{0}], ~u_{0}(h_{0})=u_{0}(-h_{0})=0~\text{and}~u_{0}(x)>0~\text{in}~(-h_{0},h_{0}), \\[2mm]
v_{0}\in\mathbb{C}[-h_{0},h_{0}], ~v_{0}(h_{0})=v_{0}(-h_{0})=0~\text{and}~v_{0}(x)>0~~\text{in}~(-h_{0},h_{0}).
\end{array} \right.
\end{eqnarray}
The impulsive function $H(x)$ is assumed to satisfy
\begin{eqnarray*}
\bf{(H):}
\left\{
\begin{array}{l}
H\in \mathbb{C}^{1}([0,\infty)), ~H(0)=0~\text{and}~H'(u)>0~\text{for}~u\in[0, \infty), \\[2mm]
\frac{H(u)}{u}~\text{is~nonincreasing}~\text{and}~0<H(u)/u\leq1~\text{for}~u>0.
\end{array} \right.
\end{eqnarray*}
It should be mentioned that \textbf{(H)} is  natural assumptions for the impulsive function. Some frequently used function forms satisfying assumption \textbf{(H)}
are the linear function $H(u)=c_{1}u$, where $c_{1}\in(0,1]$, and the Beverton-Holt function $H(u)=\frac{c_{2}u}{c_{3}+u}$, where $0<c_{2}\leq c_{3}<+\infty$. In the next content, the assumptions about the initial, kernel, impulsive, and growth functions are always satisfied without specification.

This subsection formulates a pulsed nonlocal epidemic model in a moving infected environment. It is well known that the global existence and uniqueness of the solution to model \eqref{Zhou-Lin} is the basis of this research. This issue will be investigated in the next subsection.
\subsection{Preliminaries}\label{subsection 2-2}
For the sake of later writing, some notations are given first. For any given $\mathbf{u}=(u_{1}, u_{2})$, $\mathbf{v}=(v_{1}, v_{2})$$\in \mathbb{R}^{2}$, define
\begin{equation*}
\mathbf{u}\succ(\prec,\succeq,\preceq,=)\mathbf{v} \text{~if~}u_{i}>(<,\geq,\leq,=)v_{i},~i=1,2.
\end{equation*}
For any given $T$, $\tau$, $h_{0}>0$, there exists a nonnegative integer $N$ such that $N\tau<T\leq(N+1)\tau$, and then define
\begin{equation*}
\begin{split}
\mathbb{G}^{h_{0}}_{T, \tau}=\Bigg\{g(t)\big|g(t)\in\mathbb{C}[0, T], g(0)=-h_{0}, \sup\limits_{n\tau< t_{1}<t_{2}\leq (n+1)\tau} \frac{g(t_{2})-g(t_{1})}{t_{2}-t_{1}}<0, \, \sup\limits_{N\tau< t_{1}<t_{2}\leq T} \frac{g(t_{2})-g(t_{1})}{t_{2}-t_{1}}<0 \Bigg\},\\
\mathbb{H}^{h_{0}}_{T, \tau}=\Bigg\{h(t)\big|h(t)\in\mathbb{C}[0, T], h(0)=~h_{0}, \inf\limits_{n\tau< t_{1}<t_{2}\leq (n+1)\tau} \frac{h(t_{2})-h(t_{1})}{t_{2}-t_{1}}>0, \, \inf\limits_{N\tau< t_{1}<t_{2}\leq T} \frac{h(t_{2})-h(t_{1})}{t_{2}-t_{1}}>0 \Bigg\},
\end{split}
\end{equation*}
where $n=0,1,\cdots, N-1$. For any given $(g,h)\in\mathbb{G}^{h_{0}}_{T, \tau}\times\mathbb{H}^{h_{0}}_{T, \tau}$ and $(u_{0}, v_{0})$ satisfying \eqref{initial value}, we write
\begin{equation*}
\begin{split}
&\Omega_{N\tau}^{T}=\Omega_{N\tau}^{T}(g,h)=\Big\{(x,t) | g(t)<x<h(t), (N\tau)^{+}<t\leq T  \Big\},\\
&\overline{\Omega}_{N\tau}^{T}=\overline{\Omega}_{N\tau}^{T}(g,h)=\Big\{(x,t)|  g(t)\leq x\leq h(t),(N\tau)^{+}<t\leq T \Big\},\\
&\Omega_{n\tau}^{(n+1)\tau}=\Omega_{n\tau}^{(n+1)\tau}(g,h)=\Big\{(x,t)|  g(t)<x<h(t),(n\tau)^{+}<t\leq (n+1)\tau  \Big\},\\
&\overline{\Omega}_{n\tau}^{(n+1)\tau}=\overline{\Omega}_{n\tau}^{(n+1)\tau}(g,h)=\Big\{(x,t)| g(t)\leq x\leq h(t) , (n\tau)^{+}<t\leq (n+1)\tau  \Big\},\\
&\mathbb{X}^{T}=\mathbb{X}^{T}(g,h,u_{0}, v_{0})=\bigg\{(\zeta, \eta)\big |(\zeta, \eta)\in \Big[\mathbb{C}\big(\overline{\Omega}_{n\tau}^{(n+1)\tau}\big)\cap\mathbb{C}\big(\overline{\Omega}_{N\tau}^{T}\big)\Big]^{2},(\zeta(x,0), \eta(x,0))=(u_{0}, v_{0})\text{~in~}[-h_{0}, h_{0}],\\
&~~~~~~~~~~~~~~~~~~~~~~~~~ (\zeta, \eta)\succeq \mathbf{0}\text{~in~}\Omega_{n\tau}^{(n+1)\tau}\cup\Omega_{N\tau}^{T}  \text{,~and~} \zeta(x,t)=\eta(x,t)=0\text{~for~}x\in\{g(t), h(t)\}  \text{~and~} t\in[0, T]\bigg\},
\end{split}
\end{equation*}
where $n=0,1,\cdots, N-1$. Define
\begin{equation}\label{C1C2}
C_{1}:=\max\Big\{u^{*}, \|u_{0}\|_{\infty}, \frac{a_{12}}{a_{11}}\|v_{0}\|_{\infty}\Big\}, ~~C_{2}:=\max\Big\{\|v_{0}\|_{\infty}, \frac{G(C_{1})}{a_{22}}\Big\},
\end{equation}
where $u^{*}$ is determined by $\frac{G(u)}{u}=\frac{a_{11}a_{22}}{a_{12}}$ if $\frac{a_{12}G'(0)}{a_{11}a_{22}}>1$, and $u^{*}=0$ if $\frac{a_{12}G'(0)}{a_{11}a_{22}}\leq1.$ It is not hard to obtain that
\begin{equation}\label{2-2-0}
-a_{11}C_{1}+a_{12}C_{2}\leq 0 \text{~and~}-a_{22}C_{2}+G(C_{1})\leq 0.
\end{equation}

We first prove the following estimate that  has an important function in the proof of the global solution of \eqref{Zhou-Lin}.
\begin{lemma}\label{lemma 2-1}
Assume that $h_{0}$, $\tau>0$, and that $(u_{0}, v_{0})$ satisfies \eqref{initial value}. Then, for any $T>0$ and $(g,h)\in\mathbb{G}^{h_{0}}_{T, \tau}\times\mathbb{H}^{h_{0}}_{T, \tau}$, the following problem
\begin{eqnarray}\label{2-2-1}
\left\{
\begin{array}{ll}
u_{t}= d_{1}\int\limits^{h(t)}_{g(t)}J_{1}(x-y)u(y,t)dy-d_{1}u-a_{11}u+a_{12}v,\; &\, (x,t)\in\Omega_{n\tau}^{(n+1)\tau}\cup\Omega_{N\tau}^{T}, \\[2mm]
v_{t}= d_{2}\int\limits^{h(t)}_{g(t)}J_{2}(x-y)v(y,t)dy-d_{2}v-a_{22}v+G(u),\; &\, (x,t)\in\Omega_{n\tau}^{(n+1)\tau}\cup\Omega_{N\tau}^{T}, \\[2mm]
u(x,(n\tau)^{+})= H(u(x,n\tau)), ~v(x,(n\tau)^{+})= v(x,n\tau),\; &\, x\in[g(n\tau), h(n\tau)], \\[2mm]
u(g(t),t)=u(h(t),t)= 0,u(g(t),t)= v(h(t),t)=0,\; &\, t\in(0, T], \\[2mm]
u(x,0)=u_{0}(x), ~~~~~~~~~~~~~~~~v(x,0)=v_{0}(x), \; &\,x\in[-h_{0},h_{0}]
\end{array} \right.
\end{eqnarray}
has a unique solution $(u,v)\in \mathbb{X}^{T}=\mathbb{X}^{T}(g,h,u_{0}, v_{0})$. Moreover,  $(u,v)$ satisfies
\begin{equation*}
0<u(x,t)\leq C_{1}, ~0<v(x,t)\leq C_{2}\text{~for~}(x,t)\in\Omega_{n\tau}^{(n+1)\tau}\cup\Omega_{N\tau}^{T}.
\end{equation*}
\begin{proof}
Because the initial function $(u_{0}, v_{0})$ satisfies \eqref{initial value} and the impulse function satisfies assumption \textbf{(H)}, it is not hard to derive that
\begin{eqnarray*}
\left\{
\begin{array}{l}
u(x, 0^{+})\in\mathbb{C}[-h_{0},h_{0}], ~u(\pm h_{0},0^{+})=0~\text{and}~u(x,0^{+})>0~\text{in}~(-h_{0},h_{0}), \\[2mm]
v(x, 0^{+})\in\mathbb{C}[-h_{0},h_{0}], ~v(\pm h_{0},0^{+})=0~\text{and}~v(x,0^{+})>0~~\text{in}~(-h_{0},h_{0}),
\end{array} \right.
\end{eqnarray*}
which means that $(u(x, 0^{+}), v(x, 0^{+}))$ satisfies \eqref{initial value} and can be viewed as an initial function. Then, using similar techniques as in \cite[Lemma 5.1]{Wang-Du-1}
yields that problem \eqref{2-2-1} admits a unique solution $(u,v)\in \mathbb{X}^{\tau}$ and
\begin{equation*}
\begin{split}
&0<v(x,t)\leq \max\Big\{\|v(x,0^{+})\|_{\infty}, G(M_{1})/a_{22}\Big\}\text{~for~}(x,t)\in\Omega_{0}^{\tau},\\
0<u(x,&t)\leq\max\Big\{u^{*}, \|u(x,0^{+})\|_{\infty}, \frac{a_{12}}{a_{11}}\|v(x,0^{+})\|_{\infty}\Big\}:=M_{1}\text{~for~}(x,t)\in\Omega_{0}^{\tau}.
\end{split}
\end{equation*}
The hypothesis of $0<H(u)/u\leq1$ in \textbf{(H)} shows that
\begin{equation*}
\max\Big\{u^{*}, \|u(x,0^{+})\|_{\infty}, \frac{a_{12}}{a_{11}}\|v(x,0^{+})\|_{\infty}\Big\}\leq C_{1}\text{~and~}\max\Big\{\|v(x,0^{+})\|_{\infty}, G(M_{1})/a_{22}\Big\}\leq C_{2}
\end{equation*}
which implies that the conclusion holds for $(t,x)\in\Omega_{0}^{\tau}$.

Since $(u(x,\tau^{+}), v(x, \tau^{+}))$ satisfies
\begin{eqnarray*}
\left\{
\begin{array}{l}
u(x,\tau^{+})\in\mathbb{C}[g(\tau),h(\tau)], ~u(g(\tau),\tau^{+})=u(h(\tau),\tau^{+})=0~\text{and}~u( x,\tau^{+})>0~\text{in}~(g(\tau),h(\tau)), \\[2mm]
v(x,\tau^{+})\in\mathbb{C}[g(\tau),h(\tau)], ~v(g(\tau),\tau^{+})=v(h(\tau),\tau^{+})=0~\text{and}~v( x,\tau^{+})>0~~\text{in}~(g(\tau),h(\tau)),
\end{array} \right.
\end{eqnarray*}
it can be viewed as a new initial function for $t\in(\tau^{+}, 2\tau]$. By the same procedures, it follows that problem \eqref{2-2-1}
admits a unique solution $(u,v)$ for $t\in(\tau^{+}, 2\tau]$ and
\begin{equation*}
\begin{split}
&0<v(x,t)\leq \max\Big\{\|v(x,\tau)\|_{\infty}, G(M_{2})/a_{22}\Big\}\text{~for~}(x,t)\in\Omega_{\tau}^{2\tau},\\
0<u(x,&t)\leq\max\Big\{u^{*}, \|u(x,\tau)\|_{\infty}, \frac{a_{12}}{a_{11}}\|v(x,\tau)\|_{\infty}\Big\}:=M_{2}\text{~for~}(x,t)\in\Omega_{\tau}^{2\tau}.
\end{split}
\end{equation*}
Since $(u(x,\tau), v(x,\tau))\preceq ( C_{1}, C_{2})$ in $(g(\tau),h(\tau))$, we have that
\begin{equation*}
\begin{split}
&0<v(x,t)\leq \max\Big\{C_{2}, G(M_{3})/a_{22}\Big\}\text{~for~}(x,t)\in\Omega_{\tau}^{2\tau},\\
0<u(x,&t)\leq\max\Big\{C_{1},a_{12}G(C_{1})/(a_{11}a_{22})\Big\}:=M_{3}\text{~for~}(x,t)\in\Omega_{\tau}^{2\tau}.
\end{split}
\end{equation*}
With the help of \eqref{2-2-0}, it follows that $a_{12}G(C_{1})/(a_{11}a_{22})\leq C_{1}$ and $G(C_{1})/a_{22}\leq C_{2}$. Therefore, the conclusion also holds for $(x,t)\in\Omega_{\tau}^{2\tau}$. Step by step, problem \eqref{2-2-1} has a unique solution $(u,v)\in \mathbb{X}^{T}$ and
\begin{equation*}
\mathbf{0}\prec(u(x,t),v(x,t)) \preceq(C_{1},  C_{2})\text{~for~}(x,t)\in\Omega_{n\tau}^{(n+1)\tau}\cup\Omega_{N\tau}^{T}.
\end{equation*}
This ends the proof.
\end{proof}
\end{lemma}
On the basis of the preceding lemma, we can now establish the following theorem.
\begin{theorem}\label{theorem 2-1}
For any given $h_{0}$, $\tau>0$, the solution of model \eqref{Zhou-Lin}, denoted by $(u,v,g,h)$, exists and is unique for all $t\in(0,\infty)$. Moreover, for any $T>0$, we have that
\begin{equation*}
(u,v,g,h)\in\mathbb{X}^{T}\times\mathbb{G}^{h_{0}}_{T, \tau}\times\mathbb{H}^{h_{0}}_{T, \tau}.
\end{equation*}
\begin{proof}
Because the initial function $(u_{0}, v_{0})$ satisfies \eqref{initial value} and the impulse function satisfies assumption \textbf{(H)}, it is not hard to derive that
\begin{eqnarray*}
\left\{
\begin{array}{l}
u(x,0^{+})\in\mathbb{C}[-h_{0},h_{0}], ~u(\pm h_{0},0^{+})=0~\text{and}~u( x, 0^{+})>0~\text{in}~(-h_{0},h_{0}), \\[2mm]
v(x,0^{+})\in\mathbb{C}[-h_{0},h_{0}], ~v(\pm h_{0}, 0^{+})=0~\text{and}~v( x, 0^{+})>0~~\text{in}~(-h_{0},h_{0}),
\end{array} \right.
\end{eqnarray*}
which means that $(u(x, 0^{+}), v(x, 0^{+}))$ satisfies \eqref{initial value} and can be viewed as an initial function. Then, using similar techniques as in \cite[Theorem 1.1]{Wang-Du-1} obtains that the conclusion holds for $t\in(0^{+}, \tau]$.

Since $(u(x,\tau^{+}), v(x, \tau^{+}))$ satisfies
\begin{eqnarray}\label{XINg}
\left\{
\begin{array}{l}
u(x,\tau^{+})\in\mathbb{C}[g(\tau),h(\tau)], ~u(g(\tau),\tau^{+})=u(h(\tau),\tau^{+})=0~\text{and}~u( x,\tau^{+})>0~\text{in}~(g(\tau),h(\tau)), \\[2mm]
v(x,\tau^{+})\in\mathbb{C}[g(\tau),h(\tau)], ~v(g(\tau),\tau^{+})=v(h(\tau),\tau^{+})=0~\text{and}~v( x,\tau^{+})>0~~\text{in}~(g(\tau),h(\tau)),
\end{array} \right.
\end{eqnarray}
it can be viewed as a new initial function for $t\in(\tau^{+}, 2\tau]$. Then,  the conclusion also holds for $t\in(\tau^{+}, 2\tau]$ with the help of \autoref{lemma 2-1}. By the same procedures, we can find a maximum time $T_{\max}$ such that model \eqref{Zhou-Lin} has a unique solution $(u,v,g,h)$ for $t\in(0^{+}, T_{\max})$ and $(u,v,g,h)\in\mathbb{X}^{T}\times\mathbb{G}^{h_{0}}_{T, \tau}\times\mathbb{H}^{h_{0}}_{T, \tau}$ for
$T\in(0^{+}, T_{\max})$. Next, we use a proof by contradiction to show that $T_{\max}=\infty$. Assume that $T_{\max}<\infty$. This assumption implies that there exists a nonnegative integer $k_{0}$ such that $T_{\max}\in((k_{0}\tau)^{+}, (k_{0}+1)\tau]$. We consider the following two cases.

(1) When $T_{\max}\in((k_{0}\tau)^{+}, (k_{0}+1)\tau)$, $(u(x,(k_{0}\tau)^{+}),v(x, (k_{0}\tau)^{+}))$ is viewed as a new initial vector-valued function. In virtue
of \autoref{lemma 2-1} and the similar proof method
as in \cite[Theorem 1.1]{Wang-Du-1}, one can obtain that the solution of model \eqref{Zhou-Lin}, denoted by $(u,v,g,h)$, exists and is unique for all $t\in(0,(k_{0}+1)\tau]$ and $(u,v,g,h)\in\mathbb{X}^{T}\times\mathbb{G}^{h_{0}}_{T, \tau}\times\mathbb{H}^{h_{0}}_{T, \tau}$ for $T\in(0, (k_{0}+1)\tau]$. This is a contradiction to the definition of $T_{\max}$.

(2) When $T_{\max}=(k_{0}+1)\tau$, we have that \eqref{XINg} with $\tau^{+}$ replaced by $((k_{0}+1)\tau)^{+}$ still holds. Taking $(u(x,((k_{0}+1)\tau)^{+}),v(x,((k_{0}+1)\tau)^{+}))$ to a new initial function and repeating the procedure in (1) yields that the solution of model \eqref{Zhou-Lin}, denoted by $(u,v,g,h)$, exists and is unique for all $t\in(0,(k_{0}+1)\tau]$ and $(u,v,g,h)\in\mathbb{X}^{T}\times\mathbb{G}^{h_{0}}_{T, \tau}\times\mathbb{H}^{h_{0}}_{T, \tau}$ for
$T\in(0, (k_{0}+2)\tau]$. This is also a contradiction to the definition of $T_{\max}$. The proof is completed.
\end{proof}
\end{theorem}
Lastly, we present the following comparison principle, which is used frequently in later analysis.
\begin{lemma}\label{lemma 2-0}
Assume that $(g,h)\in\mathbb{G}^{h_{0}}_{T, \tau}\times\mathbb{H}^{h_{0}}_{T, \tau}$ for some $h_{0}$, $\tau$, $T>0$, and that $c_{12}$, $c_{21}\geq 0$, $u$, $v$, $u_{t}$, $v_{t}\in \mathbb{C}\big(\overline{\Omega}_{n\tau}^{(n+1)\tau}\big)\cap\mathbb{C}\big(\overline{\Omega}_{N\tau}^{T}\big)$, and $c_{ij}\in \mathbb{L}^{\infty}\Big(\overline{\Omega}_{n\tau}^{(n+1)\tau}\cup\overline{\Omega}_{N\tau}^{T}\Big)$ for $i,~j=1,2,$ and $n=0,1,\cdots, N-1$. If
\begin{eqnarray*}
\left\{
\begin{array}{ll}
u_{t}\geq d_{1}\int\limits^{h(t)}_{g(t)}J_{1}(x-y)u(y,t)dy-d_{1}u+c_{11}u+c_{12}v,\; &\, (x,t)\in\Omega_{n\tau}^{(n+1)\tau}\cup\Omega_{N\tau}^{T}, \\[2mm]
v_{t}\geq d_{2}\int\limits^{h(t)}_{g(t)}J_{2}(x-y)v(y,t)dy-d_{2}v+c_{21}u+c_{22}v,\; &\, (x,t)\in\Omega_{n\tau}^{(n+1)\tau}\cup\Omega_{N\tau}^{T}, \\[2mm]
u(x,(n\tau)^{+})\geq H(u(x,n\tau)), ~~v(x,(n\tau)^{+})\geq v(x,n\tau),\; &\, x\in[g(n\tau), h(n\tau)], \\[2mm]
u(g(t),t), u(h(t),t)\geq 0,~~~~ u(g(t),t), v(h(t),t)\geq 0,\; &\, t\in(0, T], \\[2mm]
u(x,0)\geq 0, v(x,0)\geq 0, \; &\,x\in[-h_{0},h_{0}].
\end{array} \right.
\end{eqnarray*}
Then, $(u,v)\succeq \mathbf{0}$ in $\overline{\Omega}_{n\tau}^{(n+1)\tau}\cup\overline{\Omega}_{N\tau}^{T}$. Moreover, if $u(x,0)(v(x,0))\not \equiv 0$ in $[-h_{0}, h_{0}]$, then $u(x,t)(v(x,t))>0$ in ${\Omega}_{n\tau}^{(n+1)\tau}\cup{\Omega}_{N\tau}^{T}$.
\begin{proof}
When $(x,t)\in\Omega_{0}^{\tau}$, it follows from assumption \textbf{(H)} that $(u(x,0^{+}), v(x,0^{+}))\succeq \mathbf{0}$ in $[-h_{0},h_{0}]$. Then treat $(u(x,0^{+}), v(x,0^{+}))$ as an initial function and the conclusion holds by using a similar proof method as in \cite[Lemma 2.1]{Wang-Du-1}. For $(x,t)\in\Omega_{\tau}^{2\tau}$, we take $(u(x,\tau^{+}), v(x,\tau^{+}))$ as a new initial function, which satisfies $(u(x,\tau^{+}), v(x,\tau^{+}))\succeq \mathbf{0}$ in $[g(\tau),h(\tau)]$, and $u(x,\tau^{+})(v(x,\tau^{+}))\not \equiv 0$ in $[g(\tau),h(\tau)]$ as $u(x,0)(v(x,0))\not \equiv 0$ in $[-h_{0}, h_{0}]$. Therefore, the conclusion holds for $(x,t)\in\Omega_{\tau}^{2\tau}$ by using a similar proof method as in \cite[Lemma 2.1]{Wang-Du-1} again. By use of mathematical induction, it follows that $(u,v)\succeq \mathbf{0}$ in $\overline{\Omega}_{n\tau}^{(n+1)\tau}\cup\overline{\Omega}_{N\tau}^{T}$, and $u(x,t)(v(x,t))>0$ in ${\Omega}_{n\tau}^{(n+1)\tau}\cup{\Omega}_{N\tau}^{T}$ when $u(x,0)(v(x,0))\not \equiv 0$ in $[-h_{0}, h_{0}]$. This ends the proof.
\end{proof}
\end{lemma}
\section{\bf The associated problem with fixed boundaries}\label{Section-3}
In this section, we investigate the fixed boundary problem corresponding to model \eqref{Zhou-Lin}, that is model \eqref{Zhou-Lin} with
\begin{equation*}
-g(t)\equiv r>0\text{~and~}h(t)\equiv s>0,
\end{equation*}
which is represented by the following equations:
\begin{eqnarray}\label{Fixed}
\left\{
\begin{array}{ll}
u_{t}=d_{1}\int\limits^{s}_{r}J_{1}(x-y)u(y,t)dy-d_{1}u-a_{11}u+a_{12}v,\; &\,x\in[r,s], t\in((k\tau)^{+}, (k+1)\tau],  \\[2mm]
v_{t}=d_{2}\int\limits^{s}_{r}J_{2}(x-y)v(y,t)dy-d_{2}v-a_{22}v+G(u),\; &\, x\in[r,s],t\in((k\tau)^{+}, (k+1)\tau], \\[2mm]
u(x,(k\tau)^{+})=H(u(x,k\tau)), ~v(x,(k\tau)^{+})=v(x,k\tau),\; &\, x\in[r,s], \\[2mm]
u(x,0)=u_{0}(x), ~~~~~~~~~~~~~~~~~v(x,0)=v_{0}(x), \; &\,x\in[r,s].
\end{array} \right.
\end{eqnarray}

It turns out that the methods required here are quite different from previous works where the impulsive function $H(u)$ in \eqref{Fixed} is replaced by $u$.
A good understanding of the dynamical behaviours of \eqref{Fixed} is remarkably significant in determining the dynamical behaviours of model \eqref{Zhou-Lin}. Therefore, we next present some results for \eqref{Fixed}.
\subsection{The corresponding eigenvalue problem}
There are some works focusing on the existence of principal eigenvalues of the time-periodic cooperative problems with long-range diffusion, see, for instance \cite{Liang-Zhang-Zhao,Bao-Shen, Feng-Li-Ruan-Xin}. Motivated by these works, the existence of the principal eigenvalue of the linearised system corresponding to \eqref{Fixed} will be established next. For the specific definitions of some terms mentioned in this subsection, we refer to~\cite[Section 3.1]{Thieme} for positive perturbation, ~\cite[Definition 4.5]{Thieme-1} for essentially compact perturbation,~\cite[Section 2.1]{Feng-Li-Ruan-Xin} for positive operator, strongly positive operator, resolvent positive operator, spectral bound, real spectral bound, spectral radius, and principal eigenvalue.

It is a well-known fact that the dynamical behaviours of model \eqref{Fixed} are controlled by the principal eigenvalue of its linearised eigenvalue problem at $(0,0)$, which has the following form
\begin{eqnarray}\label{3-1-1}
\left\{
\begin{array}{ll}
\zeta_{t}=d_{1}\int\limits^{s}_{r}J_{1}(x-y)\zeta(y,t)dy-d_{1}\zeta-a_{11}\zeta+a_{12}\eta+\lambda\zeta,\; &\, x\in[r,s],t\in(0^{+}, \tau],  \\[2mm]
\eta_{t}=d_{2}\int\limits^{s}_{r}J_{2}(x-y)\eta(y,t)dy-d_{2}\eta-a_{22}\eta+G'(0)\zeta+\lambda\eta,\; &\, x\in[r,s], t\in(0^{+}, \tau],  \\[2mm]
\zeta(x,0^{+})=H'(0)\zeta(x,0), ~~~\eta(x,0^{+})=\eta(x,0),\; &\, x\in[r,s], \\[2mm]
\zeta(x,0)=\zeta(x,\tau), ~~~~~~~~~~~~~~\eta(x,0)=\eta(x,\tau), \; &\,x\in[r,s].
\end{array} \right.
\end{eqnarray}
To overcome the difficulties presented by the impulse intervention, we consider the following eigenvalue problem
\begin{eqnarray}\label{3-1-2}
\left\{
\begin{array}{ll}
-\frac{\phi_{t}}{\tau}+d_{1}\int\limits^{s}_{r}J_{1}(x-y)\phi(y,t)dy-d_{1}\phi+a_{22}\phi-\frac{\ln H'(0)}{\tau}\phi+a_{12}\varphi=\mu\phi,\; &\,x\in[r,s], t\in(0, 1],  \\[2mm]
-\frac{\varphi_{t}}{\tau}+d_{2}\int\limits^{s}_{r}J_{2}(x-y)\varphi(y,t)dy-d_{2}\varphi+a_{11}\varphi-\frac{\ln H'(0)}{\tau}\varphi+G'(0)\phi=\mu\varphi,\; &\,  x\in[r,s],t\in(0, 1], \\[2mm]
\phi(x,0)=H'(0)\phi(x,1), ~~~\varphi(x,0)=\varphi(x,1),\; &\, x\in[r,s], \\[2mm]
\end{array} \right.
\end{eqnarray}
where $\mu=a_{11}+a_{22}-\frac{\ln H'(0)}{\tau}-\lambda$. With the help of \cite[Section 3]{Zhou-Lin-Santos}, problems \eqref{3-1-2} and \eqref{3-1-1} are equivalent.

For any given $\alpha\in(0,\infty)$, let
\begin{equation*}
\mathcal{X}_{\alpha}=\big\{(\phi, \varphi)\in\mathbb{C}\big([r,s]\times[0,\alpha], \mathbb{R}^{2}\big) : \phi(x,0)=H'(0)\phi(x,\alpha),  \varphi(x,0)=\varphi(x,\alpha),  x\in  [r,s]  \big \}
\end{equation*}
with norm $\|(\phi, \varphi)\|_{\mathcal{X}_{\alpha}}=\sup\limits_{[r,s]\times[0,\alpha]}\sqrt{\phi^{2}+\varphi^{2}}$, and
\begin{equation*}
\begin{split}
\mathcal{X}^{+}_{\alpha}=\big\{(\phi, \varphi)\in\mathcal{X}: (\phi, \varphi)\succeq \mathbf{0},  (x,t)\in[r,s]\times [0,\alpha]\big \},\\
\mathcal{X}^{++}_{\alpha}=\big\{(\phi, \varphi)\in\mathcal{X}: (\phi, \varphi)\succ \mathbf{0},  (x,t)\in[r,s]\times [0,\alpha]\big \}.
\end{split}
\end{equation*}
Let $\mathbb{X}=\mathbb{C}([r,s], \mathbb{R}^{2})$ with norm $\|(\phi, \varphi)\|_{\mathbb{X}}=\sup\limits_{[r,s]}\sqrt{\phi^{2}+\varphi^{2}}$, and
\begin{equation*}
\begin{split}
\mathbb{X}^{+}=\big\{(\phi, \varphi)\in\mathbb{X}: (\phi, \varphi)\succeq \mathbf{0},  x\in [r,s] \big \},\\
\mathbb{X}^{++}=\big\{(\phi, \varphi)\in\mathbb{X}: (\phi, \varphi)\succ \mathbf{0},  x\in  [r,s] \big \}.
\end{split}
\end{equation*}
Now, define two operators $\mathcal{A}$, $\mathcal{B}: \mathcal{X}_{1}\rightarrow\mathcal{X}_{1}$ by
\begin{equation*}
\mathcal{A}[\phi, \varphi](x,t)=
\Bigg(
\begin{matrix}
d_{1}\int\limits^{s}_{r}J_{1}(x-y)\phi(y,t)dy&d_{2}\int\limits^{s}_{r}J_{1}(x-y)\varphi(y,t)dy\\
  \end{matrix}
\Bigg)
\end{equation*}
and
\begin{equation*}
\mathcal{B}[\phi, \varphi](x,t)=
\Bigg(
\begin{matrix}
-\frac{\phi_{t}}{\tau}-d_{1}\phi+a_{22}\phi-\frac{\ln H'(0)}{\tau}\phi+a_{12}\varphi &-\frac{\varphi_{t}}{\tau}-d_{2}\varphi+a_{11}\varphi-\frac{\ln H'(0)}{\tau}\varphi+G'(0)\phi\\
  \end{matrix}
\Bigg)
\end{equation*}
for $(\phi, \varphi)\in \mathcal{X}_{1}$. Then, problem \eqref{3-1-2} can be rewritten as
\begin{equation*}
(\mathcal{A}+\mathcal{B})[\phi, \varphi]:=\mathcal{C}[\phi, \varphi]=\mu (\phi, \varphi).
\end{equation*}
Let $\{\mathcal{G}_{\mu}(t, t_{0})|0\leq t_{0}\leq t\leq 1\}$ and $\{\mathcal{H}_{\mu}(t, t_{0})|0\leq t_{0}\leq t\leq 1\}$ respectively be the evolution family on $\mathbb{X}$
determined by
\begin{eqnarray*}
\left\{
\begin{array}{ll}
\frac{\phi_{t}}{\tau}=-d_{1}\phi+a_{22}\phi-\frac{\ln H'(0)}{\tau}\phi+a_{12}\varphi-\mu\phi,\; &\,x\in[r,s], t\in(0, 1],  \\[2mm]
\frac{\varphi_{t}}{\tau}=-d_{2}\varphi+a_{11}\varphi-\frac{\ln H'(0)}{\tau}\varphi+G'(0)\phi-\mu\varphi,\; &\, x\in[r,s],t\in(0, 1],\\[2mm]
\end{array} \right.
\end{eqnarray*}
and
\begin{eqnarray*}
\left\{
\begin{array}{ll}
\frac{\phi_{t}}{\tau}=d_{1}\int\limits^{s}_{r}J_{1}(x-y)\phi(t,y)dy-d_{1}\phi+a_{22}\phi-\frac{\ln H'(0)}{\tau}\phi+a_{12}\varphi-\mu\phi,\; &\,x\in[r,s], t\in(0, 1],  \\[2mm]
\frac{\varphi_{t}}{\tau}=d_{2}\int\limits^{s}_{r}J_{2}(x-y)\varphi(t,y)dy-d_{2}\varphi+a_{11}\varphi-\frac{\ln H'(0)}{\tau}\phi+G'(0)\phi-\mu\varphi,\; &\, x\in[r,s],t\in(0, 1]. \\[2mm]
\end{array} \right.
\end{eqnarray*}

With the help of \cite[Lemma 2.8]{Feng-Li-Ruan-Xin} (or \cite[Lemma 2.1]{Bao-Shen}), we have the following lemma.
\begin{lemma}\label{lemma 3-1}
The eigenvalue problem $\mathcal{B}[\phi, \varphi]=\eta(\phi, \varphi)$ has a principal eigenvalue, denoted by $\eta_{1}$, which has a unique corresponding
eigenfunction in $\mathcal{X}^{++}_{1}$, denoted by $(\phi_{1}, \varphi_{1})$, with $\|(\phi_{1}, \varphi_{1})\|_{\mathcal{X}_{1}}=1$. Furthermore, we have that
\begin{equation*}
\eta_{1}=\frac{a_{11}+a_{22}-(d_{1}+d_{2}+2\ln H'(0)/\tau)+\sqrt{(d_{1}+d_{2}+2\ln H'(0)/\tau-a_{11}-a_{22})^{2}+4a_{12}G'(0)}}{2},
\end{equation*}
and for any $(x, t)\in[r,s]\times[0,1]$,
\begin{equation*}
\phi_{1}(x, t)\equiv\phi_{1}(t)\text{~and~}\varphi_{1}(x, t)\equiv\varphi_{1}(t).
\end{equation*}
\end{lemma}
Based on the preparations  above, some properties of the operator $\mathcal{B}$ can now be given.
\begin{proposition}\label{proposition 1}
The resolvent operator $(\alpha\mathcal{I}-\mathcal{B})^{-1}$ exists for every $\alpha\in \mathbb{C}$ with Re $\alpha>\eta_{1}$. Furthermore, $\mathcal{B}$ is a resolvent positive operator and the spectral bound $s(\mathcal{B})$ of the operator $\mathcal{B}$ satisfy $s(\mathcal{B})=\eta_{1}$.
\begin{proof}
To begin with, the fundamental theorem of ordinary differential equations gives that for every $\alpha\in \mathbb{C}$ with Re $\alpha>\eta_{1}$, $(\alpha\mathcal{I}-\mathcal{B})^{-1}$ exists in $\mathcal{X}_{1}$. This implies that $s(\mathcal{B})\leq \eta_{1}$. Then, it follows from \cite[Lemma B.3]{Liang-Zhang-Zhao} that the evolution family $\{\mathcal{G}_{0}(t, 0)|t\in(0,1]\}$ generated by $\mathcal{B}$ is positive on $\mathbb{X}_{1}$ for any $t\in(0,1]$. With the help of \cite[Theorem 3.12]{Thieme}, it can be immediately obtained that $\mathcal{B}$ is a resolvent positive operator.

In order to prove that $s(\mathcal{B})=\eta_{1}$, we only need to show that $s(\mathcal{B})\geq \eta_{1}$. Arguing indirectly, assume that $\eta_{1}\in \rho(\mathcal{B})$, where $\rho(\mathcal{B})$ denotes the resolvent set of the operator $\mathcal{B}$. Then $(\mathcal{B}-\eta_{1}\mathcal{I})\mathbf{u}=\mathbf{v}$ has a unique solution $\mathbf{u}\in \mathcal{X}_{1}$ with $\mathbf{u}(x,t)\equiv\mathbf{u}(t)$ for any $\mathbf{v}\in \mathcal{X}_{1}$ with $\mathbf{v}(x,t)\equiv\mathbf{v}(t)$. In virtue of the Fredholm alternative (see e.g., \cite[Theorem 6.6]{Brezis}), it follows that $(\mathcal{B}-\eta_{1}\mathcal{I})\mathbf{u}=\mathbf{0}$ has no nontrivial solution $\mathbf{u}(t)=(u_{1}(t), u_{2}(t))$ with $u_{1}(0)=H'(0)u_{1}(1)$ and $u_{2}(0)=u_{2}(1)$. Recall that $(\eta_{1}, \phi_{1}, \varphi_{1})$ is the principal eigenpair of the eigenvalue problem $\mathcal{B}[\phi, \varphi]=\eta(\phi, \varphi)$. Then, $\mathbf{u}(t)=(\phi_{1}(t),\varphi_{1}(t))$ satisfy $(\mathcal{B}-\eta_{1}\mathcal{I})\mathbf{u}=\mathbf{0}$, which is a contradiction. Thus, $\eta_{1}\in\mathbb{C}- \rho(\mathcal{B})$, and $\eta_{1}\leq s(\mathcal{B})$. This proof is completed.
\end{proof}
\end{proposition}

With the help of the evolution family $\{\mathcal{H}_{\mu}(t, t_{0})|0\leq t_{0}\leq t\leq 1\}$, we define an operator $\mathcal{T}_{\mu}$ by
\begin{equation*}
\mathcal{T}_{\mu}\phi= e^{- \mu\tau}\mathcal{H}_{0}(1, 0)(H'(0)\phi_{1}, \phi_{2})\text{~for~}\phi=(\phi_{1}, \phi_{2})\in\mathbb{X}.
\end{equation*}
\begin{proposition}\label{proposition 2}
There exists $\mu_{0}\in (-\infty, +\infty)$ such that the spectral radius $r(\mathcal{T}_{\mu_{0}})$ of the operator $\mathcal{T}_{\mu_{0}}$ satisfy $r(\mathcal{T}_{\mu_{0}})=1$. Moreover, the operator $\mathcal{C}$ is resolvent positive and $s(\mathcal{C})=\mu_{0}$.
\begin{proof}
We consider the following equation
\begin{equation}\label{3-1-3}
\mathbf{v}=(\mu \mathcal{I}-\mathcal{C})^{-1}\mathbf{u}, ~\mu\in\rho(\mathcal{C}), ~\mathbf{u}=(u_{1}, u_{2})\in\mathcal{X}_{1}.
\end{equation}
With the help of the variation of constant formula, it can be immediately obtained that
\begin{equation}\label{3-1-4}
\mathbf{v}=e^{- \mu\tau t}\mathcal{H}_{0}(t,0)\mathbf{v}(x,0)+\tau\int^{t}_{0} e^{- \mu\tau (t-s)}\mathcal{H}_{0}(t,s)\mathbf{u}(x,s)ds.
\end{equation}
Let $\mathbf{v}=(v_{1}, v_{2})$. Noting that $v_{1}(x,0)=H'(0)v_{1}(x,1)$ and $v_{2}(x,0)=v_{2}(x,1)$ for any $x\in  [r,s]$, it follows from \eqref{3-1-4} that
\begin{equation}\label{3-1-5}
(\mathcal{I}-\mathcal{T}_{\mu})\mathbf{v}(x,0)=\tau\int^{1}_{0} e^{- \mu \tau(1-s)}\mathcal{H}_{0}(1,s)\mathbf{u}(x,s)ds\cdot \text{diag} (H'(0), 1).
\end{equation}
In virtue of \eqref{3-1-3}-\eqref{3-1-5}, it follows that if $1\in \rho(\mathcal{T}_{\mu})$, then
\begin{equation}\label{3-1-6}
\begin{split}
(\mu \mathcal{I}-\mathcal{C})^{-1}\mathbf{u}=&\tau\int^{t}_{0} e^{- \mu\tau (t-s)}\mathcal{H}_{0}(t,s)\mathbf{u}(x,s)ds\\
&+\tau e^{-\mu \tau t}\mathcal{H}_{0}(t,0)(\mathcal{I}-\mathcal{T}_{\mu})^{-1}\int^{1}_{0} e^{- \mu\tau (1-s)}\mathcal{H}_{0}(1,s)\mathbf{u}(x,s)ds\cdot \text{diag} (H'(0), 1).
\end{split}
\end{equation}
Therefore, $1\in \rho(\mathcal{T}_{\mu})$ if and only if $\mu\in \rho(\mathcal{C})$.

A simple calculation yields that the matrix
\begin{equation*}
\mathbf{A}:=
\Bigg[
\begin{matrix}
   a_{22} & a_{12}\\
   G'(0) & a_{11}\\
  \end{matrix}
\Bigg]
\end{equation*}
has a real eigenvalue
\begin{equation*}
\mu_{1}(\mathbf{A})=\frac{a_{11}+a_{22}+\sqrt{(a_{11}-a_{22})^{2}+4a_{12}G'(0)}}{2},
\end{equation*}
which has a positive eigenvector
\begin{equation*}
\mathbf{x}=(a_{12}, \mu_{1}(\mathbf{A})-a_{22})^{\top}.
\end{equation*}
In addition, assumption \textbf{(J)} yields that $(\mu_{1}-\frac{\ln H'(0)}{\tau}, \mathbf{x})$ is an eigenpair of the operator
\begin{equation*}
\mathcal{D}[\phi, \varphi](x):=
\Bigg(
\begin{matrix}
d_{1}\int\limits^{s}_{r}J_{1}(x-y)[\phi(y)-\phi(x)]dy+a_{22}\phi-\frac{\ln H'(0)}{\tau}\phi+a_{12}\varphi\\
d_{2}\int\limits^{s}_{r}J_{1}(x-y)[\varphi(y)-\varphi(x)]dy+a_{11}\varphi-\frac{\ln H'(0)}{\tau}\varphi+G'(0)\phi\\
  \end{matrix}
\Bigg)
\end{equation*}
In virtue of \cite[Lemma 5.8]{Thieme}, one can obtain that
\begin{equation*}
e^{\sigma(\tau \mathcal{D})t}=\sigma(\mathcal{H}_{0}(t, 0))-\{0\} \text{~for~all~}t\in(0,1],
\end{equation*}
where $\sigma(\cdot)$ denotes the spectrum of operator $\cdot$. Recalling that $\mathcal{H}_{0}(t, 0)$ is a positive operator, it follows from \cite[Proposition 4.1.1]{Meyre-Nieberg} that $r(\mathcal{H}_{0}(t, 0))\in \sigma(\mathcal{H}_{0}(t, 0))$. Then, we have that
\begin{equation*}
r(\mathcal{H}_{0}(t, 0))=e^{r(\tau\mathcal{D})t}\geq e^{\tau[\mu_{1}-\ln H'(0)/\tau]t} \text{~for~all~}t\in(0,1],
\end{equation*}
Therefore, we have that
\begin{equation*}
r(\mathcal{T}_{\mu_{1}})\geq H'(0) r(e^{- \mu_{1}\tau}\mathcal{H}_{0}(1, 0))\geq 1.
\end{equation*}
On the other hand, the definition of the operator $\mathcal{T}_{\mu}$ implies that $\lim\limits_{\mu\rightarrow +\infty}r(\mathcal{T}_{\mu})=0$. Since $r(\mathcal{T}_{\mu})$ is strongly decreasing about $\mu$, a unique $\mu_{0}$ can be found such that $r(\mathcal{T}_{\mu_{0}})=1$ and when $\mu>\mu_{0}$, $r(\mathcal{T}_{\mu})<r(\mathcal{T}_{\mu_{0}})=1$. This means that $\rho(\mathcal{C})$ contains a ray $(\mu_{0}, +\infty)$. Additionally, it follows from \eqref{3-1-6} that $(\mu \mathcal{I}-\mathcal{C})^{-1}$ is a positive operator. Hence, $\mathcal{C}$ is a resolvent positive operator. Recalling that $\mathcal{T}_{\mu_{0}}$ is a positive operator, it follows from \cite[Proposition 4.1.1]{Meyre-Nieberg} that $r(\mathcal{T}_{\mu_{0}})=1\in \sigma(\mathcal{T}_{\mu_{0}})$. This implies that $\mu_{0}\in \sigma(\mathcal{C})$  and then $\mu_{0}=s_{\mathbb{R}}(\mathcal{C})$. Finally, \cite[Theorem 3.5]{Thieme-1} gives that $s(\mathcal{C})= s_{\mathbb{R}}(\mathcal{C})=\mu_{0}$. This proof is complected.
\end{proof}
\end{proposition}
Now we present the main theorem of this subsection.
\begin{theorem}\label{theorem 3-1}
The eigenvalue problem \eqref{3-1-1} has a unique eigenvalue $\lambda_{1}$, and the corresponding eigenfunction satisfies $(\zeta_{1}, \eta_{1})\in\mathcal{X}^{++}_{\tau}$. Moreover, we have that
\begin{equation*}
\lambda_{1}=a_{11}+a_{22}-\frac{\ln H'(0)}{\tau}-s(\mathcal{C}).
\end{equation*}
\begin{proof}
In virtue of \autoref{proposition 1}, one can obtain that $\mathcal{B}$ is a resolvent positive operator. Additionally, $\mathcal{A}$ is a positive linear operator. Therefore, $\mathcal{C}=\mathcal{A}+\mathcal{B}$ is a positive perturbation of $\mathcal{B}$. By \autoref{proposition 2}, we now have that the operator $\mathcal{C}$ is resolvent positive. This implies that case (i) in \cite[Theorem 3.6]{Thieme} is impossible.

By standard calculation, the eigenvalue problem
\begin{eqnarray*}
\left\{
\begin{array}{ll}
\frac{\phi_{t}}{\tau}=a_{22}\phi-\frac{\ln H'(0)}{\tau}\phi+a_{12}\varphi,\; &\, x\in[r,s], t\in(0, 1],\\[2mm]
\frac{\varphi_{t}}{\tau}=a_{11}\varphi-\frac{\ln H'(0)}{\tau}\varphi+G'(0)\phi,\; &\, x\in[r,s],t\in(0, 1], \\[2mm]
\phi(x,0)=H'(0)\phi(x,1), ~~~\varphi(x,0)=\varphi(x,1),\; &\, x\in[r,s] \\[2mm]
\end{array} \right.
\end{eqnarray*}
has a principal eigenpair $(\eta_{2}, (\phi_{2}, \varphi_{2}))$. Then $\mathcal{C}[\phi_{2}, \varphi_{2}]=\eta_{2}(\phi_{2}, \varphi_{2})$ and $s(\mathcal{C})\geq \eta_{2}>s(\mathcal{B})$. This implies that case (ii) in \cite[Theorem 3.6]{Thieme} is also impossible. Hence, case (iii) in \cite[Theorem 3.6]{Thieme} will happen. Define
\begin{equation*}
\mathcal{F}_{\mu}=\mathcal{A}(\mu \mathcal{I}-\mathcal{B})^{-1}, ~~\mu>s(\mathcal{B}).
\end{equation*}
The case (iii) in \cite[Theorem 3.6]{Thieme} tells us that there exist $s(\mathcal{B})<\mu_{1}<\mu_{2}$ such that $r(\mathcal{F}_{\mu_{2}})<1\leq r(\mathcal{F}_{\mu_{1}})$.

Let $\{\mathbf{u_{n}}\}\subset \mathcal{X}_{1}\oplus i\mathcal{X}_{1}$ be a bounded sequence. Define $\mathbf{v_{n}}=(\alpha\mathcal{I}-\mathcal{B})^{-1}\mathbf{u_{n}}$, where $\alpha\in \mathbb{C}$ and $\text{Re} \alpha>\eta_{1}$. According to the boundedness of $\mathcal{B}+\frac{\partial _{t}}{\tau}$, both $\mathbf{v_{n}}$ and $\partial _{t}\mathbf{v_{n}}$ are bounded sequences in $\mathcal{X}_{1}\oplus i\mathcal{X}_{1}$. Then, the assumption \textbf{(J)} yields that $\mathcal{A}\mathbf{v_{n}}$ is equicontinuous and uniformly bounded. It follows from the Arzel\`{a}-Ascoli theorem that $\mathcal{A}\mathbf{v_{n}}$ is precompact in $\mathcal{X}$. Therefore, $\mathcal{A}(\alpha\mathcal{I}-\mathcal{B})^{-1}$ is a compact operator in $\mathcal{X}$. Noticing that $(\alpha\mathcal{I}-\mathcal{B})^{-1}$ is a linear bounded operator in $\mathcal{X}$, $(\alpha\mathcal{I}-\mathcal{B})^{-1}\mathcal{A}(\alpha\mathcal{I}-\mathcal{B})^{-1}$ is compact for any $\alpha>s(\mathcal{B})$. This implies that $\mathcal{A}$ is an essentially compact perturbator of $\mathcal{B}$. Now, applying \cite[Theorem 4.7]{Thieme-1} yields that $s(\mathcal{C})$ is a principal eigenvalue of $\mathcal{C}$ with an eigenfunction $(\zeta_{0}, \eta_{0})$ in $\mathcal{X}^{++}_{1}$. By using a similar proof method as in \cite[Theorem 2.3]{Bao-Shen}, it can be obtained that $s(\mathcal{C})$ is an algebraically simple principal eigenvalue.

Next, we prove uniqueness. Let $(\eta, (\phi, \varphi))$ be the principal eigenpairs of the operator $\mathcal{C}$, where $(\zeta, \eta)\in\mathcal{X}^{++}_{1}$. A standard calculation gives that
\begin{equation}\label{3-1-7}
\mathcal{H}_{0}(t, 0)(\zeta(x,0), \eta(x,0))=e^{ \tau\eta t}(\zeta(x,t), \eta(x,t)).
\end{equation}
The strong positivity of the eigenfunction $(\zeta, \eta)$ yields that for any given $(u_{1}(x), u_{2}(x))\in\mathbb{X}^{+}$,
\begin{equation*}
(u_{1}(x), u_{2}(x))\preceq C_{0}(\zeta(x,0), \eta(x,0)):=\frac{\|(u_{1}(x), u_{2}(x))\|_{\mathbb{X}}}{\min\Big[\min\limits_{[r,s]}\zeta(x,0), \min\limits_{[r,s]}\eta(x,0)\Big]}(\zeta(x,0), \eta(x,0)), ~~x\in[r,s].
\end{equation*}
Since the operator $\mathcal{H}_{0}(t, 0)$ is positive, we have that
\begin{equation}\label{3-1-8}
\mathcal{H}_{0}(t, 0)(u_{1}(x), u_{2}(x))\preceq C_{0} \mathcal{H}_{0}(t, 0)(\zeta(x,0), \eta(x,0)), ~~t\in(0,1].
\end{equation}
Now \eqref{3-1-7} together with \eqref{3-1-8} yields that
\begin{equation*}
\mathcal{H}_{0}(t, 0)(u_{1}(x), u_{2}(x))\preceq C_{0} e^{\tau \eta t}(\zeta(x,t), \eta(x,t)), ~~t\in(0,1].
\end{equation*}
This implies that
\begin{equation}\label{3-1-9}
\tau s(\mathcal{C})\leq \omega(\mathcal{H}_{0}(t, 0))\leq  \tau\eta,
\end{equation}
where $\omega(\mathcal{H}_{0}(t, 0))$ denotes the growth bound of the operator $\mathcal{H}_{0}(t, 0)$. Combining \eqref{3-1-9} and the definition of $s(\mathcal{C})$ gives that $\eta=s(\mathcal{C})$. Finally, from the relationship between problems \eqref{3-1-1} and \eqref{3-1-2} we see that the desired conclusion holds. The ends the proof.
\end{proof}
\end{theorem}
\begin{remark}\label{remark 3-1-1}
By the proof of \autoref{theorem 3-1}, it follows that the unique eigenvalue $\lambda_{1}$ of the eigenvalue problem \eqref{3-1-1} is an algebraically simple principal eigenvalue, and its corresponding vector-valued eigenfunction $(\zeta_{1}, \eta_{1})$ is uniquely determined in the sense that the norm is one.
\end{remark}
The theorem above exhibits that the eigenvalue problem \eqref{3-1-1} corresponding to \eqref{Fixed} has a unique principal eigenvalue $\lambda_{1}$.
The next two lemmas will provide a method for estimating the principal eigenvalue
\begin{lemma}\label{lemma 3-1-2}
If there exists a nonnegative vector-valued function $(\zeta,\eta)$ with $\zeta, \eta\not \equiv 0$ and $(\zeta_{t}, \eta_{t})\in\mathbb{C}\big([r,s]\times[0,\tau], \mathbb{R}^{2}\big)$, and a number $\bar{\lambda}$ satisfy
\begin{eqnarray}\label{3-1-10}
\left\{
\begin{array}{ll}
\zeta_{t}\leq d_{1}\int\limits^{s}_{r}J_{1}(x-y)\zeta(y,t)dy-d_{1}\zeta-a_{11}\zeta+a_{12}\eta+\bar{\lambda}\zeta,\; &\, x\in[r,s],t\in(0^{+}, \tau],  \\[2mm]
\eta_{t}\leq d_{2}\int\limits^{s}_{r}J_{2}(x-y)\eta(y,t)dy-d_{2}\eta-a_{22}\eta+G'(0)\zeta+\bar{\lambda}\eta,\; &\,x\in[r,s], t\in(0^{+}, \tau], \\[2mm]
\zeta(x,0^{+})\leq H'(0)\zeta(x,0), ~~~\eta(x,0^{+})\leq \eta(x,0),\; &\, x\in[r,s], \\[2mm]
\zeta(x,0)=\zeta(x,\tau), ~~~~~~~~~~~~~~\eta(x,0)=\eta(x,\tau), \; &\,x\in[r,s],
\end{array} \right.
\end{eqnarray}
then $\lambda_{1}\leq\bar{\lambda}$ holds. Moreover, $\lambda_{1}=\bar{\lambda}$ iff the equalities of \eqref{3-1-10} all hold.
\begin{proof}
We first prove the first assertion. Let
\begin{equation*}
\phi(x,t)=\zeta_{1}(x,\tau-t)\text{~and~}\varphi(x,t)=\eta_{1}( x,\tau-t),
\end{equation*}
where $(\zeta_{1}, \eta_{1})$ is the eigenvector pair corresponding to the unique principal eigenvalue $\lambda_{1}$ of the eigenvalue problem \eqref{3-1-1}. Then, $\phi_{1}$ and $\varphi_{1}$ satisfy the following equation
\begin{eqnarray}\label{3-1-11}
\left\{
\begin{array}{ll}
-\phi_{t}=d_{1}\int\limits^{s}_{r}J_{1}(x-y)\phi(y,t)dy-d_{1}\phi-a_{11}\phi+a_{12}\varphi+\lambda_{1}\phi,\; &\, x\in[r,s],  t\in[0, \tau^{-}),\\[2mm]
-\varphi_{t}=d_{2}\int\limits^{s}_{r}J_{2}(x-y)\varphi(y,t)dy-d_{2}\varphi-a_{22}\varphi+G'(0)\phi+\lambda_{1}\varphi,\; &\, x\in[r,s], t\in[0, \tau^{-}), \\[2mm]
\phi(x,\tau^{-})=H'(0)\phi(x,\tau), ~~~\varphi(x,\tau^{-})=\varphi(x,\tau),\; &\, x\in[r,s], \\[2mm]
\phi(x,0)=\phi(x,\tau), ~~~~~~~~~~~~~~\varphi(x,0)=\varphi(x,\tau), \; &\,x\in[r,s].
\end{array} \right.
\end{eqnarray}
By multiplying the first equations of \eqref{3-1-10} and \eqref{3-1-11} by $\phi$ and $\zeta$, respectively, it follows that
\begin{eqnarray}\label{3-1-12}
\left\{
\begin{array}{ll}
\zeta_{t}\phi \leq d_{1}\int\limits^{s}_{r}J_{1}(x-y)\zeta(y,t)dy\phi-d_{1}\zeta\phi-a_{11}\zeta\phi+a_{12}\eta\phi+\bar{\lambda}\zeta\phi, \\[2mm]
\phi_{t}\zeta=-d_{1}\int\limits^{s}_{r}J_{1}(x-y)\phi(y,t)dy\zeta+d_{1}\phi\zeta+a_{11}\phi\zeta-a_{12}\varphi\zeta-\lambda_{1}\phi\zeta.
\end{array}
\right.
\end{eqnarray}
Then, integrating the two sides of \eqref{3-1-12} on $[r, s]\times(0^{+}, \tau^{-})$ and summing the outcomes yield that
\begin{equation}\label{3-1-13}
\begin{aligned}
\int_{r}^{s}\int^{\tau^{-}}_{0^{+}}(\zeta_{t}\phi+\phi_{t}\zeta)dtdx&\leq(\bar{\lambda}-\lambda_{1})\int_{r}^{s}\int^{\tau^{-}}_{0^{+}}\zeta\phi dtdx+a_{12}\int_{r}^{s}\int^{\tau^{-}}_{0^{+}}(\eta\phi-\varphi\zeta)dtdx\\
 &+ d_{1}\int_{r}^{s}\int^{\tau^{-}}_{0^{+}}\Big[\int\limits^{s}_{r}J_{1}(x-y)\zeta(y,t)dy\phi-\int\limits^{s}_{r}J_{1}(x-y)\phi(y,t)dy\zeta\Big]dtdx.
\end{aligned}
\end{equation}
Recalling that the assumption $J_{1}(x)=J_{1}(-x)\geq 0~\text{for}~x\in\mathds{R}$ in \textbf{(J)}, one can obtain that
\begin{equation}\label{3-1-14}
\begin{aligned}
\int_{r}^{s}\int^{\tau^{-}}_{0^{+}}\Big[\int\limits^{s}_{r}J_{1}(x-y)\zeta(y,t)dy\phi-\int\limits^{s}_{r}J_{1}(x-y)\phi(y,t)dy\zeta\Big]dtdx=0.
\end{aligned}
\end{equation}
Additionally, the impulse and period conditions of \eqref{3-1-10} and \eqref{3-1-11} yield that
\begin{equation}\label{3-1-15}
\begin{aligned}
\int_{r}^{s}\int^{\tau^{-}}_{0^{+}}(\zeta_{t}\phi+\phi_{t}\zeta)dtdx&=\int_{r}^{s}\big[\zeta(x,\tau^{-})\phi( x,\tau^{-})-\zeta( x,0^{+})\phi( x,0^{+})\big]dx\\
&\geq H'(0)\int_{r}^{s}\big[\zeta( x,\tau)\phi( x,0)-\zeta( x,\tau)\phi(x,0)\big]dx\\
&=0.
\end{aligned}
\end{equation}
Inserting  \eqref{3-1-14} and \eqref{3-1-15} into \eqref{3-1-13} yields that
\begin{equation}\label{3-1-16}
\begin{aligned}
\frac{\bar{\lambda}-\lambda_{1}}{a_{12}}\int_{r}^{s}\int^{\tau^{-}}_{0^{+}}\zeta\phi dtdx+\int_{r}^{s}\int^{\tau^{-}}_{0^{+}}(\eta\phi-\varphi\zeta)dtdx\geq 0.
\end{aligned}
\end{equation}
A similar process results in
\begin{equation}\label{3-1-17}
\begin{aligned}
\frac{\bar{\lambda}-\lambda_{1}}{G'(0)}\int_{r}^{s}\int^{\tau^{-}}_{0^{+}}\eta\varphi dtdx+\int_{r}^{s}\int^{\tau^{-}}_{0^{+}}(\zeta\varphi-\phi\eta)dtdx\geq 0.
\end{aligned}
\end{equation}
Adding \eqref{3-1-16} and \eqref{3-1-17} gives that
\begin{equation}\label{3-1-18}
\begin{aligned}
(\bar{\lambda}-\lambda_{1})\int_{r}^{s}\int^{\tau^{-}}_{0^{+}}\Big[\frac{\zeta\phi}{a_{12}} +\frac{\eta\varphi }{G'(0)}\Big]dtdx\geq 0.
\end{aligned}
\end{equation}
Noting that $(\zeta_{1}, \eta_{1})\in\mathcal{X}^{++}_{\tau}$ and $\zeta, \eta\geq(\not \equiv) 0$, it follows from \eqref{3-1-18} that $\bar{\lambda}\geq \lambda_{1}$.

The proof of the remaining conclusion is shown next. When the equalities of \eqref{3-1-10} all hold, the conclusion $\bar{\lambda}= \lambda_{1}$ can be obtained since the eigenvalue of the eigenvalue problem \eqref{3-1-1} is unique. Lastly, when  $\bar{\lambda}= \lambda_{1}$, we prove that the equalities of \eqref{3-1-10} must all hold. Arguing indirectly, assume that there exists an inequality holding in \eqref{3-1-10}. Then, the same method as for the first assertion yields that $\bar{\lambda}> \lambda_{1}$, which is a contradiction to $\bar{\lambda}= \lambda_{1}$. This proof is completed.
\end{proof}
\end{lemma}
\begin{lemma}\label{lemma 3-1-3}
If there exists a nonnegative vector-valued function $(\zeta,\eta)$ with $\zeta, \eta\not \equiv 0$ and $(\zeta_{t}, \eta_{t})\in\mathbb{C}\big([0,\tau]\times[r,s], \mathbb{R}^{2}\big)$, and a number $\underline{\lambda}$ satisfies
\begin{eqnarray}\label{3-1-19}
\left\{
\begin{array}{ll}
\zeta_{t}\geq d_{1}\int\limits^{s}_{r}J_{1}(x-y)\zeta(y,t)dy-d_{1}\zeta-a_{11}\zeta+a_{12}\eta+\underline{\lambda}\zeta,\; &\,x\in[r,s], t\in(0^{+}, \tau], \\[2mm]
\eta_{t}\geq d_{2}\int\limits^{s}_{r}J_{2}(x-y)\eta(y,t)dy-d_{2}\eta-a_{22}\eta+G'(0)\zeta+\underline{\lambda}\eta,\; &\, x\in[r,s],t\in(0^{+}, \tau],\\[2mm]
\zeta(x,0^{+})\geq H'(0)\zeta(x,0), ~~~\eta(x,0^{+})\geq \eta(x,0),\; &\, x\in[r,s], \\[2mm]
\zeta(x,0)=\zeta(x,\tau), ~~~~~~~~~~~~~~\eta(x,0)=\eta(x,\tau), \; &\,x\in[r,s],
\end{array} \right.
\end{eqnarray}
then $\lambda_{1}\geq\underline{\lambda} $ holds. Moreover, $\lambda_{1}=\underline{\lambda} $ iff the equalities of \eqref{3-1-10} all hold.
\begin{proof}
This proof is similar to \autoref{lemma 3-1-2}, and we omit it here.
\end{proof}
\end{lemma}

Subsequently, we will investigate the monotonicity and continuity of the eigenvalue $\lambda_{1}$ of problem \eqref{3-1-1} with respect to the intensity of the impulse function $H$ and the length of the interval $[r,s]$.
\begin{lemma}\label{lemma 3-1-4}
Assume that $[r,s]=[-l, l]$. Let $H'(0)=z$. Denote the eigenvalue of \eqref{3-1-1} by $\lambda_{1}(z)$ and $\lambda_{1}(l)$ in order to emphasise its dependency on $z$ and $l$. Then, we have the following conclusions:
\begin{enumerate}
\item[$(1)$]
$\lambda_{1}(z)$ is strongly decreasing and continuous wrt. $z\in(0, 1]$.
\item[$(2)$]
$\lambda_{1}(l)$ is strongly decreasing and continuous wrt. $l\in(0, \infty)$.
\end{enumerate}
\begin{proof}
(1) First we show that $\lambda_{1}(z)$ is monotone wrt. $z$. For any given $z_{1}, z_{2}\in (0, 1]$ and $z_{1}<z_{2}$, let $(\lambda_{1}(z_{2}), \zeta_{2}, \eta_{2})$ be the principal eigenpair of \eqref{3-1-1}. According to problem \eqref{3-1-1}, we can obtain that
\begin{eqnarray*}
\left\{
\begin{array}{ll}
\frac{\partial\zeta_{2}}{\partial t}=d_{1}\int\limits^{s}_{r}J_{1}(x-y)\zeta_{2}(y,t)dy-d_{1}\zeta_{2}-a_{11}\zeta_{2}+a_{12}\eta_{2}+\lambda(z_{2})\zeta_{2},\; &\,x\in[r,s], t\in(0^{+}, \tau], \\[2mm]
\frac{\partial\eta_{2}}{\partial t}=d_{2}\int\limits^{s}_{r}J_{2}(x-y)\eta_{2}(y,t)dy-d_{2}\eta_{2}-a_{22}\eta_{2}+G'(0)\zeta_{2}+\lambda(z_{2})\eta_{2},\; &\,x\in[r,s], t\in(0^{+}, \tau], \\[2mm]
\zeta_{2}(x,0^{+})>z_{1}\zeta_{2}(x,0), ~~~\eta_{2}(x,0^{+})=\eta_{2}(x,0),\; &\, x\in[r,s], \\[2mm]
\zeta_{2}(x,\tau)=\zeta_{2}(x,\tau), ~~~~~~~~~\eta_{2}(x ,0)=\eta_{2}(x,\tau), \; &\,x\in[r,s].
\end{array} \right.
\end{eqnarray*}
Then, \autoref{lemma 3-1-3} yields that $\lambda(z_{2})<\lambda(z_{1})$.

Subsequently,we  show that $\lambda_{1}(z)$ is continuous wrt.  $z$. To do this, we will prove that for any given $z_{1}, z_{2}\in(0, 1]$, if $z_{1}
\leq z_{2}$, then
\begin{equation}\label{3-1-20}
\lambda_{1}(z_{1})-\lambda_{1}(z_{2})\leq \frac{\ln z_{2}-\ln z_{1}}{\tau},
\end{equation}
while if $z_{1}\geq z_{2}$, then
\begin{equation}\label{3-1-21}
\lambda_{1}(z_{1})-\lambda_{1}(z_{2})\geq \frac{\ln z_{2}-\ln z_{1}}{\tau}.
\end{equation}
Let $(\lambda_{1}(z_{1}), \zeta_{1}, \eta_{1})$ be the principal eigenpair of \eqref{3-1-1}. Now, we construct the following functions
\begin{eqnarray*}
\phi(x,t)=
\left\{
\begin{array}{ll}
z_{1}/z_{2} \zeta_{1}(x,t),\; &\, x\in[r,s],t=0, \\[2mm]
\zeta_{1}(x,t),\; &\, x\in[r,s],t=0^{+}, \\[2mm]
e^{\frac{\ln(z_{1}/z_{2})t}{\tau}}\zeta_{1}(x,t),\; &\, x\in[r,s], t\in(0^{+}, \tau],
\end{array} \right.
\end{eqnarray*}
and
\begin{eqnarray*}
\psi(x,t)=
\left\{
\begin{array}{ll}
z_{1}/z_{2}\eta_{1}(x, t),\; &\, x\in[r,s],t=0, \\[2mm]
\eta_{1}(x, t),\; &\,x\in[r,s], t=0^{+}, \\[2mm]
e^{\frac{\ln(z_{1}/z_{2})t}{\tau}}\eta_{1}(x, t),\; &\,x\in[r,s], t\in(0^{+}, \tau].
\end{array} \right.
\end{eqnarray*}
Then, the functions $\phi$ and $\psi$ satisfy
\begin{eqnarray}\label{3-1-22}
\left\{
\begin{array}{ll}
\frac{\partial\phi}{\partial t}=d_{1}\int\limits^{s}_{r}J_{1}(x-y)\phi(y,t)dy-d_{1}\phi-a_{11}\phi+a_{12}\psi+\big(\lambda(z_{1})+\frac{\ln z_{1}-\ln z_{2}}{\tau}\big)\phi,\; &\, (x,t)\in\Omega^{0}, \\[2mm]
\frac{\partial\psi}{\partial t}=d_{2}\int\limits^{s}_{r}J_{2}(x-y)\psi(y,t)dy-d_{2}\psi-a_{22}\psi+G'(0)\phi+\big(\lambda(z_{1})+\frac{\ln z_{1}-\ln z_{2}}{\tau}\big)\psi,\; &\, (x,t)\in\Omega^{0}, \\[2mm]
\phi(x,0^{+})=z_{2}\phi(x,0), ~~~\psi(x,0^{+})=z_{2}/z_{1}\psi(x,0),\; &\, x\in[r,s], \\[2mm]
\phi(x,0)=\phi(x,\tau), ~~~~~~~~\psi(x,0)=\psi(x,\tau), \; &\,x\in[r,s],
\end{array} \right.
\end{eqnarray}
where $\Omega^{0}=[r,s]\times(0^{+}, \tau]$. When $z_{1}\leq z_{2}$, applying \autoref{lemma 3-1-3} to \eqref{3-1-22} yields that
\begin{equation*}
\lambda_{1}(z_{2})\geq \lambda_{1}(z_{1})+\frac{\ln z_{1}-\ln z_{2}}{\tau},
\end{equation*}
which implies that \eqref{3-1-20} holds. When $z_{1}\geq z_{2}$, applying \autoref{lemma 3-1-2} to \eqref{3-1-22} yields that
\begin{equation*}
\lambda_{1}(z_{2})\leq \lambda_{1}(z_{1})+\frac{\ln z_{1}-\ln z_{2}}{\tau},
\end{equation*}
which implies that \eqref{3-1-21} also holds. We have now proved the desired continuity.

(2) This proof is similar to \cite[Lemma 3.8]{Wang-Du-1}, and we omit it here.
\end{proof}
\end{lemma}

Consider the following eigenvalue problem
\begin{eqnarray}\label{ODE IM}
\left\{
\begin{array}{ll}
\phi_{t}=a_{12}\varphi-a_{11}\phi+\mu\phi,~~~~~~~~~~~~t\in(0^{+}, \tau], \\[2mm]
\varphi_{t}=G'(0)\phi-a_{22}\varphi+\mu\varphi, ~~~~~~~~t\in(0^{+}, \tau], \\[2mm]
\phi(0)=\phi(\tau), \varphi(0)=\varphi(\tau),\\[2mm]
\phi(0^{+})=H'(0)\phi(0), \varphi(0^{+})=\varphi(0).
\end{array} \right.
\end{eqnarray}
By the standard method (see, e.g., the proof of \cite[Lemma 3.3]{Zhou-Lin-Pedersen}), it follows that problem \eqref{ODE IM}
has a unique principal eigenvalue $\mu_{1}$, which has a strongly positive eigenfunction pair $(\phi_{1},\varphi_{1})$.
At the end of this subsection, we will establish a relationship between the principal eigenvalues of \eqref{3-1-1} and \eqref{ODE IM}.

\begin{lemma}\label{lemma 3-1-5}
Assume that $[r,s]=[-l, l]$. Let $\lambda_{1}(l)$ and $\mu_{1}$ denote the principal eigenvalues of \eqref{3-1-1} and \eqref{ODE IM}, respectively. Then, we have that
\begin{equation*}
\lim\limits_{l\rightarrow +\infty}\lambda_{1}(l)=\mu_{1}.
\end{equation*}
\begin{proof}
For any given $l\in \mathds{R}$, we choose $(\zeta, \eta)=(\phi_{1},\varphi_{1})$, where $(\phi_{1},\varphi_{1})$ is the positive eigenfunction pair corresponding to $\mu_{1}$. From assumption \textbf{(J)}, it follows that
\begin{equation*}
\begin{aligned}
\zeta_{t}&= d_{1}\int\limits^{+\infty}_{-\infty}J_{1}(x-y)[\zeta(t)-\zeta(t)]dy-a_{11}\zeta+a_{12}\eta+\mu_{1}\zeta\\
&\geq d_{1}\int\limits^{-l}_{l}J_{1}(x-y)\zeta(t)dy-d_{1}\zeta-a_{11}\zeta+a_{12}\eta+\mu_{1}\zeta,\\
\end{aligned}
\end{equation*}
and
\begin{equation*}
\begin{aligned}
\eta_{t}&= d_{2}\int\limits^{+\infty}_{-\infty}J_{2}(x-y)[\eta(t)-\eta(t)]dy-a_{22}\eta+G'(0)\zeta+\mu_{1}\eta\\
&\geq d_{2}\int\limits^{-l}_{l}J_{2}(x-y)\eta(t)dy-d_{2}\eta-a_{22}\eta+G'(0)\zeta+\mu_{1}\eta.\\
\end{aligned}
\end{equation*}
\autoref{lemma 3-1-3} gives that $\lambda_{1}(l)\geq \mu_{1}$. Then, \autoref{lemma 3-1-4}\textcolor{blue}{(2)} yields that $\lim\limits_{l\rightarrow +\infty}\lambda_{1}(l)\geq\mu_{1}$.

Subsequently, we show that $\mu_{1}=\inf\limits_{l\in \mathds{R}}\lambda_{1}(l)$. To do this, we first define
\begin{equation*}
\theta_{l}(x):=\max\{0, l-|x|\}, ~~~~x\in\mathds{R}.
\end{equation*}
It follows from the proof of \cite[Lemma 3.9(2)]{Wang-Du-1} that for any sufficiently small $\epsilon>0$, there exists a sufficiently large $l^{*}>0$ such that for any $l\geq l^{*}$, we have that
\begin{equation*}
\int\limits_{-l}^{l}J_{i}(x-y)\theta_{l}(y)dy-\theta_{l}(x)>-\frac{\epsilon}{\max\{d_{1}, d_{2}\}}\theta_{l}(x)\text{~for~}x\in[-l,l], i=1,2.
\end{equation*}
For any $l\geq l^{*}$, we choose $(\zeta, \eta)=(\phi_{1}\theta_{l},\varphi_{1}\theta_{l})$. Then, a standard calculation yields that
\begin{equation*}
\begin{aligned}
\zeta_{t}&= d_{1}\zeta-d_{1}\zeta-a_{11}\zeta+a_{12}\eta+\mu_{1}\zeta\\
&< d_{1}\int\limits^{-l}_{l}J_{1}(x-y)\zeta(t)dy-d_{1}\zeta-a_{11}\zeta+a_{12}\eta+(\mu_{1}+\epsilon)\zeta,\\
\end{aligned}
\end{equation*}
and
\begin{equation*}
\begin{aligned}
\eta_{t}&= d_{2}\eta-d_{2}\eta-a_{22}\eta+G'(0)\zeta+\mu_{1}\eta\\
&< d_{2}\int\limits^{-l}_{l}J_{2}(x-y)\eta(t)dy-d_{2}\eta-a_{22}\eta+G'(0)\zeta+(\mu_{1}+\epsilon)\eta.\\
\end{aligned}
\end{equation*}
With the help of \autoref{lemma 3-1-2}, it follows that $\lambda_{1}(l)< \mu_{1}+\epsilon$. This implies that $\mu_{1}=\inf\limits_{l\in \mathds{R}}\lambda_{1}(l)$. This end the proof.
\end{proof}
\end{lemma}
\subsection{The asymptotic profile in a fixed domain }
In this subsection, we still use some notations from \autoref{subsection 2-2} with the free boundaries $g(t)$ and $h(t)$ replaced by the fixed boundaries $r$ and $s$.
By a similar discussion as in \autoref{lemma 2-1}, for any given $\tau>0$, $[r,s]\in \mathbb{R}$, and $(u_{0}(x), v_{0}(x))$ satisfying \eqref{initial value} with $[-h_{0}, h_{0}]$ replaced by $[r,s]$, problem \eqref{Fixed} has a unique nonnegative solution defined for all $t>0$. In this subsection, we always denote this solution by $(u(x,t), v(x,t))$ or $(u, v)$.
Moreover, for any $T>0$, one can obtain that
\begin{equation*}
(u, v)\in \Big[\mathbb{C}\big(\overline{\Omega}_{n\tau}^{(n+1)\tau}\big)\cap\mathbb{C}\big(\overline{\Omega}_{N\tau}^{T}\big)\Big]^{2}
\end{equation*}
and
\begin{equation*}
0<u(x,t)\leq C_{1}, ~0<v(x,t)\leq C_{2}\text{~for~}(x,t)\in\Omega_{n\tau}^{(n+1)\tau}\cup\Omega_{N\tau}^{T},
\end{equation*}
where $C_{1}$ and $C_{2}$ are defined in \eqref{C1C2}. The corresponding steady state problem for \eqref{Fixed} can be written as
\begin{eqnarray}\label{3-2-1}
\left\{
\begin{array}{ll}
U_{t}=d_{1}\int\limits^{s}_{r}J_{1}(x-y)U(y,t)dy-d_{1}U-a_{11}U+a_{12}V,\; &\,x\in[r,s], t\in(0^{+}, \tau], \\[2mm]
V_{t}=d_{2}\int\limits^{s}_{r}J_{2}(x-y)V(y,t)dy-d_{2}V-a_{22}V+G(U),\; &\, x\in[r,s], t\in(0^{+}, \tau], \\[2mm]
U(x,0^{+})=H(U(x,0)), ~~~~~~~~~~~~~V(x,0^{+})=V(x,0),\; &\, x\in[r,s], \\[2mm]
U(x,0)=U(x,\tau), ~~~~~~~~~~~~~~~~~~~~~~V(x,0)=V(x,\tau), \; &\,x\in[r,s].
\end{array} \right.
\end{eqnarray}

Before discussing the asymptotic profile of model \eqref{Fixed}, we first give the definitions of the ordered super- and sub-solutions of problems \eqref{Fixed}
and \eqref{3-2-1}, and then prove that the comparison principles of problems \eqref{Fixed} and \eqref{3-2-1} are also valid. They form the basis of this subsection.
\begin{definition}
For any given $T>0$, nonnegative vector-valued functions $(\overline{u}, \overline{v})$ and $(\underline{u}, \underline{v})$ with
\begin{equation*}
(\overline{u}_{t}, \overline{v}_{t}), (\underline{u}_{t}, \underline{v}_{t}) \in \Big[\mathbb{C}\big(\overline{\Omega}_{n\tau}^{(n+1)\tau}\big)\cap\mathbb{C}\big(\overline{\Omega}_{N\tau}^{T}\big)\Big]^{2}
\end{equation*}
are called the ordered super- and sub-solutions of \eqref{Fixed}, respectively, if $(\overline{u}, \overline{v})$ and $(\underline{u}, \underline{v})$ satisfy
\begin{eqnarray*}
\left\{
\begin{array}{ll}
\overline{u}_{t}\geq d_{1}\int\limits^{s}_{r}J_{1}(x-y)\overline{u}(y,t)dy-d_{1}\overline{u}-a_{11}\overline{u}+a_{12}\overline{v},\; &\, (x,t)\in\Omega_{n\tau}^{(n+1)\tau}\cup\Omega_{N\tau}^{T}, \\[2mm]
\overline{v}_{t}\geq d_{2}\int\limits^{s}_{r}J_{2}(x-y)\overline{v}(y,t)dy-d_{2}\overline{v}-a_{22}\overline{v}+G(\overline{u}),\; &\, (x,t)\in\Omega_{n\tau}^{(n+1)\tau}\cup\Omega_{N\tau}^{T}, \\[2mm]
\overline{u}(x,(n\tau)^{+})\geq H(\overline{u}(x,n\tau)), ~\overline{v}(x,(n\tau)^{+})\geq\overline{v}(x,n\tau),\; &\, x\in[r,s], \\[2mm]
\overline{u}(x,0)\geq u_{0}(x), ~~~~~~~~~~~~~~~~~\overline{v}(x,0)\geq v_{0}(x), \; &\,x\in[r,s],
\end{array} \right.
\end{eqnarray*}
and
\begin{eqnarray*}
\left\{
\begin{array}{ll}
\underline{u}_{t}\leq d_{1}\int\limits^{s}_{r}J_{1}(x-y)\underline{u}(y,t)dy-d_{1}\underline{u}-a_{11}\underline{u}+a_{12}\underline{v},\; &\, (x,t)\in\Omega_{n\tau}^{(n+1)\tau}\cup\Omega_{N\tau}^{T}, \\[2mm]
\underline{v}_{t}\leq d_{2}\int\limits^{s}_{r}J_{2}(x-y)\underline{v}(y,t)dy-d_{2}\underline{v}-a_{22}\underline{v}+G(\underline{u}),\; &\, (x,t)\in\Omega_{n\tau}^{(n+1)\tau}\cup\Omega_{N\tau}^{T}, \\[2mm]
\underline{u}(x,(n\tau)^{+})\leq H(\underline{u}(x,n\tau)), ~\underline{v}(x,(n\tau)^{+})\leq\underline{v}(x,n\tau),\; &\, x\in[r,s], \\[2mm]
\underline{u}(x,0)\leq u_{0}(x), ~~~~~~~~~~~~~~~~~\underline{v}(x,0)\leq v_{0}(x), \; &\,x\in[r,s],
\end{array} \right.
\end{eqnarray*}
respectively.
\end{definition}
\begin{remark}\label{remark 3-2-1}
Denote the ordered super- and sub-solutions of problem \eqref{3-2-1} by the pairs $(\overline{U}, \overline{V})$ and $(\underline{U}, \underline{V})$, respectively. Similarly, an upper (lower) solution of problem \eqref{3-2-1} can also be defined by substituting the initial condition $(\overline{u}(x,0),\overline{v}(x,0))\succeq(u_{0}(x), v_{0}(x))$ $((\underline{u}(x,0),\underline{v}(x,0))\preceq(u_{0}(x), v_{0}(x)))$ with the periodic condition $(\overline{U}(x,0), \overline{V}(x,0))\succeq( \overline{U}(x,\tau), \overline{V}(x,\tau))$ $((\underline{U}(x,0),\underline{V}(x,0))\preceq (\underline{U}(x,\tau), \overline{V}(x,\tau))$.
\end{remark}
\begin{lemma}\label{lemma 3-2-1}
Assume that $r$, $s$, $\tau$, $T>0$ and $c_{12}$, $c_{21}\geq 0$, that $u$, $v$, $u_{t}$, $v_{t}\in \mathbb{C}\big(\overline{\Omega}_{n\tau}^{(n+1)\tau}\big)\cap\mathbb{C}\big(\overline{\Omega}_{N\tau}^{T}\big)$, and that $c_{ij}\in \mathbb{L}^{\infty}\Big(\overline{\Omega}_{n\tau}^{(n+1)\tau}\cup\overline{\Omega}_{N\tau}^{T}\Big)$ for $i,~j=1,2,$ and $n=0,1,\cdots, N-1$. If
\begin{eqnarray*}
\left\{
\begin{array}{ll}
u_{t}\geq d_{1}\int\limits^{s}_{r}J_{1}(x-y)u(y,t)dy-d_{1}u+c_{11}u+c_{12}v,\; &\, (x,t)\in\Omega_{n\tau}^{(n+1)\tau}\cup\Omega_{N\tau}^{T}, \\[2mm]
v_{t}\geq d_{2}\int\limits^{s}_{r}J_{2}(x-y)v(y,t)dy-d_{2}v+c_{21}u+c_{22}v,\; &\, (x,t)\in\Omega_{n\tau}^{(n+1)\tau}\cup\Omega_{N\tau}^{T}, \\[2mm]
u(x,(n\tau)^{+})\geq 0, ~~~~~~~~~~~~~~~~~v(x,(n\tau)^{+})\geq 0,\; &\, x\in[r, s], \\[2mm]
u(x,0)\geq 0, ~~~~~~~~~~~~~~~~~~~~~~~~v(x,0)\geq 0, \; &\,x\in[r,s],
\end{array} \right.
\end{eqnarray*}
then, $(u,v)\succeq \mathbf{0}$ in $\overline{\Omega}_{n\tau}^{(n+1)\tau}\cup\overline{\Omega}_{N\tau}^{T}$. Moreover, if $u(x,0)\not \equiv 0$ $(v(x,0)\not \equiv 0)$ in $[-h_{0}, h_{0}]$, then $u(x,t)>0$ ($v(x,t)>0)$ in ${\Omega}_{n\tau}^{(n+1)\tau}\cup{\Omega}_{N\tau}^{T}$.
\begin{proof}
By the similarity to the proof of \autoref{lemma 2-0}, the desired conclusion can be obtained. Since it is actually much simpler, we omit the details here.
\end{proof}
\end{lemma}
\begin{remark}\label{remark 3-2-2}
Note that the fixed boundary problem \eqref{Fixed} and the corresponding steady state problem \eqref{3-2-1} have a very similar structure. Therefore, when the initial conditions $u(x,0)\geq 0$ and $v(x,0)\geq 0$ are replaced by the periodic conditions $U(x,0)\geq U(x,\tau)$ and $V(x,0)\geq V(x,\tau)$, the conclusions of \autoref{lemma 3-2-1} also hold for problem \eqref{3-2-1}.
\end{remark}
Now, we present the major results of the subsection.
\begin{theorem}\label{theorem 3-2-1}
Let $\lambda_{1}$ be the unique algebraically simple principal eigenvalue of \eqref{3-1-1}. Then the conclusions below hold:
\begin{enumerate}
\item[$(1)$]
If $\lambda_{1}\geq 0$, then $(0,0)$ is the unique solution to problem \eqref{3-2-1}. Moreover, the solution of problem \eqref{Fixed} satisfies
\begin{equation*}
\lim\limits_{t\rightarrow+\infty}(u(x,t),v(x,t))=(0,0) \text{~uniformly~for~}x\in[r,s].
\end{equation*}
\item[$(2)$]
If $\lambda_{1}< 0$, then the problem \eqref{3-2-1} has a unique solution $(U, V)$. Moreover, the solution of problem \eqref{Fixed} satisfies
\begin{equation*}
\lim\limits_{m\rightarrow+\infty}\big(u(x,t+m\tau),v(x,t+m\tau)\big)=\big(U(x,t), V(x,t)\big) \text{~uniformly~for~}(x,t)\in[r,s]\times[0, +\infty).
\end{equation*}
\end{enumerate}
\begin{proof}
(1) We first prove the first assertion. Clearly, $\mathbf{0}$ is a nonnegative solution of problem \eqref{3-2-1}. Next, we prove that it is the only one.
Arguing indirectly, suppose that $(U, V)$ is a strongly positive solution pair of \eqref{3-2-1}. By assumptions \textbf{(G)} and \textbf{(H)}, it follows that
\begin{equation}\label{3-2-2}
(G(U),H(U))\prec(G'(0)U, H'(0)U).
\end{equation}
Substituting \eqref{3-2-2} into \eqref{3-2-1} yields that
\begin{eqnarray}\label{3-2-3}
\left\{
\begin{array}{ll}
U_{t}=d_{1}\int\limits^{s}_{r}J_{1}(x-y)U(y,t)dy-d_{1}U-a_{11}U+a_{12}V,\; &\, x\in[r,s], t\in(0^{+}, \tau],\\[2mm]
V_{t}<d_{2}\int\limits^{s}_{r}J_{2}(x-y)V(y,t)dy-d_{2}V-a_{22}V+G'(0)U,\; &\,x\in[r,s], t\in(0^{+}, \tau], \\[2mm]
U(x,0^{+})<H'(0)U(x,0), ~~~~~~~~~~~~~V(x,0^{+})=V(x,0),\; &\, x\in[r,s], \\[2mm]
U(x,0)=U(x,\tau), ~~~~~~~~~~~~~~~~~~~~~~~~V(x,0)=V(x,\tau), \; &\,x\in[r,s].
\end{array} \right.
\end{eqnarray}
In virtue of \autoref{lemma 3-1-2}, we have from \eqref{3-2-3} that $\lambda_{1}< 0$, which is a contradiction to $\lambda_{1}\geq 0$.

The second assertion is proved next. Its proof is divided into the two situations below.

\textbf{Case 1:} $\lambda_{1}> 0$.

It follows from \autoref{remark 3-1-1} that the unique principal eigenvalue $\lambda_{1}$ has a unique eigenfunction pair $(\zeta_{1}, \eta_{1})$.
With the help of this eigenfunction pair, let
\begin{equation*}
\overline{u}(x,t)=Ke^{-\lambda_{1}t}\zeta_{1}(x,t) \text{~and~}\overline{v}(x,t)=Ke^{-\lambda_{1}t}\eta_{1}(x,t).
\end{equation*}
where $K$ is a sufficiently large positive constant such that $K(\zeta_{1}(x,0), \eta_{1}(x,0))\succeq(u_{0}(x), v_{0}(x))$. By a  standard calculation,
it can be obtained that the function pair $(\overline{u}, \overline{v})$ is a supersolution of \eqref{Fixed}. Now applying \autoref{lemma 3-2-1}
to $(\overline{u}-u, \overline{v}-v)$, it follows that
\begin{equation*}
(u, v)\preceq (\overline{u}, \overline{v}) \text{~for~all~}(x,t)\in [r,s]\times \mathds{R}^{+}.
\end{equation*}
This combined with the condition $\lambda_{1}> 0$ yields that
\begin{equation*}
\lim\limits_{t\rightarrow+\infty}(u(x,t),v(x,t))=(0,0) \text{~uniformly~for~}x\in[r,s].
\end{equation*}

\textbf{Case 2:} $\lambda_{1}=0$.

By a standard calculation, the function pair $(C_{1}, C_{2})$ is a subsolution of \eqref{Fixed}, where $C_{1}$ and $C_{2}$ are defined in \eqref{C1C2}. By selecting $(C_{1}, C_{2})$
as an initial iteration, we can obtain the iteration sequence $\big\{(\overline{u}^{(j)}, \overline{v}^{(j)})\big\}_{j=1}^{+\infty}$ from the following law
\begin{eqnarray}\label{3-2-4}
\left\{
\begin{array}{ll}
\frac{\partial \overline{u}^{(j)}}{\partial t}-d_{1}\int\limits^{s}_{r}J_{1}(x-y)\overline{u}^{(j)}dy+m_{1}\overline{u}^{(j)}=a_{12}\overline{v}^{(i-1)},\; &\, (x,t)\in\Omega^{k}, \\[2mm]
\frac{\partial \overline{v}^{(j)}}{\partial t}-d_{2}\int\limits^{s}_{r}J_{1}(x-y)\overline{v}^{(j)}dy+m_{2}\overline{v}^{(j)}=G(\overline{u}^{(i-1)}),\; &\, (x,t)\in\Omega^{k}, \\[2mm]
\overline{u}^{(j)}(x,0)=\overline{u}^{(j-1)}(x,\tau), \overline{u}^{(j)}(x,(k\tau)^{+})=H(\overline{u}^{(j-1)}(x,(k+1)\tau)),\; &\,x\in[r,s],\\[2mm]
\overline{v}^{(j)}(x,0)=\overline{v}^{(j-1)}(x,\tau), \overline{v}^{(j)}(x,(k\tau)^{+})=\overline{v}^{(j-1)}(x,(k+1)\tau), \; &\, x\in[r,s],
\end{array} \right.
\end{eqnarray}
where $\Omega^{k}=[r,s]\times ((k\tau)^{+}, (k+1)\tau]$, $m_{1}=d_{1}+a_{11}$, and $m_{2}=d_{2}+a_{22}$. Next, we prove that the obtained sequence $\big\{(\overline{u}^{(j)}, \overline{v}^{(j)})\big\}_{j=0}^{+\infty}$ is monotonically decreasing. Let
\begin{equation*}
\phi^{(j)}=\overline{u}^{(j)}-\overline{u}^{(j+1)}\text{~and~}\psi^{(j)}=\overline{v}^{(j)}-\overline{v}^{(j+1)}.
\end{equation*}
Then, the initial conditions satisfy
\begin{eqnarray*}
\left\{
\begin{array}{ll}
\phi^{(0)}(x,0)=\overline{u}^{(0)}(x,0)-\overline{u}^{(1)}(x,0)=\overline{u}^{(0)}(x,0)-\overline{u}^{(0)}(x,\tau)=0,~~x\in[r,s],\\
\psi^{(0)}(x, 0)=\overline{v}^{(0)}(x,0)-\overline{v}^{(1)}(x,0)=\overline{v}^{(0)}(x,0)-\overline{v}^{(0)}(x,\tau)=0,~~~x\in[r,s].
\end{array} \right.
\end{eqnarray*}
The impulse conditions satisfy
\begin{eqnarray*}
\left\{
\begin{array}{ll}
\phi^{(0)}(x,0^{+})=\overline{u}^{(0)}(x,0^{+})-\overline{u}^{(1)}(x,0^{+})=\overline{u}^{(0)}(x,0^{+})-H(\overline{u}^{(0)}(x,\tau))\geq 0,~~x\in[r,s],\\
\psi^{(0)}(x,0^{+})=\overline{v}^{(0)}(x,0^{+})-\overline{v}^{(1)}(x,0^{+})=\overline{v}^{(0)}(x,0^{+})-\overline{v}^{(0)}(x,\tau)=0,~~~~~~~~x\in[r,s].
\end{array} \right.
\end{eqnarray*}
The equations satisfy
\begin{eqnarray*}
\left\{
\begin{array}{ll}
\phi^{(0)}_{t}\geq d_{1}\int\limits^{s}_{r}J_{1}(x-y)\phi^{(0)}(y,t)dy-d_{1}\phi^{(0)}-a_{11}\phi^{(0)},\; &\, (x,t)\in\Omega^{0}, \\[2mm]
\psi^{(0)}_{t}\geq d_{2}\int\limits^{s}_{r}J_{1}(x-y)\psi^{(0)}(y,t)dy-d_{2}\psi^{(0)}-a_{22}\phi^{(0)},\; &\, (x,t)\in\Omega^{0}. \\[2mm]
\end{array} \right.
\end{eqnarray*}
Based on these, \autoref{lemma 3-2-1} yields that $(\phi^{(0)}, \psi^{(0)})\succeq \mathbf{0}$ for $x\in[r,s]$ and $t\in[0,\tau]$. Taking $(\phi^{(0)}(x,\tau^{+}), \psi^{(0)}\\
(x,\tau^{+}))$ as a fresh start function for $t\in(\tau^{+}, 2\tau]$, a similar discussion yields that $(\phi^{(0)}, \psi^{(0)})\succeq \mathbf{0}$ for $(x,t)\in[r,s]\times(\tau, 2\tau]$.
From the method of induction, $(\phi^{(0)}, \psi^{(0)})\succeq \mathbf{0}$ for $(x,t)\in\Omega^{k}$. Using the method of induction again,
the sequence $(\phi^{(j)}, \psi^{(j)})\succeq \mathbf{0}$ for $(x,t)\in\Omega^{k}$, which implies that $\big\{(\overline{u}^{(j)}, \overline{v}^{(j)})\big\}_{j=0}^{+\infty}$
is monotonically decreasing. Noting that $(\overline{u}^{(0)}, \overline{v}^{(0)})$ is a subsolution of \eqref{Fixed}, applying \autoref{lemma 3-2-1} to $(\overline{u}^{(0)}-u, \overline{v}^{(0)}-v)$ yields that
\begin{equation}\label{3-2-5}
(u, v)\preceq (\overline{u}^{(0)}, \overline{v}^{(0)}) \text{~for~all~}(x,t)\in \Omega^{k}.
\end{equation}
Then, it follows from \eqref{3-2-5} that
\begin{equation}\label{3-2-6}
(u(x,\tau), v(x,\tau))\preceq (\overline{u}^{(0)}(x,\tau), \overline{v}^{(0)}(x,\tau)) \text{~for~all~}x\in[r,s].
\end{equation}
This together with $(\overline{u}^{(0)}(x,\tau), \overline{v}^{(0)}(x,\tau))=(\overline{u}^{(1)}(x,0), \overline{v}^{(1)}(x,0))$ yields that
\begin{equation}\label{3-2-7}
(u(x,\tau), v(x,\tau))\preceq (\overline{u}^{(1)}(x,0), \overline{v}^{(1)}(x,0)) \text{~for~all~}x\in[r,s].
\end{equation}
Recalling that $H$ is a monotonically increasing function with respect to $u$, it follows from \eqref{3-2-6}, the third equation of \eqref{Fixed},
and the last two equations of \eqref{3-2-3} that
\begin{equation}\label{3-2-8}
(u(x,\tau^{+}), v(x,\tau^{+}))\preceq (H(u(x,\tau)), v(x,\tau))\preceq (H(\overline{u}^{(0)}(x,\tau)), \overline{v}^{(0)}(x,\tau)) \preceq (\overline{u}^{(1)}(x,0^{+}), \overline{v}^{(1)}(x,0^{+}))
\end{equation}
for all $x\in[r,s]$. Additionally, \eqref{3-2-5} and the first two equations of \eqref{Fixed} and \eqref{3-2-3} yield that
\begin{eqnarray}\label{3-2-9}
\left\{
\begin{array}{ll}
\zeta_{t}\geq d_{1}\int\limits^{s}_{r}J_{1}(x-y)\zeta(y,t)dy-d_{1}\zeta-a_{11}\zeta,\; &\, (x, t)\in\Omega^{0}, \\[2mm]
\eta _{t}\geq d_{2}\int\limits^{s}_{r}J_{1}(x-y)\eta(y,t)dy-d_{2}\eta-a_{22}\eta,\; &\, (x, t)\in\Omega^{0}, \\[2mm]
\end{array} \right.
\end{eqnarray}
where $(\zeta(x,t), \eta(x,t))= \big(\overline{u}^{(1)}(x,t)-u(x,t+\tau), \overline{v}^{(1)}(x,t)-v(x,t+\tau)\big)$. Based on \eqref{3-2-7}-\eqref{3-2-9}, applying \autoref{lemma 3-2-1} to $(\overline{u}^{(1)}(x,t)-u(x,t+\tau), \overline{v}^{(1)}(x,t)-v(x,t+\tau))$ yields that
\begin{equation*}
u(x,t+\tau)\leq \overline{u}^{(1)}(x,t)\text{~and~}v(x,t+\tau) \leq\overline{v}^{(1)}(x,t)
\end{equation*}
for all $(x,t)\in[r,s]\times[0,\tau]$. By the method of induction, it follows that
\begin{equation*}
u(x,t+\tau)\leq \overline{u}^{(1)}(x,t)\text{~and~}v(x,t+\tau) \leq \overline{v}^{(1)}(x,t)
\end{equation*}
for all $(x,t)\in[r,s]\times[0,+\infty).$ Using mathematical induction again we arrive at
\begin{equation}\label{3-2-10}
u(x,t+m\tau)\leq \overline{u}^{(m)}(x,t)\text{~and~}v(x,t+m\tau) \leq \overline{v}^{(m)}(x,t)
\end{equation}
for all $(x,t)\in[r,s]\times[0,\tau]$, where $m$ is an arbitrary natural number. Recalling that $(0,0)$ is the unique solution to problem \eqref{3-2-1} when $\lambda_{1}=0$, it follows from \eqref{3-2-10} and the nonnegativity of the solution to problem \eqref{Fixed} that
\begin{equation*}
\lim\limits_{m\rightarrow+\infty}\big(u(x,t+m\tau),v(x,t+m\tau)\big)=(0,0) \text{~for~all~}(x,t)\in[r,s]\times[0,+\infty).
\end{equation*}
which implies that
\begin{equation*}
\lim\limits_{t\rightarrow+\infty}(u(x,t),v(x,t))=(0,0) \text{~uniformly~for~}x\in[r,s].
\end{equation*}

(2) We first prove the first assertion. To begin with, we construct the following functions
\begin{equation*}
(\overline{U}(x,t), \overline{V}(x,t))=(C_{1}, C_{2}) \text{~for~all~}(x,t)\in[r,s]\times[0,\tau],
\end{equation*}
where $C_{1}$ and $C_{2}$ are defined in \eqref{C1C2}. For any $\lambda_{1}<0$, there exist $\theta>0$ such that $\lambda_{1}+\theta<0$. On this basis, define
\begin{eqnarray*}
\underline{U}(x,t)=
\left\{
\begin{array}{ll}
\epsilon \zeta_{1}(x,t),\; &\,x\in[r,s], t=0, \\[2mm]
\epsilon e^{(\lambda_{1}+\theta)\tau} \zeta_{1}(x,t),\; &\, x\in[r,s],t=0^{+}, \\[2mm]
\epsilon e^{(\lambda_{1}+\theta)(\tau-t)}\zeta_{1}(x,t),\; &\, x\in[r,s],t\in(0^{+}, \tau],
\end{array} \right.
\end{eqnarray*}
and
\begin{eqnarray*}
\underline{V}(x,t)=
\left\{
\begin{array}{ll}
\epsilon \eta_{1}(x,t),\; &\, x\in[r,s],t=0, \\[2mm]
\epsilon e^{(\lambda_{1}+\theta)\tau} \eta_{1}(x,t),\; &\, x\in[r,s],t=0^{+}, \\[2mm]
\epsilon e^{(\lambda_{1}+\theta)(\tau-t)}\eta_{1}(x,t),\; &\, x\in[r,s],t\in(0^{+}, \tau],
\end{array} \right.
\end{eqnarray*}
where $(\zeta_{1}, \eta_{1})$ is the unique eigenfunction pair corresponding to the principal eigenvalue $\lambda_{1}$, and $\epsilon$ is a positive number such that
\begin{equation*}
e^{(\lambda_{1}+\theta)\tau}H'(0)\epsilon\zeta_{1}(x,0)-H(\epsilon\zeta_{1}(x,0))\leq 0
\end{equation*}
and
\begin{equation*}
\frac{G'(0)\epsilon e^{(\lambda_{1}+\theta)(\tau-t)} \zeta_{1}(x,t)-G(\epsilon e^{(\lambda_{1}+\theta)(\tau-t)} \zeta_{1}(x,t))}{\epsilon}\leq \theta \min\limits_{\Omega^{0}}\eta_{1}(x,t).
\end{equation*}
A standard calculation shows that the constructed $(\overline{U}, \overline{V})$ and $(\underline{U}, \underline{V})$ are the ordered upper and lower solutions of the steady state problem \eqref{3-2-1}, respectively.

The two sequences of iteration $\big\{(\overline{U}^{(j)}, \overline{V}^{(j)})\big\}_{j=0}^{+\infty}$ and $\big\{(\underline{U}^{(j)}, \underline{V}^{(j)})\big\}_{j=0}^{+\infty}$ may then be derived from the iteration law \eqref{3-2-4} by choosing $(\overline{U}^{(0)}, \overline{V}^{(0)})=(\overline{U}, \overline{V})$ and $(\underline{U}^{(0)}, \underline{V}^{(0)})=(\underline{U}, \underline{V})$, respectively. Using a similar discussion as in the case 2 in (1), it follows that
\begin{equation*}
\underline{U}\leq \underline{U}^{(m-1)}\leq \underline{U}^{(m)}\leq \overline{U}^{(m)} \leq\overline{U}^{(m-1)}\leq \overline{U}\text{~and~}\underline{V}\leq \underline{V}^{(m-1)}\leq \underline{V}^{(m)}\leq \overline{V}^{(m)} \leq\overline{V}^{(m-1)}\leq \overline{V},
\end{equation*}
where $m\in\mathbb{Z}^{+}$. By the monotone boundedness theorem, it follows that
\begin{equation*}
\lim\limits_{m\rightarrow+\infty}(\underline{U}^{(m)}, \underline{V}^{(m)})=(U_{1}, V_{1})\text{~and~}\lim\limits_{m\rightarrow+\infty}(\overline{U}^{(m)}, \overline{V}^{(m)})=(U_{2}, V_{2}),
\end{equation*}
where $(U_{1}, V_{1})$ and $(U_{2}, V_{2})$ are the solutions of the steady state problem \eqref{3-2-1}. By  standard methods  (see, e.g., the proof of \cite[Theorem 4.2]{Zhou-Lin-Santos}), it follows that $(U_{1}, V_{1})$ and $(U_{2}, V_{2})$ are the maximum and minimum solutions of the steady state problem \eqref{3-2-1}, respectively.

It remains to show that the positive periodic solution of problem \eqref{3-2-1} is unique. Arguing indirectly, suppose that $(\tilde{U}, \tilde{V})$ and $(\ddot{U}, \ddot{V})$ are different solutions of problem \eqref{3-2-1}. Since the two different solutions are both strongly positive,
\begin{equation*}
\omega_{0}:=\min\Big\{\omega\geq 1 \big| \omega(\tilde{U}, \tilde{V})\succeq(\ddot{U}, \ddot{V})\text{~for~}(x,t)\in[r,s]\times[0,\tau] \Big\}
\end{equation*}
is well-defined and finite. Additionally, we have that $\omega_{0}(\tilde{U}, \tilde{V})\succeq(\ddot{U}, \ddot{V})$ for $(x,t)\in[r,s]\times[0,\tau]$ and there exists some $( x_{0},t_{0})\in[r,s]\times(0^{+},\tau]$ such that
\begin{equation*}
\mbox{(i)}~\omega_{0}\tilde{U}(x_{0},t_{0}) =\ddot{U}(x_{0},t_{0}) \text{~or~} \mbox{(ii)}~\omega_{0}\tilde{V}(x_{0},t_{0}) =\ddot{V}(x_{0},t_{0}).
\end{equation*}
To simplify notation, let
\begin{equation*}
P(x,t)=\omega_{0}\tilde{U}(x,t)-\ddot{U}(x,t) \text{~and~} \mbox{(ii)}~Q(x,t)=\omega_{0}\tilde{V}(x,t)-\ddot{V}(x,t) .
\end{equation*}
If $\mbox{(i)}$ holds, then
\begin{equation*}
\begin{aligned}
0\geq P_{t}( x_{0},t_{0})&=d_{1}\int\limits^{s}_{r}J_{1}(x_{0}-y)P(y,t_{0})dy-(d_{1}+a_{11})P( x_{0},t_{0})+a_{12}Q( x_{0},t_{0})\\
&\geq a_{12}Q( x_{0},t_{0}).
\end{aligned}
\end{equation*}
This implies that the above inequality can hold only if $Q( x_{0},t_{0})=0$, that is, $\mbox{(ii)}$ holds. We claim $\omega_{0}=1$. If not, $\omega_{0}>1$. From the condition $\omega_{0}>1$ and assumption \textbf{(G)}, it follows that $\omega_{0}G(\tilde{U})>G(\omega_{0}\tilde{U})$. If $\mbox{(ii)}$ holds, then
\begin{equation*}
\begin{aligned}
0\geq Q_{t}( x_{0},t_{0})&=d_{1}\int\limits^{s}_{r}J_{1}(x_{0}-y)Q(y,t_{0})dy-(d_{2}+a_{22})Q( x_{0},t_{0})+\omega_{0}G(\tilde{U}( x_{0},t_{0}))-G(\ddot{U}( x_{0},t_{0}))\\
&\geq \omega_{0}G(\tilde{U}( x_{0},t_{0}))-G(\ddot{U}( x_{0},t_{0}))> G( \omega_{0}\tilde{U}( x_{0},t_{0}))-G(\ddot{U}( x_{0},t_{0}))=0,
\end{aligned}
\end{equation*}
which is impossible. Hence, $\omega_{0}=1$, which means that $(\tilde{U}, \tilde{V})\succeq(\ddot{U}, \ddot{V})$ for $(x,t)\in[r,s]\times [0,\tau]$. By exchanging $(\tilde{U}, \tilde{V})$ and $(\ddot{U}, \ddot{V})$, the conclusion that $(\ddot{U}, \ddot{V})\succeq(\tilde{U}, \tilde{V})$ for $(x,t)\in[r,s]\times [0,\tau]$ may also be obtained. As a result, $(\ddot{U}, \ddot{V})=(\tilde{U}, \tilde{V})$ for $x\in[r,s]$ and $t\in[0,\tau]$. This completes the proof of the first assertion.

Next, we prove the second assertion. A standard calculation yields that $(\overline{U}, \overline{V})$ is a supersolution of \eqref{Fixed}. Recall that $(u_{0}(x), v_{0}(x))\succeq, \not\equiv\mathbf{0}$ for $x\in [r,s]$. We assume, without loss of generality, that $(u_{0}(x), v_{0}(x))\succ \mathbf{0}$ for $x\in [r,s]$. Otherwise, \autoref{lemma 3-2-1} yields that $(u(x,\tau), v(x,\tau))\succ \mathbf{0}$ for $x\in [r,s]$, and it can be seen as a new initial function. By decreasing the $\epsilon$ if necessary, it follows from a standard calculation that $(\underline{U}, \underline{V})$ is a subsolution of \eqref{Fixed}. Then. \autoref{lemma 3-2-1} yields that the solution $(u, v)$ of \eqref{Fixed} satisfies the following estimate
\begin{equation*}
\underline{U}\leq u \leq \overline{U}\text{~and~} \underline{V}\leq v\leq \overline{V}
\end{equation*}
for $(x,t)\in[r,s]\times[0, +\infty)$. By a similar discussion as in the case 2 in (1), we have that
\begin{equation}\label{3-2-11}
\underline{U}^{(m)}(x,t)\leq u(x,t+m\tau)\leq \overline{U}^{(m)}(x,t)\text{~and~}\underline{V}^{(m)}(x,t)\leq v(x,t+m\tau)\leq \overline{V}^{(m)}(x,t)
\end{equation}
for $(x,t)\in\Omega^{k}$, where $m=0, 1, 2, \cdots$. With the help of the first assertion, we have that
\begin{equation}\label{3-2-12}
\lim\limits_{m\rightarrow+\infty}(\underline{U}^{(m)}, \underline{V}^{(m)})=\lim\limits_{m\rightarrow+\infty}(\overline{U}^{(m)}, \overline{V}^{(m)})=(U, V) \text{~for~}(x,t)\in[r,s]\times[0, \tau].
\end{equation}
Noting that $(U, V)$ is periodic, \eqref{3-2-12} also holds for $(x,t)\in[r,s]\times[0, +\infty)$. Hence, it follows from \eqref{3-2-11} that
\begin{equation*}
\lim\limits_{m\rightarrow+\infty}\big(u(x,t+m\tau),v(x,t+m\tau)\big)=\big(U(x,t), V(x,t)\big) \text{~uniformly~for~}(x,t)\in[r,s]\times[0, +\infty).
\end{equation*}
This proof is finished.
\end{proof}
\end{theorem}
\autoref{theorem 3-2-1} above gives the asymptotical profiles of problem \eqref{Fixed} when the infection environment $[r,s]$ is fixed. Next, we will analyse the situation when the infection environment is increasing and eventually goes to infinity.
\begin{lemma}\label{lemma 3-2-2}
Assume that $[r,s]=[-l,l]$. Let $\lambda_{1}(\infty)$ denote the eigenvalue of \eqref{3-1-1} as $l$ goes to infinity. If $\lambda_{1}(\infty)<0$, then
\begin{equation*}
\lim\limits_{l\rightarrow\infty}(U_{l}(x,t),V_{l}(x,t))=(U_{1}(t),V_{1}(t)) \text{~locally~uniformly~in~}\mathds{R},
\end{equation*}
where $(U_{l},V_{l})$ denotes the solution to problem \eqref{3-2-1} with $[r, s]$ replaced by $[-l,l]$, and $(U_{1},V_{1})$ stands for the unique solution to problem
\begin{eqnarray}\label{3-24}
\left\{
\begin{array}{ll}
\zeta_{t}=a_{12}\eta-a_{11}\zeta,~~~~~~~~t\in(0^{+}, \tau], \\[2mm]
\eta_{t}=G(\zeta)-a_{22}\eta, ~~~~~~~~t\in(0^{+}, \tau], \\[2mm]
\big(\zeta(0),\eta(0) \big)=\big(\zeta(\tau), \eta(\tau)\big),\\[2mm]
\big(\zeta(0^{+}), \eta(0^{+})\big)=\big(H(\zeta(0)),\eta(0)\big).
\end{array} \right.
\end{eqnarray}
\begin{proof}
It is well known that the asymptotic profiles of problem \eqref{3-24} are determined by the eigenvalue $\mu_{1}$ of the following eigenvalue problem
\begin{eqnarray*}
\left\{
\begin{array}{ll}
\phi_{t}=a_{12}\varphi-a_{11}\phi+\mu\phi,~~~~~~~~~~~~t\in(0^{+}, \tau], \\[2mm]
\varphi_{t}=G'(0)\phi-a_{22}\varphi+\mu\varphi, ~~~~~~~~t\in(0^{+}, \tau], \\[2mm]
\big(\phi(0), \varphi(0)\big)=\big(\phi(\tau),\varphi(\tau)\big),\\[2mm]
\big(\phi(0^{+}),\varphi(0^{+})\big)=\big(H'(0)\phi(0),\varphi(0)\big).
\end{array} \right.
\end{eqnarray*}
Now, we claim that if $\lambda_{1}(\infty)<0$, then $\mu_{1}<0$. Otherwise, assume that $\mu_{1}\geq0$ and its corresponding nonnegative eigenfunction is $(\phi_{1}(t), \varphi_{1}(t))$. Next, it follows to check that for any given $l>0$, $(\zeta, \eta, \underline{\lambda})=(\phi_{1}, \varphi_{1}, 0)$ satisfies \eqref{3-1-19}. By standard computation, for $x\in[-l,l]$,
\begin{equation*}
d_{1}\int_{-l}^{l}J_{1}(x-y)\zeta(y,t)dy-d_{1}\zeta-a_{11}\zeta+a_{12}\eta-\zeta_{t}\leq, \not\equiv  -a_{11}\phi_{1}+a_{12}\varphi_{1}-\phi_{1t}\leq 0,
\end{equation*}
and
\begin{equation*}
d_{1}\int_{-l}^{l}J_{2}(x-y)\eta(y,t)dy-d_{2}\eta-a_{22}\eta+G'(0)\zeta-\eta_{t}\leq, \not\equiv  -a_{22}\varphi_{1}+G'(0)\phi_{1}-\varphi_{1t}\leq 0.
\end{equation*}
Then, \autoref{lemma 3-1-3} yields that $\lambda_{1}(l)>0$ for any $l>0$, which implies that $\lambda_{1}(\infty)\geq0$. This contradicts with the condition $\lambda_{1}(\infty)<0$. Therefore, the claim holds  by contraposition.

By using a similar proof method as in  the first assertion of \autoref{theorem 3-2-1}\textcolor{blue}{(2)}, it follows that when $\mu_{1}<0$, \eqref{3-24} has a unique positive solution $(U_{1},V_{1})$. Finally, the desired conclusion can be obtained by standard methods  (see, e.g., the proof of \cite[Lemma 3.11]{Wang-Du-1}). This ends the proof.
\end{proof}
\end{lemma}
\section{\bf The asymptotic profile in a moving environment}\label{Section-4}
This section first presents some basic knowledge. Based on this, the vanishing-spreading dichotomy is proved and sufficient conditions used for determining spreading or vanishing are provided. Throughout this section, we always denote the solution of \eqref{Zhou-Lin} by $(u, v, g, h)$ or $(u(x,t), v(x,t),\\ g(t), h(t))$.
\subsection{Some introductory results}
In here, the definition of the super- and sub-solutions of model \eqref{Zhou-Lin}, the comparison principle applicable to model \eqref{Zhou-Lin}, and the definition of the spreading and vanishing of model \eqref{Zhou-Lin} will be presented in turn.
\begin{definition}\label{definition 4-1}
If for any $T>0$, $\overline{g}, \overline{h}\in\mathbb{C}[0, T]\cap\big[\mathbb{C}^{1}(n\tau, (n+1)\tau]\cap \mathbb{C}^{1}(N\tau, T]\big]$, $\overline{u}, \overline{v}, \overline{u}_{t}, \overline{u}_{t}\in\mathbb{C}\big(\overline{\Omega}_{n\tau}^{(n+1)\tau}\big)\cap\mathbb{C}\big(\overline{\Omega}_{N\tau}^{T}\big)$ satisfy $\overline{u}, \overline{v}\geq 0$ and
\begin{eqnarray}\label{4-1-1}
\left\{
\begin{array}{ll}
\overline{u}_{t}\geq d_{1}\int\limits^{\overline{h}(t)}_{\overline{g}(t)}J_{1}(x-y)\overline{u}(y,t)dy-d_{1}\overline{u}-a_{11}\overline{u}+a_{12}\overline{v},\; &\, (x,t)\in\Omega_{n\tau}^{(n+1)\tau}\cup\Omega_{N\tau}^{T}, \\[2mm]
\overline{v}_{t}\geq d_{2}\int\limits^{\overline{h}(t)}_{\overline{g}(t)}J_{2}(x-y)\overline{v}(y,t)dy-d_{2}\overline{v}-a_{22}\overline{v}+G(\overline{u}),\; &\, (x,t)\in\Omega_{n\tau}^{(n+1)\tau}\cup\Omega_{N\tau}^{T}, \\[2mm]
\overline{g}'(t)\leq-\mu_{1}\int\limits^{\overline{h}(t)}_{\overline{g}(t)}\int\limits^{\overline{g}(t)}_{-\infty}J_{1}(x-y)\overline{u}(x,t)dydx\\
~~~~~~~~~~~~~~~~~~~~~~~~~~~-\mu_{2}\int\limits^{\overline{h}(t)}_{\overline{g}(t)}\int\limits^{\overline{g}(t)}_{-\infty}J_{2}(x-y)\overline{v}(x,t)dydx,\; &\, t\in(0, T],\\[2mm]
\overline{h}'(t)\geq\mu_{1}\int\limits^{\overline{h}(t)}_{\overline{g}(t)}\int\limits^{+\infty}_{\overline{h}(t)}J_{1}(x-y)\overline{u}(x,t)dydx\\
~~~~~~~~~~~~~~~~~~~~~~~~~~~+\mu_{2}\int\limits^{\overline{h}(t)}_{\overline{g}(t)}\int\limits^{+\infty}_{\overline{h}(t)}J_{2}(x-y)\overline{v}(x,t)dydx,\; &\, t\in(0, T],\\[2mm]
\overline{u}(\overline{g}(t),t), \overline{u}(\overline{h}(t),t)\geq0,~~~~~~~~ \overline{v}(\overline{g}(t),t), \overline{v}(\overline{h}(t),t)\geq 0,\; &\, t\in(0, T],\\[2mm]
\overline{u}(x,0)\geq u_{0}(x), \overline{v}(x,0)\geq v_{0}(x),~~ -\overline{g}(0), \overline{h}(0)\geq h_{0}, \; &\,x\in[-h_{0},h_{0}], \\[2mm]
\overline{u}(x,(n\tau)^{+})\geq H(\overline{u}(x,n\tau)), ~~~~~\overline{v}(x,(n\tau)^{+})\geq\overline{v}(x,n\tau),\; &\, x\in(\overline{g}(n\tau), \overline{h}(n\tau)),
\end{array} \right.
\end{eqnarray}
then the quadruple $(\overline{u}, \overline{v}, \overline{g}, \overline{h})$ is called a supersolution of \eqref{Zhou-Lin}.
\end{definition}
\begin{remark}\label{remark 4-1-1}
By reversing all the inequalities in \eqref{4-1-1}, a subsolution of model \eqref{Zhou-Lin} can also be defined. In this section,
denote the super- and sub-solutions of model \eqref{Zhou-Lin} by the quadruples $(\overline{u}, \overline{v}, \overline{g}, \overline{h})$ and $(\underline{u}, \underline{v}, \underline{g}, \underline{h})$, respectively.
\end{remark}
\begin{lemma}\label{lemma 4-1-1}
Denote the unique solution of model \eqref{Zhou-Lin} by $(u,v,g,h)$. If the quadruples $(\overline{u}, \overline{v}, \overline{g}, \overline{h})$ and $(\underline{u}, \underline{v}, \underline{g}, \underline{h})$ are the ordered super- and sub-solutions of model \eqref{Zhou-Lin}, respectively, then
\begin{eqnarray*}
\left\{
\begin{array}{l}
\overline{g}(t)\leq g(t)\leq\underline{g}(t), ~~~~~~~~~~~ \underline{h}(t)\leq h(t)\leq \overline{h}(t)~~~~~~~~~~\text{~for~}t\in\mathds{R}^{+}, \\[2mm]
\underline{u}(x,t)\leq u(x,t)\leq\overline{u}(x,t), \underline{v}(x,t)\leq v(x,t)\leq\overline{v}(x,t) \text{~for~}x\in[g(t), h(t)], t\in\mathds{R}^{+}.
\end{array} \right.
\end{eqnarray*}
\begin{proof}
Here, we only need to prove that the conclusion holds for the supersolution,  due to the fact that a similar analysis yields that the conclusion also holds for the subsolution.

By using the assumption ($\mathbf{H}$), at the pulse moment $t=0^{+}$, we have that
\begin{eqnarray*}
\left\{
\begin{array}{l}
u(x,0^{+})\in\mathbb{C}[-h_{0},h_{0}], ~u(\pm h_{0},0^{+})=0~\text{and}~u(x,0^{+})>0~\text{in}~(-h_{0},h_{0}), \\[2mm]
v(x,0^{+})\in\mathbb{C}[-h_{0},h_{0}], ~v(\pm h_{0},0^{+})=0~\text{and}~v(x,0^{+})>0~~\text{in}~(-h_{0},h_{0}).
\end{array} \right.
\end{eqnarray*}
By viewing $(u(x,0^{+}), v(x,0^{+}))$ as an initial function, the standard discussion (see, e.g., the proof of \cite[Theorem 3.1]{Cao-Du-Li-Li}) yields that the conclusion holds for $t\in(0^{+},\tau]$, that is, the conclusion holds for $t\in(0,\tau]$. From mathematical induction, it follows that the conclusion holds for $t\in(0,\infty)$. The proof is completed.
\end{proof}
\end{lemma}

\autoref{theorem 2-1} and \autoref{lemma 2-0} show that $(u(x,t), v(x,t))\succ \mathbf{0}$ for $x\in(g(t), h(t))$ and $t\in[0, \infty)$.
This implies that $h'(t), -g'(t)> 0$ for $t\in[0, \infty)$. Therefore,
\begin{equation*}
\lim\limits_{t\rightarrow\infty} g(t)=g_{\infty}\in [-\infty, -h_{0})\text{~and~}\lim\limits_{t\rightarrow\infty} h(t)=h_{\infty}\in (h_{0}, +\infty]
\end{equation*}
are always well-defined.
\begin{definition}\label{definition 4-2}
The model \eqref{Zhou-Lin} is called as vanishing if
\begin{equation*}
h_{\infty}-g_{\infty}<\infty~\text{and}~\lim\limits_{t\rightarrow\infty}\|u(x,t)+v(x,t)\|_{\mathbb{C}(\mathds{R})}=0,
\end{equation*}
and spreading if
\begin{equation*}
h_{\infty}-g_{\infty}=\infty~\text{and}~\lim\limits_{t\rightarrow\infty}(u(x,t), v(x,t))=(U_{1}(t),V_{1}(t)) \text{~locally~uniformly~in~}\mathds{R},
\end{equation*}
where $(U_{1},V_{1})$ denotes the solution of \eqref{3-24}.
\end{definition}
\subsection{Spreading-vanishing dichotomy}
This subsection will prove that when the time $t$ tends to infinity, model \eqref{Zhou-Lin} is either vanishing or spreading.
This conclusion will be discussed separately for the three cases: $\lambda_{1}(-\infty, +\infty)>0$, $=0$, and $<0$.
\begin{lemma}\label{lemma 4-2-2}
Let $\lambda_{1}(\infty)=\lambda_{1}(-\infty, +\infty)$. If $\lambda_{1}(\infty)> 0$, then a constant $k_{0}>0$,
which depends only the eigenvalue $\lambda_{1}(\infty)$, can be found such that
\begin{equation*}
\lim\limits_{t\rightarrow\infty}e^{k_{0}t}\|u(x,t)+v(x,t)\|_{\mathbb{C}(\mathds{R})}=0.
\end{equation*}
Moreover, we have that $(g_{\infty}, h_{\infty})$ is a finite interval.
\begin{proof}
The first conclusion is proved first. Define, for $x\in[g(t), h(t)]$ and $t\in [0,+\infty)$,
\begin{eqnarray}\label{3-4}
\begin{array}{l}
\overline{u}(x,t):=Me^{-\sigma t}\zeta(x,t), ~~~~\overline{v}(x,t):=Me^{-\sigma t}\eta(x,t),
\end{array}
\end{eqnarray}
where $M$ and $\sigma$ are positive constants to be determined later, and $(\zeta(x,t), \eta(x,t))$ is the eigenfunction of \eqref{3-1-1},
in which $[r,s]$ is substituted by $(-\infty, +\infty)$. Recalling that $(\zeta(x,t), \eta(x,t))\in\mathcal{X}^{++}_{\tau}$, $M$ can be determined such that $M(\zeta(x,0), \eta(x,0))\succeq(u_{0}(x), v_{0}(x))$ for $x\in [-h_{0}, h_{0}]$. Then, check the impulsive conditions. Clearly, assumption ($\mathbf{H}$) gives that
\begin{eqnarray*}
\begin{array}{ll}
\overline{u}(x,(k\tau)^{+})=Me^{-\sigma k\tau}\zeta(x,(k\tau)^{+})=Me^{-\sigma k\tau}H'(0)\zeta(x,k\tau)>H(\overline{u}(x,k\tau)),\; &\, x\in\mathds{R}, \\[2mm]
\overline{v}(x,(k\tau)^{+})=Me^{-\sigma k\tau}\eta(x,(k\tau)^{+})=Me^{-\sigma k\tau}\zeta(x,k\tau)=\overline{v}(x,k\tau),\; &\, x\in\mathds{R}. \\[2mm]
\end{array}
\end{eqnarray*}
Finally, assumption ($\mathbf{G}$) yields that
\begin{equation*}
\begin{aligned}
&\overline{u}_{t}-d_{1}\int^{h(t)}_{g(t)}J_{1}(x-y)\overline{u}(y,t)dy+d_{1}\overline{u}+a_{11}\overline{u}-a_{12}\overline{v}\\
\geq&~Me^{-\sigma t}\big[-\sigma\zeta+\zeta_{t}-d_{1}\int^{+\infty}_{-\infty}J_{1}(x-y)[\zeta(t,y)-\zeta(t,x)]dy+a_{11}\zeta-a_{12}\eta\big] \\
\geq&~Me^{-\sigma t}\big[\lambda_{1}(\infty)\zeta-\sigma\zeta\big]\geq 0,
\end{aligned}
\end{equation*}
and
\begin{equation*}
\begin{aligned}
&\overline{v}_{t}-d_{2}\int^{h(t)}_{g(t)}J_{2}(x-y)\overline{v}(y,t)dy+d_{2}\overline{v}+a_{22}\overline{v}-G(\overline{u})\\
\geq&~Me^{-\sigma t}\big[-\sigma\eta+\eta_{t}-d_{2}\int^{+\infty}_{-\infty}J_{2}(x-y)[\eta(y,t)-\eta(x,t)]dy+a_{22}\eta-G'(0)\zeta\big] \\
\geq&~Me^{-\sigma t}\big[\lambda_{1}(\infty)\eta-\sigma\eta\big]\geq 0,
\end{aligned}
\end{equation*}
for $x\in(-\Gamma(t), \Gamma(t))$ and $t\in((k\tau)^{+}, (k+1)\tau]$ provided that $\sigma :=\lambda_{1}(\infty)/2$. Now, applying \autoref{lemma 3-2-1} to $(\overline{u}-u, \overline{v}-v)$ yields that
\begin{equation}\label{4-2-2}
\begin{aligned}
u\leq \overline{u}\text{~and~}v\leq \overline{v} \text{~for~} (t,x)\in\mathds{R}^{+}\times\mathds{R}.
\end{aligned}
\end{equation}
Therefore, taking $k_{0}:=\sigma/2$, it follows from \eqref{4-2-2} that
\begin{equation*}
\lim\limits_{t\rightarrow\infty}e^{k_{0}t}\|u(x,t)+v(x,t)\|_{\mathbb{C}(\mathds{R})}=0.
\end{equation*}
This proves the first assertion.

Next, we prove the second assertion. The first assertion tells us that a constant $T_{0}$ can be found such that for any $t>T_{0}$, it follows that
\begin{equation*}
e^{k_{0}t}u(x,t)\leq \frac{1}{4(\mu_{1}+\mu_{2})}\text{~and~}e^{k_{0}t}v(x,t)\leq \frac{1}{4(\mu_{1}+\mu_{2})}
\end{equation*}
for $(x,t)\in[g(t), h(t)]\times(T_{0}, \infty)$. This combined with assumption ($\mathbf{J}$) yields that
\begin{equation*}
\begin{aligned}
h'(t)-g'(t)=&\mu_{1}\int\limits^{h(t)}_{g(t)}\int\limits^{+\infty}_{h(t)}J_{1}(x-y)u(x,t)dydx+\mu_{2}\int\limits^{h(t)}_{g(t)}\int\limits^{+\infty}_{h(t)}J_{2}(x-y)u(x,t)dydx\\
&+\mu_{1}\int\limits^{h(t)}_{g(t)}\int\limits^{g(t)}_{-\infty}J_{1}(x-y)u(x,t)dydx+\mu_{2}\int\limits^{h(t)}_{g(t)}\int\limits^{g(t)}_{-\infty}J_{2}(x-y)v(x,t)dydx\\
\leq&~2(\mu_{1}+\mu_{2})\int\limits^{h(t)}_{g(t)}u(x,t)dx+2(\mu_{1}+\mu_{2})\int\limits^{h(t)}_{g(t)}v(x,t)dx\\
\leq&~[h(t)-g(t)]e^{-k_{0}t}
\end{aligned}
\end{equation*}
for $t> T_{0}$. Then, a standard calculation gives us
\begin{equation*}
h(t)-g(t)\leq [h(T)-g(T)]e^{\frac{e^{-k_{0}T}}{k_{0}}-\frac{e^{-k_{0}t}}{{k_{0}}}},
\end{equation*}
which implies that
\begin{equation*}
h_{\infty}-g_{\infty}\leq [h(T)-g(T)]e^{\frac{e^{-k_{0}T}}{k_{0}}}<\infty.
\end{equation*}
Therefore, the second assertion is valid. The proof is complected.
\end{proof}
\end{lemma}
\autoref{lemma 4-2-2} shows that model \eqref{Zhou-Lin} is vanishing when $\lambda_{1}(\infty)> 0$. Subsequently, we will discuss the case $\lambda_{1}(\infty)= 0$.
\begin{lemma}\label{lemma 4-2-3}
Assume that the region at infinity is particularly unfavourable for the survival of the infectious agents, that is, $h_{\infty}-g_{\infty}< \infty$. If $\lambda_{1}(\infty)= 0$,
then we have that
\begin{equation*}
\lim\limits_{t\rightarrow\infty}\|u(x,t)+v(x,t)\|_{\mathbb{C}[g(t), h(t)]}=0.
\end{equation*}
\begin{proof}
Let $(\tilde{u}, \tilde{v})$ be the unique solution of the problem
\begin{eqnarray*}
\left\{
\begin{array}{ll}
\tilde{u}_{t}=-a_{11}\tilde{u}+a_{12}\tilde{v}, ~~~~~~~~~~~~~~~~~~~~~~~~t\in((k\tau)^{+}, (k+1)\tau], \\[2mm]
\tilde{v}_{t}=-a_{22}\tilde{v}+G(\tilde{u}),  ~~~~~~~~~~~~~~~~~~~~~~~~t\in((k\tau)^{+}, (k+1)\tau], \\[2mm]
\tilde{u}(0)=C_{1}, ~~~~~~~~~~~~~~~~~~\tilde{v}(0)=C_{2},  \\[2mm]
\tilde{u}((k\tau)^{+})=H(\tilde{u}(k\tau)), ~\tilde{v}((k\tau)^{+})=\tilde{v}(k\tau),
\end{array} \right.
\end{eqnarray*}
where $C_{1}$ and $C_{2}$ are defined in \eqref{C1C2}. Then, applying \autoref{lemma 3-2-1} to $(\tilde{u}-u, \tilde{v}-v)$ yields that
\begin{equation}\label{4-2-3}
\begin{aligned}
u\leq \tilde{u}\text{~and~}v\leq  \tilde{v}\text{~for~} (t,x)\in\mathds{R}^{+}\times\mathds{R}.
\end{aligned}
\end{equation}
From the condition $\lambda_{1}(\infty)= 0$, it follows from \autoref{lemma 3-1-5} that $\mu_{1}=0$, where $\mu_{1}$ is the principal eigenvalue of \eqref{ODE IM}.
Since $\mu_{1}=0$, $(0,0)$ is the only solution of problem \eqref{3-24}. Finally, the upper and lower solutions method show that
\begin{equation*}
\lim\limits_{t\rightarrow\infty}\tilde{u}(t)=\lim\limits_{t\rightarrow\infty}\tilde{v}(t)=0.
\end{equation*}
This together with \eqref{4-2-3} yields that
\begin{equation*}
\lim\limits_{t\rightarrow\infty}\|u(x,t)+v(x,t)\|_{\mathbb{C}(\mathds{R})}=0.
\end{equation*}
This ends the proof.
\end{proof}
\end{lemma}
Next, we will discuss the dynamic behavior of model \eqref{Zhou-Lin} when $\lambda_{1}(\infty)<0$. To do this, we start by proving the following two lemmas.

\begin{lemma}\label{lemma 4-2-1}
If $(g_{\infty}, h_{\infty})$ is a finite interval, then
\begin{equation*}
\lim\limits_{t\rightarrow\infty}u(x,t)=\lim\limits_{t\rightarrow\infty}v(x,t)=0 \text{~uniformly~for~}x\in[g(t), h(t)]
\end{equation*}
and $\lambda_{1}(g_{\infty}, h_{\infty})\geq 0$, where $\lambda_{1}(g_{\infty}, h_{\infty})$ is the eigenvalue of problem \eqref{3-1-1} in which $[r,s]$ is substituted by $[g_{\infty}, h_{\infty}]$.
\begin{proof}
Begin by showing the second assertion. Arguing indirectly, suppose that $\lambda_{1}(g_{\infty}, h_{\infty})<0$. \autoref{lemma 3-1-4}\textcolor{blue}{(2)}
gives that a sufficiently small positive constant $\epsilon$ can be found such that
\begin{equation*}
\lambda_{1}(g_{\infty}+\epsilon, h_{\infty}-\epsilon)<0.
\end{equation*}
The assumption that $f\in \mathbb{C}(\mathbb{R})$ and $f(0)>0$ in \textbf{(G)} yields that there exist sufficiently small positive
constants $\epsilon_{0}, \delta_{0}\leq \frac{h_{0}}{2}$ such that
\begin{eqnarray*}
J_{i}(x)>\delta_{0}\text{~for~}x\in(-3\epsilon_{0}, 3\epsilon_{0}), i=1,2.
\end{eqnarray*}
Then, for given $\epsilon_{1}<\min\{h_{0}, 2\epsilon_{0}\}$, a sufficiently large positive integer $k_{0}$ can be found such that
\begin{equation}\label{4-2-1}
(g_{\infty}, h_{\infty}-\epsilon_{1})\prec(g(t), h(t))\prec(g_{\infty}+\epsilon_{1}, h_{\infty})\text{~for~}t\geq k_{0}\tau.
\end{equation}
Now, we consider the following problem
\begin{eqnarray*}
\left\{
\begin{array}{ll}
\tilde{u}_{t}=d_{1}\int\limits^{h_{\infty}-\epsilon_{1}}_{g_{\infty}+\epsilon_{1}}J_{1}(x-y)\tilde{u}(y,t)dy-d_{1}\tilde{u}-a_{11}\tilde{u}+a_{12}\tilde{v},\; &\, x\in\Gamma, t\in((k\tau)^{+}, (k+1)\tau], \\[2mm]
\tilde{v}_{t}=d_{2}\int\limits^{h_{\infty}-\epsilon_{1}}_{g_{\infty}+\epsilon_{1}}J_{2}(x-y)\tilde{v}(y,t)dy-d_{2}\tilde{v}-a_{22}\tilde{v}+G(\tilde{u}),\; &\, x\in\Gamma, t\in((k\tau)^{+}, (k+1)\tau], \\[2mm]
\tilde{u}(x,(k\tau)^{+})=H(\tilde{u}(x,k\tau)), ~~~~\tilde{v}(x,(k\tau)^{+})=\tilde{v}(x,k\tau),\; &\, x\in\Gamma, \\[2mm]
\tilde{u}(x,k_{0}\tau)=u(x,k_{0}\tau), ~~~~~~~~~~~~\tilde{v}(x,k_{0}\tau)=v(x,k_{0}\tau), \; &\,x\in\Gamma, k=k_{0}, k_{0}+1, \cdots,
\end{array} \right.
\end{eqnarray*}
where $\Gamma=[g_{\infty}+\epsilon_{1}, h_{\infty}-\epsilon_{1}]$. It follows from \eqref{4-2-1} that
\begin{equation*}
\int\limits^{h_{\infty}-\epsilon_{1}}_{g_{\infty}+\epsilon_{1}}J_{1}(x-y)u(y,t)dy\leq \int\limits^{h(t)}_{g(t)}J_{1}(x-y)u(y,t)dy \text{~for~} x\in\Gamma , t\in[k_{0}\tau,\infty) .
\end{equation*}
Therefore, from \autoref{lemma 3-2-1} we have that
\begin{equation*}
\tilde{u}(x,t)\leq u(x,t)\text{~and~} \tilde{v}(x,t)\leq v(x,t)\text{~for~} x\in\Gamma \text{~and~}t\in[k_{0}\tau,\infty).
\end{equation*}
Noting that $\lambda_{1}(g_{\infty}+\epsilon_{1}, h_{\infty}-\epsilon_{1})<0$, it follows from \autoref{theorem 3-2-1}\textcolor{blue}{(2)} that
\begin{equation*}
\lim\limits_{m\rightarrow+\infty}\big(\tilde{u}(x,t+m\tau),\tilde{v}(x,t+m\tau)\big)=\big(U^{\epsilon_{1}}(x,t), V^{\epsilon_{1}}(x,t)\big) \text{~uniformly~for~}(x,t)\in\Gamma \times [0, \tau],
\end{equation*}
where $(U^{\epsilon_{1}}, V^{\epsilon_{1}})$ is the only nonnegative periodic solution of \eqref{3-2-1} in which $[r,s]$ is substituted by $[g_{\infty}+\epsilon_{1}, h_{\infty}-\epsilon_{1}]$. Then, a integer $k_{1}\geq k_{0}$ can be found such that for any $k\geq k_{1}$, we have that
\begin{equation*}
\frac{1}{4}\big(U^{\epsilon_{1}}(x,t), V^{\epsilon_{1}}(x,t)\big)\preceq(\tilde{u}(x,t+k\tau), \tilde{u}(x,t+k\tau))\preceq (u(x,t+k\tau), v(x,t+k\tau))
\end{equation*}
for $x\in\Gamma$ and $t\in[0, \tau]$. To simplify notation, let
\begin{eqnarray*}
\vartheta<\min\bigg\{\min\limits_{\Gamma\times[0, \tau]}U^{\epsilon_{1}}(x,t),\min\limits_{\Gamma\times[0, \tau]}V^{\epsilon_{1}}(x,t)\bigg\}.
\end{eqnarray*}
For $t\geq k_{1}\tau$, we have that
\begin{equation*}
\begin{split}
h'(t)&=\mu_{1}\int\limits^{h(t)}_{g(t)}\int\limits^{+\infty}_{h(t)}J_{1}(x-y)u(x,t)dydx+\mu_{2}\int\limits^{h(t)}_{g(t)}\int\limits^{+\infty}_{h(t)}J_{2}(x-y)v(x,t)dydx\\
&\geq\mu_{1}\int\limits^{h(t)}_{h(t)-2\epsilon_{0}}\int\limits^{h(t)+\epsilon_{0}}_{h(t)}J_{1}(x-y)u(x,t)dydx
+\mu_{2}\int\limits^{h(t)}_{h(t)-2\epsilon_{0}}\int\limits^{h(t)+\epsilon_{0}}_{h(t)}J_{2}(x-y)v(x,t)dydx\\
&\geq\mu_{1}\epsilon_{0}\delta_{1}\int\limits^{h(t)}_{h(t)-2\epsilon_{0}}u(x,t)dx+\mu_{2}\epsilon_{0}\delta_{1}\int\limits^{h(t)}_{h(t)-2\epsilon_{0}}v(x,t)dx\\
&\geq\mu_{1}\epsilon_{0}\delta_{1}\int\limits^{h_{\infty}-\epsilon_{1}}_{h_{\infty}-2\epsilon_{0}}u(x,t)dx
+\mu_{2}\epsilon_{0}\delta_{1}\int\limits^{h_{\infty}-\epsilon_{1}}_{h_{\infty}-2\epsilon_{0}}v(x,t)dx\\
&\geq (2\epsilon_{0}-\epsilon_{1})(\mu_{1}+\mu_{2})\epsilon_{0}\delta_{1}\vartheta>0,
\end{split}
\end{equation*}
which implies that $h_{\infty}=+\infty$. A similar process yields that $g_{\infty}=-\infty$. Therefore, we have that $(h_{\infty}, g_{\infty})=\mathds{R}$.
This contradicts the assumption that $(g_{\infty}, h_{\infty})$ is a finite interval. The proof of the second assertion is now completed.

The first assertion is shown next. Let $(\bar{u}(x,t), \bar{v}(x,t))$ be the solution of \eqref{Fixed} in which $[r,s]$ and $(\bar{u}(x,0), \bar{v}(x,0))$ are substituted by $[g_{\infty}, h_{\infty}]$ and $(\|u_{0}\|_{\infty},\|v_{0}\|_{\infty} )$, respectively. \autoref{lemma 3-2-1} yields that
\begin{equation*}
\mathbf{0}\preceq(u(x,t), v(x,t))\preceq (\bar{u}(x,t), \bar{v}(x,t)) \text{~for~} x\in[g_{\infty}, h_{\infty}], t\in[0,\infty).
\end{equation*}
With the help of \autoref{theorem 3-2-1}\textcolor{blue}{(1)}, it follows that
\begin{equation*}
\lim\limits_{t\rightarrow+\infty}\|\bar{u}(x,t)+\bar{v}(x,t)\|_{\mathbb{C}[g_{\infty}, h_{\infty}]} =0,
\end{equation*}
and hence the first assertion holds. The proof is finished.
\end{proof}
\end{lemma}
\autoref{lemma 4-2-1} implies that if $\lambda_{1}(\infty)<0$, then $h_{\infty}-g_{\infty}=\infty$. The following \autoref{lemma 4-2-4} will show us that if $\lambda_{1}(\infty)<0$ and $h_{\infty}-g_{\infty}=\infty$, then $h_{\infty}=-g_{\infty}=\infty$.
\begin{lemma}\label{lemma 4-2-4}
If $\lambda_{1}(\infty)<0$, then $-g_{\infty}<+\infty$ iff $h_{\infty}<+\infty$.
\begin{proof}
Without loss of generality, we suppose by contradiction that $-g_{\infty}<\infty$ and $h_{\infty}=\infty$. Noting that the principal eigenvalue $\lambda_{1}$ has translation invariance with respect to the space variable, it follows from $\lambda_{1}(\infty)<0$ that for any given sufficiently small negative number $\epsilon$, a sufficiently large $\tilde{l}$ can be found such that $\lambda_{1}(g_{\infty}+\epsilon, \tilde{l})<0$. For such $\epsilon$ and $\tilde{l}$, a sufficiently large $k_{0}$ can be found such that
\begin{equation*}
h(t)>\tilde{l}, ~~g(t)<g_{\infty}+\epsilon,~~    t\geq k_{0}\tau .
\end{equation*}
Consider the following problem%
\begin{eqnarray*}
\left\{
\begin{array}{ll}
\bar{u}_{t}=d_{1}\int\limits^{\tilde{l}}_{g_{\infty}+\epsilon}J_{1}(x-y)\bar{u}(y,t)dy-d_{1}\bar{u}-a_{11}\bar{u}+a_{12}\bar{v},\; &\,x\in(g_{\infty}+\epsilon,\tilde{l}), t\in((k\tau)^{+}, (k+1)\tau], \\[2mm]
\bar{v}_{t}=d_{2}\int\limits^{\tilde{l}}_{g_{\infty}+\epsilon}J_{2}(x-y)\bar{v}(y,t)dy-d_{2}\bar{v}-a_{22}\bar{v}+G(\bar{u}),\; &\, x\in(g_{\infty}+\epsilon,\tilde{l}),t\in((k\tau)^{+}, (k+1)\tau], \\[2mm]
\bar{u}(x,(k\tau)^{+})=H(\bar{u}(x,k\tau)), ~~~~~\bar{v}(x,(k\tau)^{+})=\bar{v}(x,k\tau),\; &\, x\in(g_{\infty}+\epsilon,\tilde{l}), \\[2mm]
\bar{v}(x,k_{0}\tau)=u(x,k_{0}\tau), ~~~~~~~~~~~~~ v(x,k_{0}\tau)=u(x,k_{0}\tau), \; &\,x\in(g_{\infty}+\epsilon,\tilde{l}), k=k_{0}, k_{0}+1,\cdots.
\end{array} \right.
\end{eqnarray*}
Then, $(u(x,t), v(x,t))\succeq(\bar{u}(x,t), \bar{u}(x,t))$ for $x\in [g_{\infty}+\epsilon,\tilde{l}]$ and $t\geq k_{0}\tau $ by \autoref{lemma 3-2-1}. Using methods similar to the proof of \autoref{lemma 4-2-1}, there exist small constant $c^{*}$ and large integer $k_{1}>k_{0}$ such that $g'(t)\leq c^{*}<0$ for $t\geq k_{1}\tau$. Therefore, we have that $-g_{\infty}=\infty$, which shows a contradiction with the assumption $-g_{\infty}<\infty$. The proof is completed.
\end{proof}
\end{lemma}
Now, we begin to analyse the dynamical behaviour of model \eqref{Zhou-Lin} when $\lambda_{1}(\infty)<0$.
\begin{lemma}\label{lemma 4-2-5}
If $\lambda_{1}(g_{\infty}, h_{\infty})<0$, then
\begin{equation*}
\lim\limits_{m\rightarrow \infty}(u(x, t+m\tau), v(x, t+m\tau))=(U_{1}(t), V_{1}(t))\text{~in~}[\mathbb{C}_{\text{loc}}(\mathds{R}\times [0,\tau])]^{2}
\end{equation*}
where $(U_{1},V_{1})$ denotes the unique solution to problem \eqref{3-24}, and $(h_{\infty}, g_{\infty})=\mathds{R}$.
\begin{proof}
In virtue of \autoref{lemma 4-2-1}, it follows from $\lambda_{1}(\infty)<0$ that $(h_{\infty}, g_{\infty})$ is an infinite interval. Then, \autoref{lemma 4-2-4} yields that $(h_{\infty}, g_{\infty})=\mathds{R}$.

Next verify the first claim by showing
\begin{equation}\label{4-3-1}
\limsup\limits_{m\rightarrow \infty}(u(x,t+m\tau), v(x,t+m\tau))\preceq(U_{1}(t), V_{1}(t))
\end{equation}
uniformly in $\mathds{R}\times [0,\tau]$ and
\begin{equation}\label{4-3-2}
\liminf\limits_{m\rightarrow \infty}(u(x,t+m\tau), v(x,t+m\tau))\succeq(U_{1}(t), V_{1}(t))
\end{equation}
locally uniformly in $[0,\tau]\times\mathds{R}$. We first prove that \eqref{4-3-1} holds. From \autoref{lemma 3-2-1}, it follows that
\begin{equation}\label{4-3-3}
\begin{aligned}
u\leq \tilde{u}\text{~and~}v\leq \tilde{v}\text{~for~} (x,t)\in \mathds{R}\times \mathds{R}^{+},
\end{aligned}
\end{equation}
where $(\tilde{u}, \tilde{v})$ is defined in the proof of \autoref{lemma 4-2-3}. Therefore, it follows from \eqref{4-3-3} that
\begin{equation}\label{4-3-4}
\begin{aligned}
\limsup\limits_{m\rightarrow \infty}(u(x,t+m\tau), v(x,t+m\tau))\preceq \lim\limits_{m\rightarrow \infty}(\tilde{u}(t+m\tau), \tilde{v}(t+m\tau))
\end{aligned}
\end{equation}
uniformly for $x\in[g(t+m\tau), h(t+m\tau)]$ and $t\in\mathds{R}^{+}$. Using  methods similar to the proof of \autoref{theorem 3-2-1}\textcolor{blue}{(2)}, we find that
\begin{equation}\label{4-3-5}
\begin{aligned}
\lim\limits_{m\rightarrow \infty}(\tilde{u}(t+m\tau), \tilde{v}(t+m\tau))=(U_{1}(t), V_{1}(t))
\end{aligned}
\end{equation}
uniformly in $[0,\tau]$. Hence, \eqref{4-3-4} and \eqref{4-3-5} yield that \eqref{4-3-1} holds.

We then prove that \eqref{4-3-2} also holds. Since $\lambda_{1}(\infty)<0$, a sufficiently large $\tilde{l}$ can be found such that $\lambda_{1}(-\tilde{l}, \tilde{l})<0$. Noticing that $h_{\infty}=-g_{\infty}=\infty$, there exists a large $k_{0}$ such that $(h(t),- g(t))\succeq (\tilde{l}, \tilde{l})$ for $t\geq k_{0}\tau$. Consider the problem
\begin{eqnarray}\label{4-3-6}
\left\{
\begin{array}{ll}
\bar{u}_{t}=d_{1}\int\limits^{\tilde{l}}_{-\tilde{l}}J_{1}(x-y)\bar{u}(y,t)dy-d_{1}\bar{u}-a_{11}\bar{u}+a_{12}\bar{v},\; &\,x\in[-\tilde{l},\tilde{l}], t\in((k\tau)^{+}, (k+1)\tau], \\[2mm]
\bar{v}_{t}=d_{2}\int\limits^{\tilde{l}}_{-\tilde{l}}J_{2}(x-y)\bar{v}(y,t)dy-d_{2}\bar{v}-a_{22}\bar{v}+G(\bar{u}),\; &\, x\in[-\tilde{l},\tilde{l}], t\in((k\tau)^{+}, (k+1)\tau], \\[2mm]
\bar{u}(x,(k\tau)^{+})=H(\bar{u}(x,k\tau)), ~~~~~\bar{v}(x,(k\tau)^{+})=\bar{v}(x,k\tau),\; &\, x\in[-\tilde{l},\tilde{l}], \\[2mm]
\bar{v}(x,k_{0}\tau)=u(x,k_{0}\tau), ~~~~~~~~~~~~~ v(x,k_{0}\tau)=u(x,k_{0}\tau), \; &\,x\in[-\tilde{l},\tilde{l}], k=k_{0}, k_{0}+1,\cdots.
\end{array} \right.
\end{eqnarray}
Then, it follows from \autoref{lemma 3-2-1} that
\begin{equation}\label{4-3-7}
u(x,t)\geq \bar{u}(x,t)\text{~and~}v(x,t)\geq \bar{v}(x,t)\text{~for~}x\in[-\tilde{l},\tilde{l}]\text{~and~}t\geq k_{0}\tau.
\end{equation}
Since $\lambda_{1}(-\tilde{l}, \tilde{l})<0$, \autoref{theorem 3-2-1}\textcolor{blue}{(2)} gives us that \eqref{4-3-6}
has a unique positive steady state $(U_{\tilde{l}}(x,t), V_{\tilde{l}}(x,t))$ defined in problem \eqref{3-2-1} with $[r,s]$ replaced by $[-\tilde{l},\tilde{l}]$, and
\begin{equation}\label{4-3-8}
\lim\limits_{m\rightarrow+\infty}\big(\bar{u}(x,t+m\tau),\bar{u}(x,t+m\tau)\big)=\big(U_{\tilde{l}}(x,t), V_{\tilde{l}}(x,t)\big)
\end{equation}
uniformly in $[0, \tau]\times [-\tilde{l},\tilde{l}]$. From \eqref{4-3-7} and \eqref{4-3-8}, it follows that
\begin{equation}\label{4-3-9}
\liminf\limits_{m\rightarrow \infty}(u(x,t+m\tau), v(x,t+m\tau))\succeq\big(U_{\tilde{l}}(x,t), V_{\tilde{l}}(x,t)\big)
\end{equation}
uniformly in $[0, \tau]\times [-\tilde{l},\tilde{l}]$. Finally, when $\tilde{l}\rightarrow \infty$, it follows from \autoref{lemma 3-2-2} that
\begin{equation*}
\liminf\limits_{m\rightarrow \infty}(u(x,t+m\tau), v(x,t+m\tau))\succeq(U_{1}(t), V_{1}(t))
\end{equation*}
locally uniformly in $\mathds{R}\times [0,\tau]$. This completes the proof.
\end{proof}
\end{lemma}
Based on \textcolor{blue}{Lemmas} \ref{lemma 4-2-2}, \ref{lemma 4-2-3} and \ref{lemma 4-2-5}, we now present the main result of this subsection.
\begin{theorem}\label{theorem 4-3-a}
The model \eqref{Zhou-Lin} is either vanishing or spreading.
\end{theorem}
\subsection{Spreading-vanishing criteria}
In this subsection, we will look for criteria guaranteeing vanishing or spreading for \eqref{Zhou-Lin}. For convenience, let $\lambda_{1}(h_{0})$ and $\lambda_{1}(\infty)$ denote the principal eigenvalues of \eqref{3-1-1} with $[r,s]$ replaced by $[-h_{0}, h_{0}]$ and $(-\infty, +\infty)$, respectively.

Based on \autoref{lemma 4-2-2} and \autoref{lemma 4-2-3}, the following conclusion can be obtained directly.
\begin{theorem}\label{theorem 4-3-b}
If $\lambda_{1}(\infty)\geq 0$, then extinction happens, that is,
\begin{equation*}
\lim\limits_{t\rightarrow\infty}\|u(x,t)+v(x,t)\|_{\mathbb{C}[g(t), h(t)]}=0.
\end{equation*}
Additionally, if $\lambda_{1}(\infty)>0$, then vanishing happens.
\end{theorem}
Next, we will consider the case of $\lambda_{1}(\infty)<0$. Its proof will be divided into two cases: $\lambda_{1}(h_{0})\leq 0$ and $\lambda_{1}(h_{0})>0$.
\begin{theorem}\label{theorem 4-3-2}
Assume that $\lambda_{1}(\infty)<0$. If $\lambda_{1}(h_{0})\leq 0$, then spreading happens.
\begin{proof}
Suppose by contradiction that $h_{\infty}<+\infty$ or $-g_{\infty}>-\infty$. Then, \autoref{lemma 4-2-4} yields that $(g_{\infty}, h_{\infty})$
is a finite interval. By using \autoref{lemma 4-2-1}, it follows that $\lambda_{1}(g_{\infty}, h_{\infty})\geq 0$. Noticing that both $-g(t)$ and $h(t)$
are strictly monotonically increasing, \autoref{lemma 3-1-4}\textcolor{blue}{(2)} gives that
\begin{equation*}
\lambda_{1}(h_{0})>\lambda_{1}(g_{\infty},h_{\infty})\geq 0.
\end{equation*}
which shows a contradiction with the condition $\lambda_{1}(h_{0})\leq 0$. This proof is completed.
\end{proof}
\end{theorem}

Next, we will explore the case of $\lambda_{1}(h_{0})>0$. To do this, we first provide the following lemma.

\begin{lemma}\label{lemma 4-3-4}
Assume that $\lambda_{1}(\infty)<0$. If $\lambda_{1}(h_{0})>0$ and the initial value $(u_{0}(x), v_{0}(x))$ of model \eqref{Zhou-Lin} is sufficiently small, then vanishing happens.
\begin{proof}
A suitable super-solution of model \eqref{Zhou-Lin} is constructed to prove this conclusion. Since $\lambda_{1}(h_{0})>0$, \autoref{lemma 3-1-4}\textcolor{blue}{(2)} gives that a positive number $\epsilon$ can be found such that
\begin{equation*}
\lambda_{1}(h_{1}):=\lambda_{1}(h_{1}, h_{1}):=\lambda_{1}(h_{0}+\epsilon,h_{0}+\epsilon)>0.
\end{equation*}
Denote by $(\zeta, \eta)$ the eigenfunction of \eqref{3-1-1} in which $[r,s]$ is substituted by $[-h_{1}, h_{1}]$. For $x\in[-h_{1}, h_{1}]$ and $t>0$, define
\begin{eqnarray*}
\begin{array}{l}
\overline{h}(t):=h_{1}-(h_{1}-h_{0})e^{-\sigma t},~~\overline{g}(t):=-\overline{h}(t), \\[2mm]
\overline{u}(x,t):=Me^{-\sigma t}\zeta(x,t), ~~~~~~~\overline{v}(x,t):=Me^{-\sigma t}\eta(x,t),
\end{array}
\end{eqnarray*}
where the nonnegative constants $M$ and $\sigma$ are determined later.

We start by examining some conditions and determining $M$ and $\sigma$ before applying \autoref{lemma 4-1-1}. Clearly,
\begin{equation*}
(\overline{u}(x,t), \overline{v}(x,t))\succeq \mathbf{0}\text{~for~}(x,t)\in \{\overline{g}(t), \overline{h}(t)\}\times (0,\infty).
\end{equation*}
A standard calculation yields that, for $x\in(\overline{g}(t), \overline{h}(t))$ and $t\in((k\tau)^{+}, (k+1)\tau]$,
\begin{equation*}
\begin{aligned}
&\overline{u}_{t}-d_{1}\int^{h(t)}_{g(t)}J_{1}(x-y)\overline{u}(y,t)dy+d_{1}\overline{u}+a_{11}\overline{u}-a_{12}\overline{v}\\
\geq&~Me^{-\sigma t}\big[-\sigma\zeta+\zeta_{t}-d_{1}\int^{h_{1}}_{-h_{1}}J_{1}(x-y)\zeta(y,t)dy+d_{1}\zeta+a_{11}\zeta-a_{12}\eta\big] \\
\geq&~Me^{-\sigma t}\big[\lambda_{1}(h_{1})\zeta-\sigma\zeta\big]\geq 0,
\end{aligned}
\end{equation*}
and
\begin{equation*}
\begin{aligned}
&\overline{v}_{t}-d_{2}\int^{h(t)}_{g(t)}J_{2}(x-y)\overline{v}(y,t)dy+d_{2}\overline{v}+a_{22}\overline{v}-G(\overline{u})\\
\geq&~Me^{-\sigma t}\big[-\sigma\eta+\eta_{t}-d_{2}\int^{h_{1}}_{-h_{1}}J_{2}(x-y)\eta(y,t)dy+d_{2}\eta+a_{22}\eta-G'(0)\zeta\big] \\
\geq&~Me^{-\sigma t}\big[\lambda_{1}(h_{1})\eta-\sigma\eta\big]\geq 0,
\end{aligned}
\end{equation*}
provided that $\sigma:=\frac{\lambda_{1}(h_{1})}{3}$. Based on assumption ($\mathbf{H}$), we have that for $x\in(\overline{g}(k\tau), \overline{h}(k\tau))$
\begin{eqnarray*}
\overline{u}(x,(k\tau)^{+})=Me^{-\sigma k\tau}\zeta(x,(k\tau)^{+})=Me^{-\sigma k\tau}H'(0)\zeta(x,k\tau)\geq H(\overline{u}(x,k\tau)),
\end{eqnarray*}
and
\begin{eqnarray*}
\overline{v}(x,(k\tau)^{+})=Me^{-\sigma k\tau}\eta(x,(k\tau)^{+})=Me^{-\sigma k\tau}\zeta(x,k\tau)=\overline{v}(x,k\tau).
\end{eqnarray*}
Additionally, assumption \textbf{(J)} yields that
\begin{equation*}
\begin{aligned}
&\mu_{1}\int\limits^{\overline{h}(t)}_{\overline{g}(t)}\int\limits^{+\infty}_{\overline{h}(t)}J_{1}(x-y)\overline{u}(x,t)dydx
+\mu_{2}\int\limits^{\overline{h}(t)}_{\overline{g}(t)}\int\limits^{+\infty}_{\overline{h}(t)}J_{2}(x-y)\overline{v}(x,t)dydx\\
\leq&~\mu_{1}\int\limits^{\overline{h}(t)}_{\overline{g}(t)}\overline{u}(x,t)dx+\mu_{2}\int\limits^{\overline{h}(t)}_{\overline{g}(t)}\overline{v}(x,t)dx \\
\leq&~Me^{-\sigma t}\Big[\mu_{1}\int\limits^{h_{1}}_{-h_{1}}\max\limits_{[0,\tau]\times [-h_{1}, h_{1}]}\zeta(x,t)dx
+\mu_{2}\int\limits^{h_{1}}_{-h_{1}}\max\limits_{[0,\tau]\times [-h_{1}, h_{1}]}\eta(x,t)dx\Big] \\
\leq &~(h_{1}-h_{0})\sigma e^{-\sigma t}=\overline{h}'(t)\text{~~for~}t>0,
\end{aligned}
\end{equation*}
provided that
\begin{eqnarray*}
M:=\frac{(h_{1}-h_{0})\sigma}{\mu_{1}\int\limits^{h_{1}}_{-h_{1}}\max\limits_{[-h_{1}, h_{1}]\times [0,\tau]}\zeta(x,t)dx+\mu_{2}\int\limits^{h_{1}}_{-h_{1}}\max\limits_{[-h_{1}, h_{1}]\times [0,\tau]}\eta(x,t)dx}.
\end{eqnarray*}
Similarly, we have that
\begin{eqnarray*}
 -\mu_{1}\int\limits^{h(t)}_{g(t)}\int\limits^{g(t)}_{-\infty}J_{1}(x-y)\overline{u}(x,t)dydx-\mu_{2}\int\limits^{h(t)}_{g(t)}\int\limits^{g(t)}_{-\infty}J_{2}(x-y)\overline{v}(x,t)dydx\geq \overline{g}'(t)\text{~~for~}t>0.
\end{eqnarray*}
For the initial value, if we choose $(u_{0}(x), v_{0}(x))$ sufficiently small such that
\begin{eqnarray*}
\Big(\|u_{0}(x)\|_{\mathbb{C}[-h_{0}, h_{0}]}, \|v_{0}(x)\|_{\mathbb{C}[-h_{0}, h_{0}]}\Big)\preceq M\Big(\min\limits_{[-h_{0}, h_{0}]}\zeta(x,0), \min\limits_{ [-h_{0}, h_{0}]}\eta(x,0)\Big) ,
\end{eqnarray*}
then
\begin{eqnarray*}
\big(u_{0}(x), v_{0}(x)\big)\preceq \big(\overline{u}(x,0), \overline{v}(x,0)\big)\text{~~for~}x\in[-h_{0}, h_{0}].
\end{eqnarray*}

Now it is clear that the quadruple $(\overline{u}, \overline{v}, \overline{g}, \overline{h})$ is a super-solution of \eqref{Zhou-Lin}. Therefore, by \autoref{lemma 4-1-1}, it follows that
\begin{eqnarray*}
(g_{\infty}, h_{\infty})\subset \Big(\lim\limits_{t\rightarrow+\infty}\overline{g}(t), \lim\limits_{t\rightarrow+\infty}\overline{h}(t) \Big)\subset(-h_{1},h_{1}),
\end{eqnarray*}
which implies that $(g_{\infty}, h_{\infty})$ is a finite interval. Therefore, \autoref{lemma 4-2-1} yields that vanishing happens. This ends the proof.
\end{proof}
\end{lemma}
\begin{remark}\label{remark 4-3-1}
In the proof of \autoref{lemma 4-3-4}, it follows from observing the definition of $M$ that when $\mu_{2}=\rho\mu_{1}$, where $\rho$ is a positive constant, we have that $M\rightarrow\infty$ as $\mu_{1}\rightarrow 0$.
It means that, if $\lambda_{1}(\infty)<0$, $\lambda_{1}(h_{0})>0$, and $\mu_{2}=\rho\mu_{1}$, for any $(u_{0}(x), v_{0}(x))$, a constant $\underline{\mu}>0$ can be found such that when $\mu_{1}\in(0, \underline{\mu}]$ vanishing happens.
\end{remark}
\begin{lemma}\label{lemma 4-3-6}
Assume that $\lambda_{1}(\infty)<0$. If $\lambda_{1}(h_{0})>0$, then there exists $l^{*}>0$ such that $\lambda_{1}(g_{\infty}, h_{\infty})<0$ if $h_{\infty}-g_{\infty}>2l^{*}$.
\begin{proof}
Since $\lambda_{1}(\infty)<0$ and $\lambda_{1}(h_{0})>0$, it follows from \autoref{lemma 3-1-4}\textcolor{blue}{(2)}
that there exists a constant $l^{*}>0$ such that $\lambda_{1}(-l^{*}, l^{*})=0$. Again using \autoref{lemma 3-1-4}\textcolor{blue}{(2)} yields that when $h_{\infty}-g_{\infty}>2l^{*}$, $\lambda_{1}(g_{\infty}, h_{\infty})<0$. The proof is completed.
\end{proof}
\end{lemma}
\begin{lemma}\label{lemma 4-3-3}
Suppose that $\lambda_{1}(\infty)<0$. If $\lambda_{1}(h_{0})>0$, then a positive constant $\overline{\mu}$ can be found such that if $\mu_{1}+\mu_{2}>\overline{\mu}$
spreading occurs.
\begin{proof}
Denote the unique solution of model \eqref{Zhou-Lin} by $(u^{\bm{\mu}}, v^{\bm{\mu}}, g^{\bm{\mu}}, h^{\bm{\mu}})$ in order to
emphasise the dependency on $\bm{\mu}:=(\mu_{1}, \mu_{2})$. First of all, we claim that $\bm{\tilde{\mu}}\succ\bm{0}$ can be found such that for some large $t_{0}>0$,
\begin{eqnarray}\label{4-3-10}
\lambda_{1}(g^{{\bm{\tilde{\mu}}}}(t_{0}), h^{{\bm{\tilde{\mu}}}}(t_{0}))<0.
\end{eqnarray}
If not, suppose that for any $\bm{\mu}\succ\bm{0}$ and $t>0$, $\lambda_{1}(g^{\bm{\mu}}(t), h^{\bm{\mu}}(t))\geq0$. We will derive a contradiction.

By \autoref{theorem 2-1}, we have that $-g^{\bm{\mu}}(t)$ and $h^{\bm{\mu}}(t)$ are strongly increasing with respect to $t>0$. Additionally, it follows from \autoref{lemma 4-1-1} that $-g^{\bm{\mu}}(t)$ and $h^{\bm{\mu}}(t)$ are nondecreasing with respect to $\bm{\mu}\succ \bm{0}$. Hence,
\begin{eqnarray*}
G_{\infty}:=\lim\limits_{\bm{\mu}\rightarrow\infty}\lim\limits_{t\rightarrow\infty} g^{\bm{\mu}}(t)\text{~~and~~}H_{\infty}:=\lim\limits_{\bm{\mu}\rightarrow\infty}\lim\limits_{t\rightarrow\infty} h^{\bm{\mu}}(t)
\end{eqnarray*}
are well defined. By \autoref{lemma 4-3-6}, it follows that $H_{\infty}-G_{\infty}\leq2l^{*}$. Assumption \textbf{(J)} yields that the sufficiently small $\epsilon_{0}>0$ and $\delta_{0}>0$ can be found such that
\begin{eqnarray*}
J_{i}(x)\geq \delta_{0}>0\text{~for~}x\in[-\epsilon_{0}, \epsilon_{0}]\text{~and~}i=1,2.
\end{eqnarray*}
Then for the fixed $\epsilon_{0}$, we can find $\bm{\mu}_{0}\succ \bm{0}$ and $k_{0}\in \mathds{N}^{+}$ such that
\begin{eqnarray*}
h^{\bm{\mu}}(t)+\frac{\epsilon_{0}}{4}>H_{\infty}\text{~for~}\bm{\mu}\succeq\bm{\mu}_{0}\text{~and~}t\geq k_{0}\tau.
\end{eqnarray*}
By integrating both sides of the fourth equation in the model over $[k_{0}\tau, (k_{0}+1)\tau]$, it follows that
\begin{equation*}
\begin{aligned}
&h^{\bm{\mu}}((k_{0}+1)\tau)-h^{\bm{\mu}}(k_{0}\tau)\\
&=\int_{k_{0}\tau}^{(k_{0}+1)\tau}\Bigg[ \mu_{1}\int\limits^{h^{\bm{\mu}}(t)}_{g^{\bm{\mu}}(t)}\int\limits^{+\infty}_{h^{\bm{\mu}}(t)}J_{1}(x-y)u(x,t)dydx
+\mu_{2}\int\limits^{h^{\bm{\mu}}(t)}_{g^{\bm{\mu}}(t)}\int\limits^{+\infty}_{h^{\bm{\mu}}(t)}J_{2}(x-y)v(x,t)dydx\Bigg]dt\\
&\geq\int_{k_{0}\tau}^{(k_{0}+1)\tau}\Bigg[ \mu_{1}\int\limits^{h^{\bm{\mu}_{0}}(t)}_{g^{\bm{\mu}_{0}}(t)}\int\limits^{+\infty}_{h^{\bm{\mu}_{0}}(t)
+\frac{\epsilon_{0}}{4}}J_{1}(x-y)u(x,t)dydx+\mu_{2}\int\limits^{h^{\bm{\mu}_{0}}(t)}_{g^{\bm{\mu}_{0}}(t)}\int\limits^{+\infty}_{h^{\bm{\mu}_{0}}(t)
+\frac{\epsilon_{0}}{4}}J_{2}(x-y)v(x,t)dydx\Bigg]dt\\
&\geq\int_{k_{0}\tau}^{(k_{0}+1)\tau}\Bigg[ \mu_{1}\int\limits^{h^{\bm{\mu}_{0}}(t)}_{h^{\bm{\mu}_{0}}(t)-\frac{\epsilon_{0}}{2}}\int\limits^{h^{\bm{\mu}_{0}}(t)
+\frac{\epsilon_{0}}{2}}_{h^{\bm{\mu}_{0}}(t)+\frac{\epsilon_{0}}{4}}J_{1}(x-y)u(x,t)dydx
+\mu_{2}\int\limits^{h^{\bm{\mu}_{0}}(t)}_{h^{\bm{\mu}_{0}}(t)-\frac{\epsilon_{0}}{2}}\int\limits^{h^{\bm{\mu}_{0}}(t)
+\frac{\epsilon_{0}}{2}}_{h^{\bm{\mu}_{0}}(t)+\frac{\epsilon_{0}}{4}}J_{2}(x-y)v(x,t)dydx\Bigg]dt\\
&\geq\frac{\epsilon_{0}\delta_{0}}{4}\int_{k_{0}\tau}^{(k_{0}+1)\tau}\Bigg[ \mu_{1}\int\limits^{h^{\bm{\mu}_{0}}(t)}_{h^{\bm{\mu}_{0}}(t)-\frac{\epsilon_{0}}{2}}u(x,t)dx
+\mu_{2}\int\limits^{h^{\bm{\mu}_{0}}(t)}_{h^{\bm{\mu}_{0}}(t)-\frac{\epsilon_{0}}{2}}v(x,t)dx\Bigg]dt\\
&\geq \frac{(\mu_{1}+\mu_{2})\delta_{0}\tau A\epsilon^{2}_{0}}{8},
\end{aligned}
\end{equation*}
where
\begin{eqnarray*}
A:=\min\Bigg\{\min\limits_{[h^{\bm{\mu}_{0}}(t)-\frac{\epsilon_{0}}{2}, h^{\bm{\mu}_{0}}(t)] \times [k_{0}\tau, (k_{0}+1)\tau] }u(x,t), \min\limits_{[h^{\bm{\mu}_{0}}(t)-\frac{\epsilon_{0}}{2}, h^{\bm{\mu}_{0}}(t)]\times [k_{0}\tau, (k_{0}+1)\tau]}v(x,t)\Bigg\}.
\end{eqnarray*}
This implies that
\begin{eqnarray*}
\mu_{1}+\mu_{2}\leq \frac{8(h^{\bm{\mu}}((k_{0}+1)\tau)-h^{\bm{\mu}}(k_{0}\tau))}{\delta_{0}\tau A\epsilon^{2}_{0}}\leq  \frac{32l^{*}}{\delta_{0}\tau A\epsilon^{2}_{0}}<\infty.
\end{eqnarray*}
which shows a contradiction with the assumption. Therefore, $\eqref{4-3-10}$ holds. Finally, it follows from \autoref{lemma 4-2-5} that the conclusion holds. The proof is completed.
\end{proof}
\end{lemma}
Now we consider the case $\lambda_{1}(h_{0})>0$.
\begin{theorem}\label{theorem 4-3-3}
Assume that $\lambda_{1}(\infty)<0$. If $\lambda_{1}(h_{0})>0$, then there exist $\underline{\mu}>0$ and $\overline{\mu}>0$ such that vanishing occurs when $\mu_{1}+\mu_{2}\leq \underline{\mu}$ and spreading happens when $\mu_{1}+\mu_{2}>\overline{\mu}$. Additionally, if $\mu_{2}=\rho\mu_{1}$ where $\rho$ is a positive constant, then a positive constant $\mu^{*}$ can be found such that vanishing occurs when $\mu_{1}\leq \mu^{*}$ and spreading happens when $\mu_{1}>\mu^{*}$.
\begin{proof}
We first prove the first assertion. In the proof of \autoref{lemma 4-3-4}, we have from observing the definition of $M$ that as $\mu_{1}+\mu_{2}\rightarrow 0$ $M\rightarrow\infty$ . This implies that for any given pair of initial function $(u_{0}(x), v_{0}(x))$ of model \eqref{Zhou-Lin}, a sufficiently small positive constant $\underline{\mu}$ can be found such that when $\mu_{1}+\mu_{2}\leq \underline{\mu}$ vanishing occurs. The remainder of the first conclusion is a direct result of \autoref{lemma 4-3-3}. For the proof of the second assertion, it can be obtained by a method similar to that in \cite[Theorem 3.14]{Cao-Du-Li-Li}, and we are omitting it here. This proof is completed.
\end{proof}
\end{theorem}
\section{\bf Numerical simulation}\label{Section-5}
In this section, we will present a numerical example to demonstrate the theoretical outcomes and further understand the effect of the impulse intervention on the spread of the faecal-oral diseases.

\begin{figure}[!httb]
  \centering
  \includegraphics[scale=0.5]{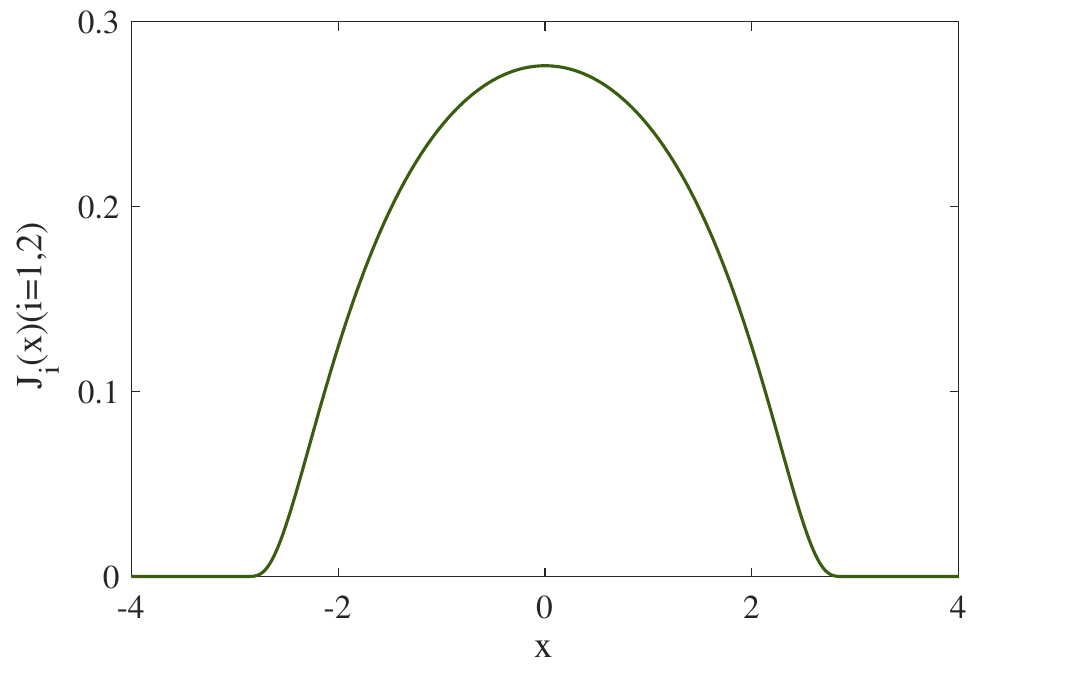}
    \caption{Dispersal kernel function $J_{i}(x)(i=1,2)$ defined by \eqref{5-1}.}\label{figure 1}
\end{figure}

Since the relationship between the infectious agents and the infected individuals is cooperative, we only give the numerical simulation of the infectious agents in order to save space. Moreover, due to the reason that the exact $h(t)$ of model \eqref{Zhou-Lin} cannot be easily obtained, we regard the simulated $h(t)$ as the replacement of the exact $h(t)$ to make the estimation of $\lambda_{1}$ easier. We choose the dispersal kernel functions as follows
\begin{eqnarray}\label{5-1}
J_{i}(x)=
\left\{
\begin{array}{l}
ke^{\frac{1}{ |\frac{x}{3}|^{2}-1}}, ~~|x|<3, \\[2mm]
0,~~~~~~~~~~~|x|\geq 3,
\end{array} \right.
\end{eqnarray}
where $i=1,2$ and $\frac{1}{k}=\int_{-\infty}^{+\infty}e^{\frac{1}{ |\frac{x}{3}|^{2}-1}}dx$. The diagram for $J_{i}(x)(i=1,2)$ is shown as \autoref{figure 1}.
The initial infection region is chosen as $[-2,2]$ and the initial functions are as follows
\begin{equation*}
u_{0}(x)=3\cos\Big(\frac{\pi x}{4}\Big),~v_{0}(x)=\cos\Big(\frac{\pi x}{4}\Big), ~x\in[-2, 2].
\end{equation*}
The other parameters of model \eqref{Zhou-Lin} are given in the following example.
\begin{exm}\label{exm1}
Fix $d_{1}=0.10$, $d_{2}=0.10$, $\tau=1$, $a_{11}=0.35$, $a_{12}=0.11$, $a_{22}=0.10$, $\mu_{1}=20$, $\mu_{2}=200$, and $G(u)=0.5u/(10+u)$. $H(u)$ is chosen as $u$ and $\frac{0.1u}{10+u}$, respectively.
\end{exm}

\begin{figure}[!ht]
\centering
\subfigure{ {
\includegraphics[width=0.40\textwidth]{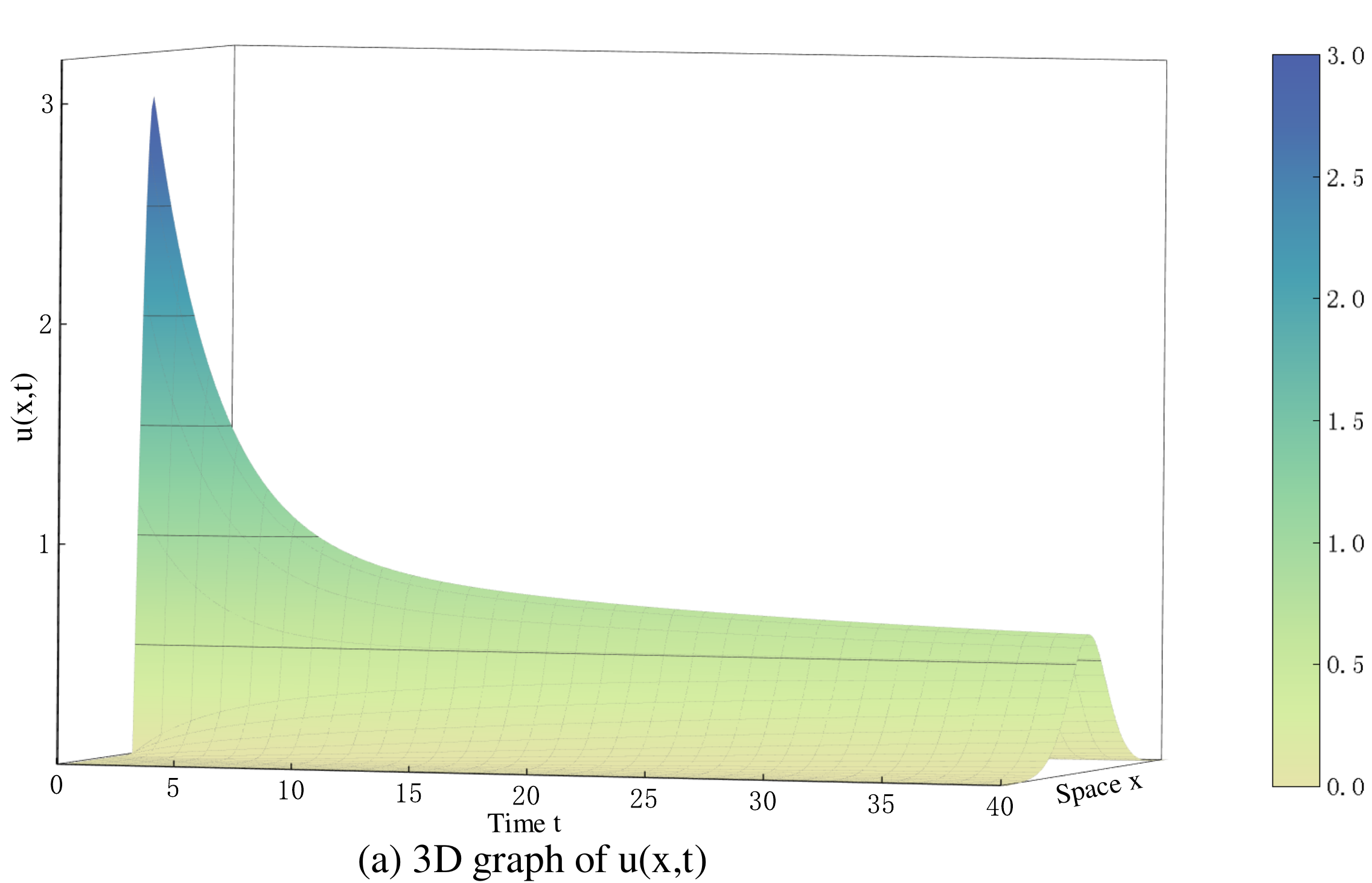}
} }
\subfigure{ {
\includegraphics[width=0.40\textwidth]{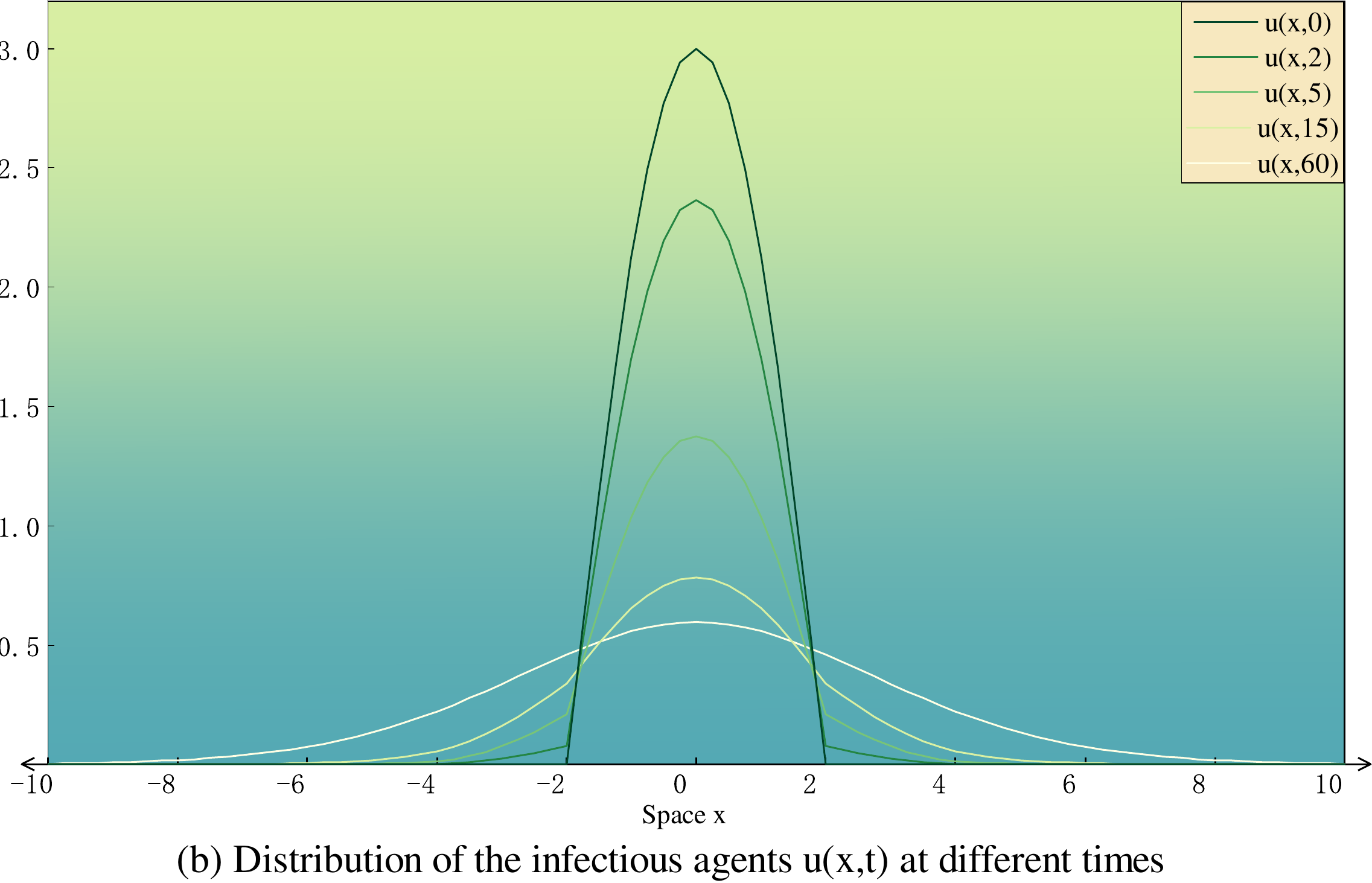}
} }
\subfigure{ {
\includegraphics[width=0.40\textwidth]{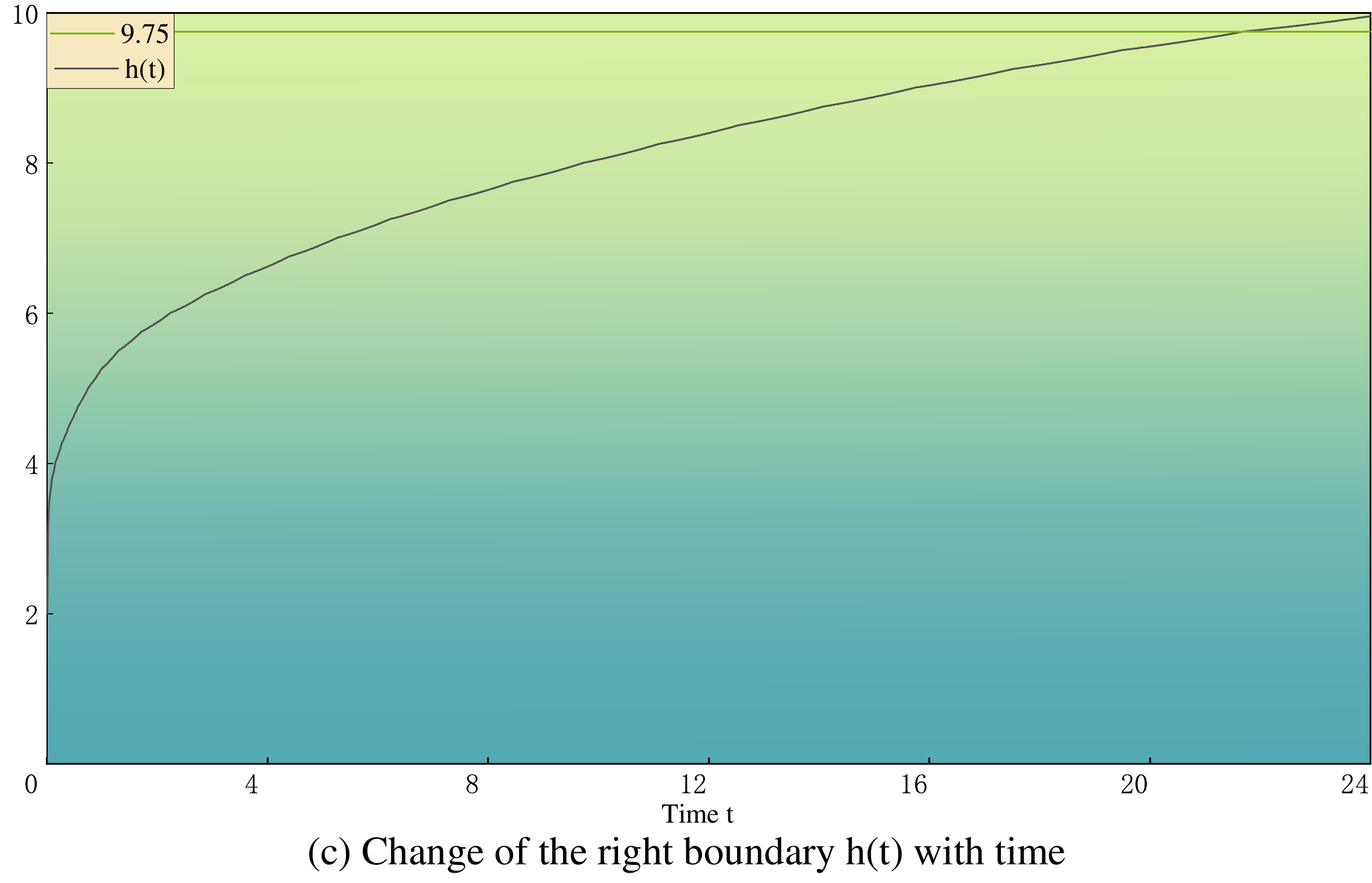}
} }
\subfigure{ {
\includegraphics[width=0.4\textwidth]{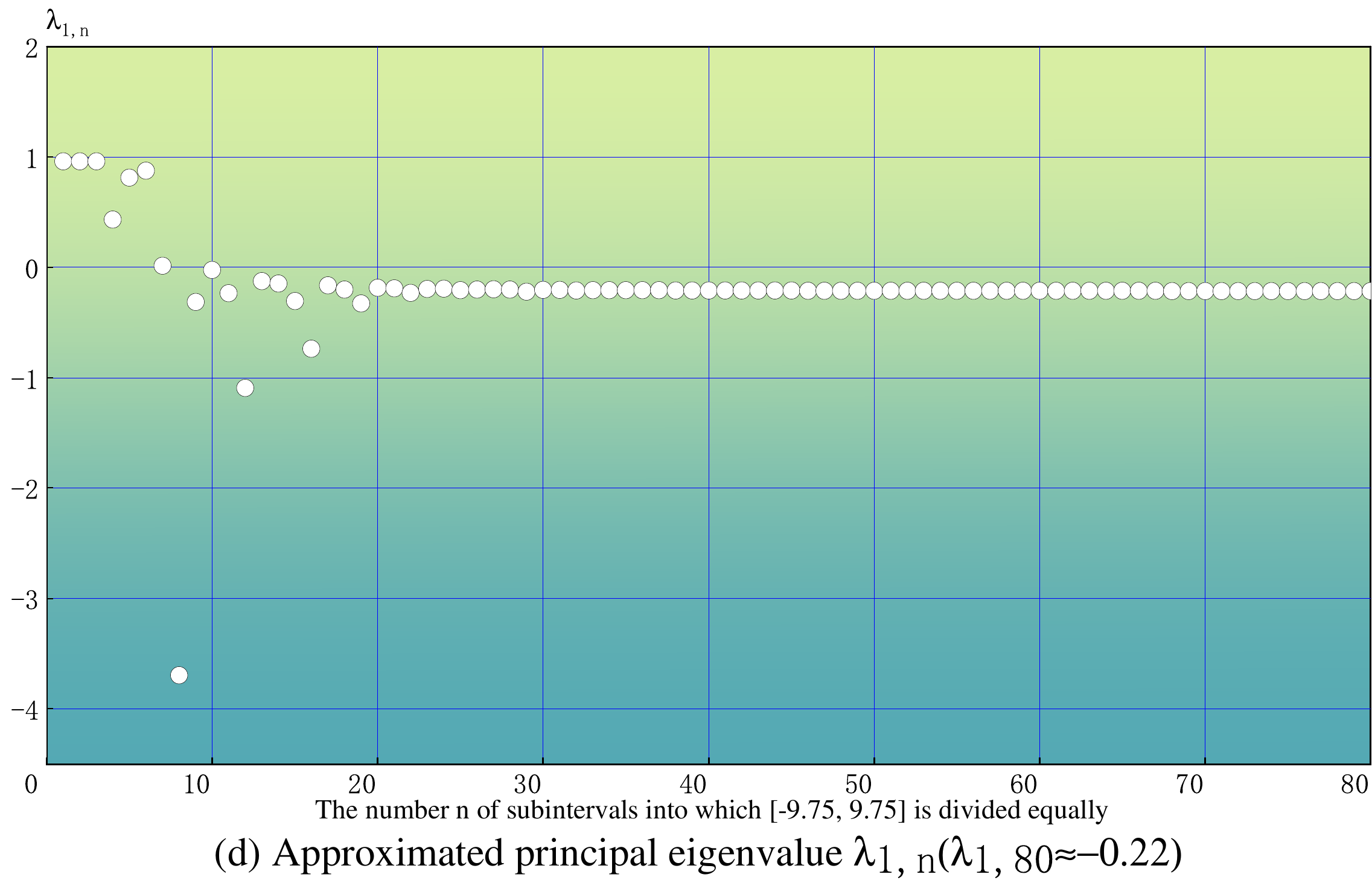}
} }
\caption{When $H(u)=u$, which implies no pulse, $u$ converges to a steady state and $\lambda_{1}<0$.}
\label{A}
\end{figure}
\begin{figure}[!ht]
\centering
\subfigure{ {
\includegraphics[width=0.40\textwidth]{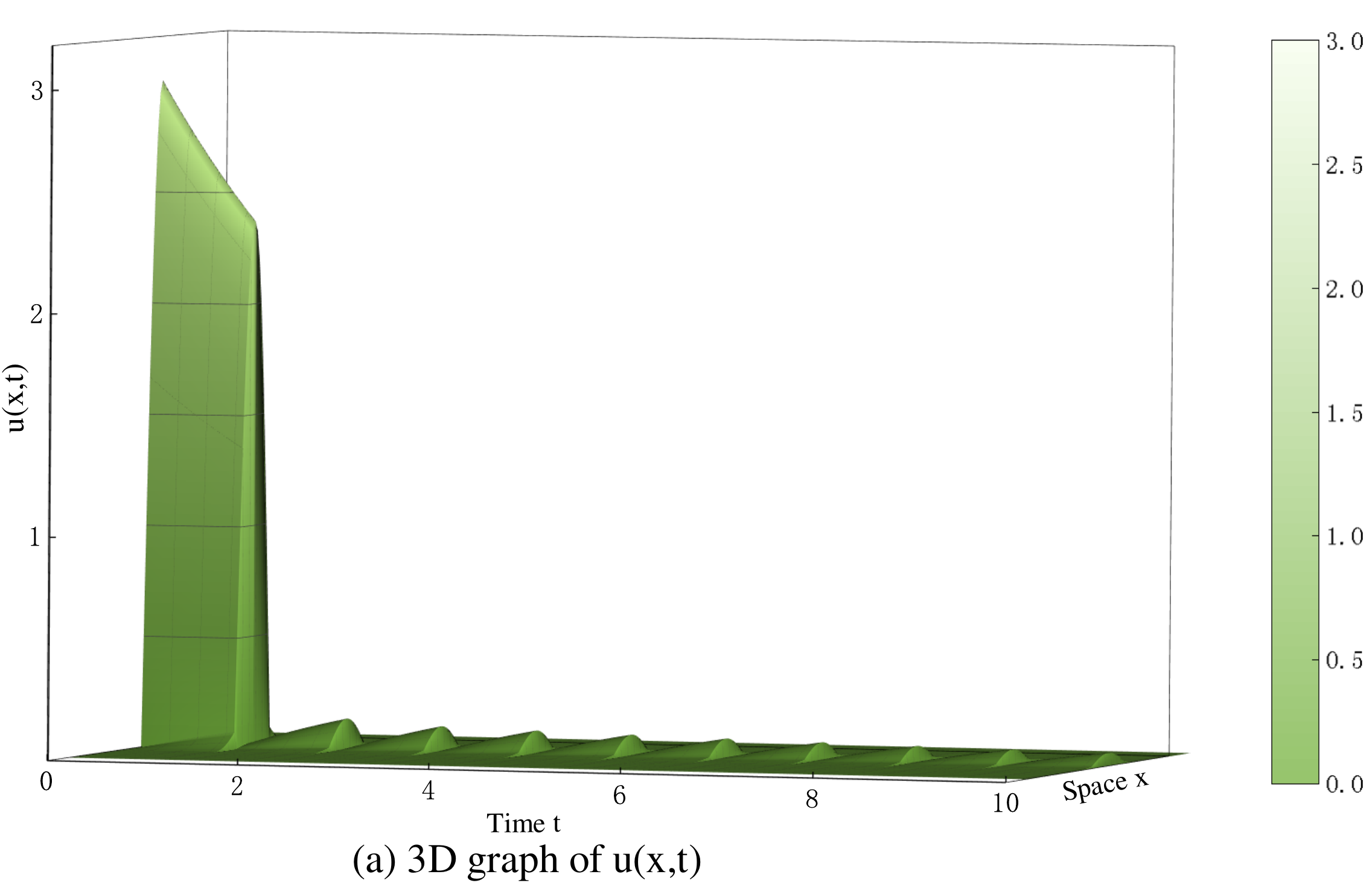}
} }
\subfigure{ {
\includegraphics[width=0.40\textwidth]{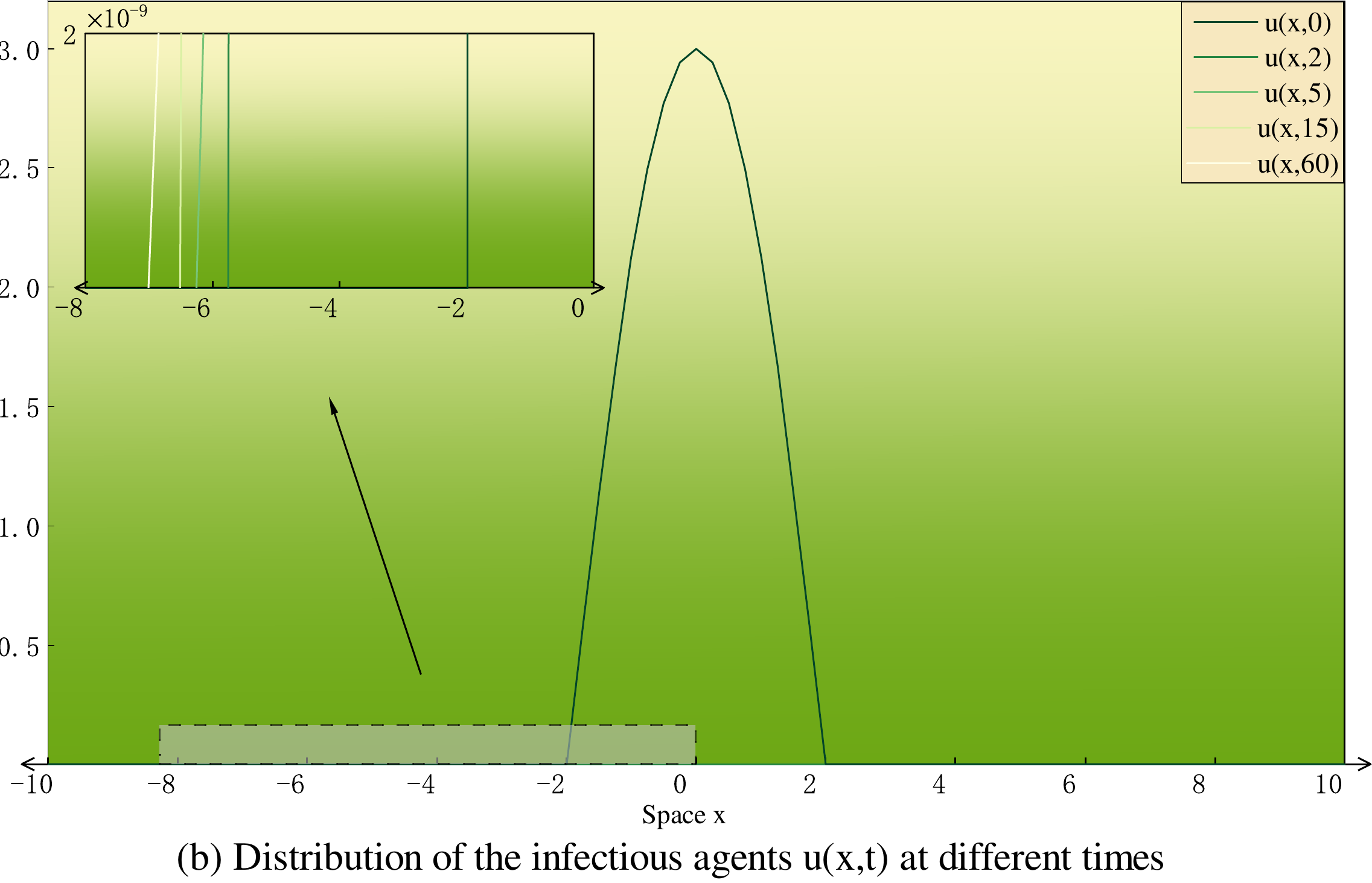}
} }
\subfigure{ {
\includegraphics[width=0.40\textwidth]{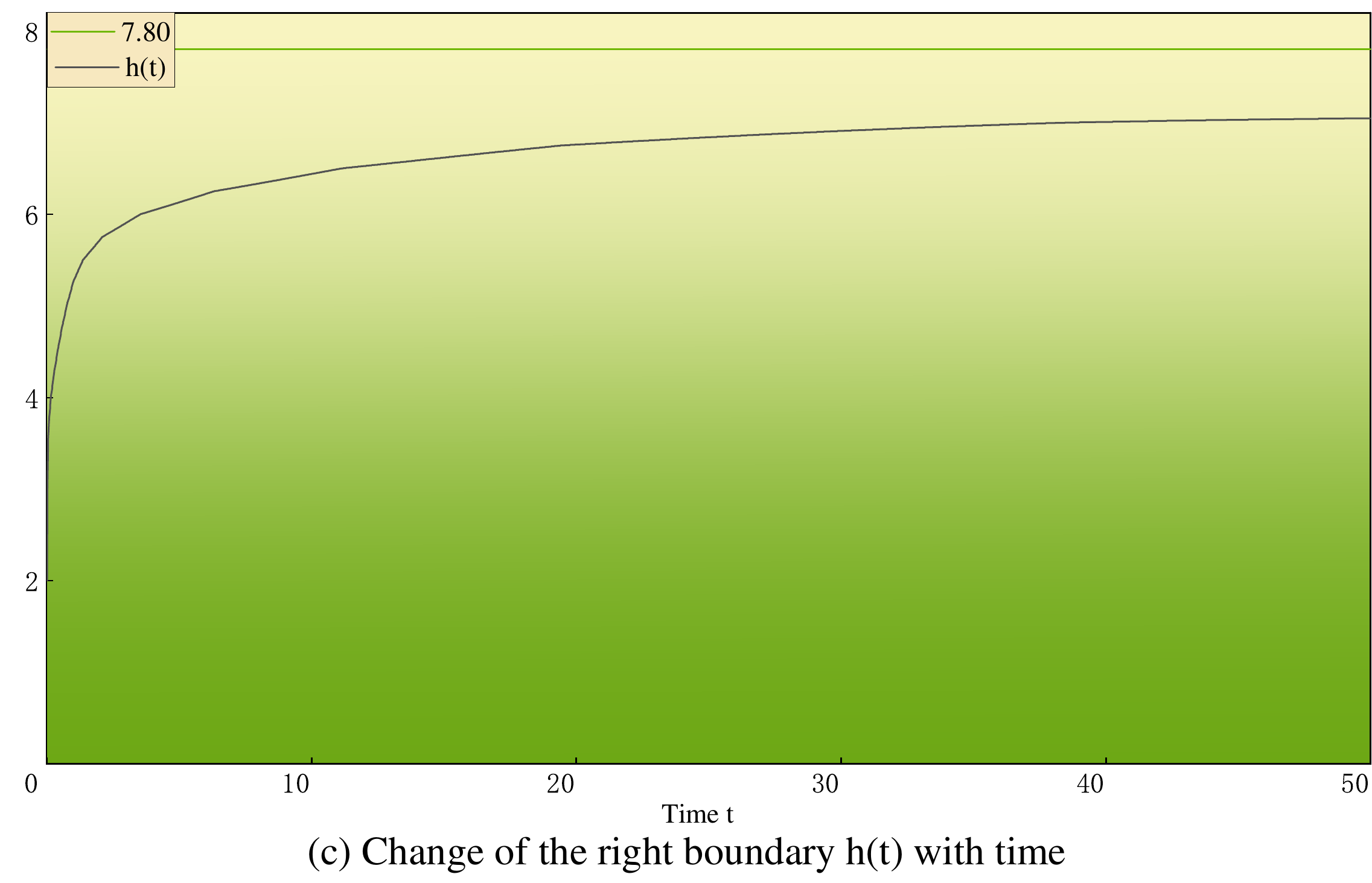}
} }
\subfigure{ {
\includegraphics[width=0.4\textwidth]{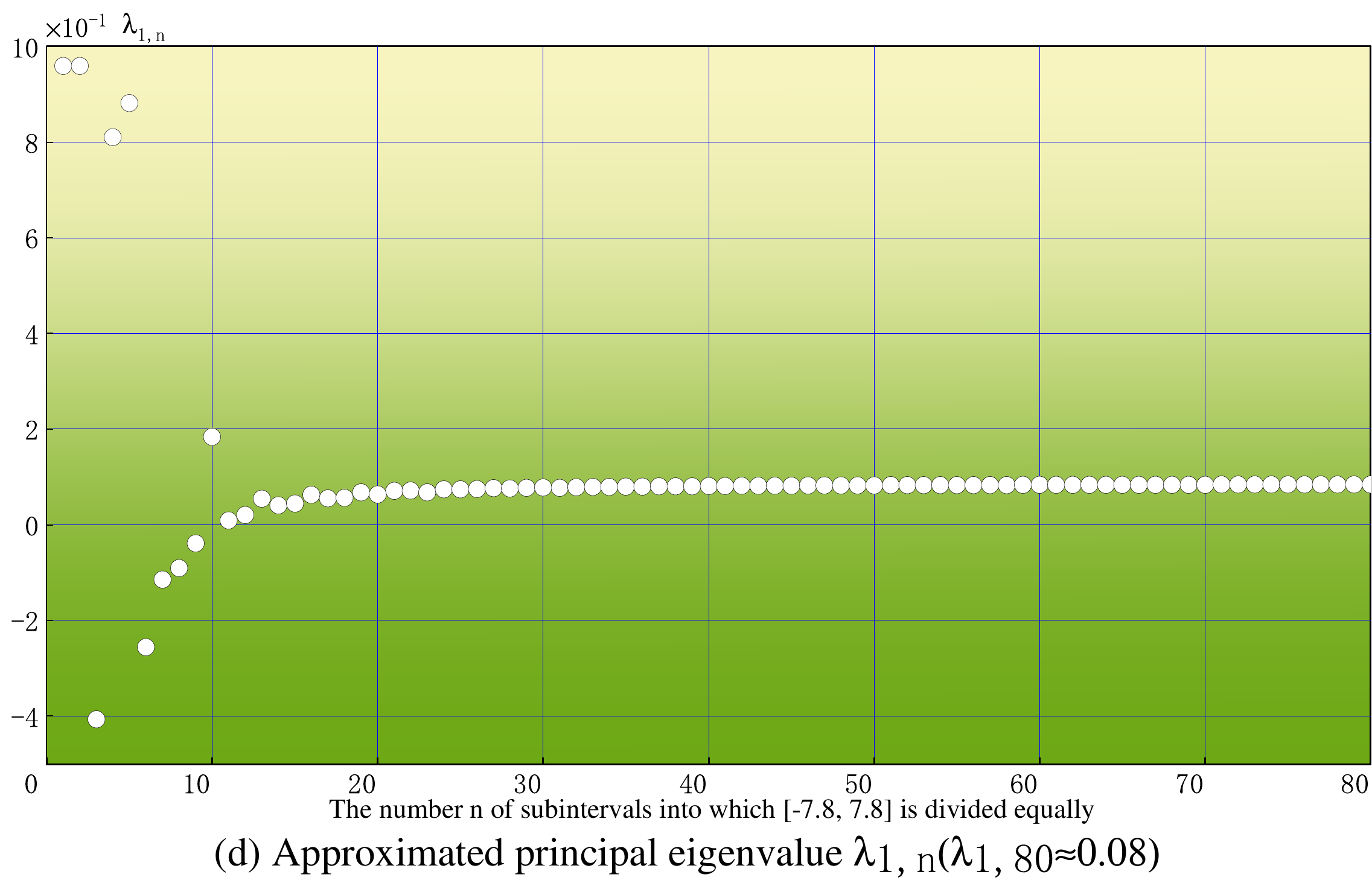}
} }
\caption{When $H(u)=\frac{0.1u}{10+u}$, $u$ decays to $0$ and $\lambda_{1}>0$.}
\label{B}
\end{figure}

When there is no the impulsive intervention, it can be seen from \autoref{A}\textcolor[rgb]{0.00,0.00,1.00}{(c)} that $h_{\infty}>h(24)>9.75$. By using the theory in \cite{Chatelin} and the method in \cite{Kuniya-Wang, Wang-Yang}, we obtain \autoref{A}\textcolor[rgb]{0.00,0.00,1.00}{(d)}. This combined with \autoref{lemma 3-1-4}\textcolor{blue}{(2)} yields that
\begin{equation*}
\lambda_{1}(-g_{\infty}, h_{\infty})<\lambda_{1}(-9.75, 9.75)\approx -0.22<0.
\end{equation*}
Then, it follows from \autoref{lemma 4-2-5} that $(U_{1}(t), V_{1}(t))$ is globally asymptotically stable. Actually, it follows from \autoref{A}\textcolor[rgb]{0.00,0.00,1.00}{(a)} that the infectious agents indeed converge to a heterogeneous stable state, which is consistent with the conclusion of \autoref{lemma 4-2-5}.

When the impulse function is taken as $H(u)=\frac{0.1u}{10+u}$, one can obtain from \autoref{B}\textcolor[rgb]{0.00,0.00,1.00}{(c)} that $h_{\infty}<7.8$. Similarly, \autoref{B}\textcolor[rgb]{0.00,0.00,1.00}{(d)} can be obtained. This and \autoref{lemma 3-1-4} tell us that
\begin{equation*}
\lambda_{1}(H'(0), (-g_{\infty}, h_{\infty}))= \lambda_{1}(0.01, (-g_{\infty}, h_{\infty}))>\lambda_{1}(1, (-7.8, 7.8))\approx 0.08>0.
\end{equation*}
Then, we have from \autoref{lemma 4-2-2} that the disease-free equilibrium point $(0,0)$ is globally asymptotically stable. Actually, it is not hard to see from \autoref{B}\textcolor[rgb]{0.00,0.00,1.00}{(a)} that the infectious agents continuously converge to zero as time $t$ increases, which is consistent with the conclusion obtained.

It can be obtained from comparing \autoref{A} and \autoref{B} that the implementation of the impulsive intervention can reduce the number of the infectious agents, slow down the expansion speed of the diseases, decrease the infected region of the diseases, and extinguish the diseases which are originally persistent. Therefore, the impulsive intervention is beneficial in combating the diseases transmitted by the faecal-oral route.
\section{\bf Conclusion and future work}\label{Section-6}
In this paper, we develop an impulsive nonlocal model to capture the effect of the combination of  periodic spraying of disinfectant and  nonlocal diffusion of individuals on the diseases transmitted by the faecal-oral route. Our model extends the model of chang and Du in~\cite{Chang-Du}, which does not take into account the pulse intervention. With the introduction of the pulse intervention, the model is not continuous at the pulse point and the corresponding steady state is governed by a periodic parabolic problem, hence the techniques used here are rather different from those in~\cite{Chang-Du}.

Based on the work of Wang and Du in~\cite{Wang-Du-1}, \autoref{theorem 2-1} shows that the model has a unique nonnegative global classical solution. Then, the corresponding periodic eigenvalue problem \eqref{3-1-1} is presented in order to study the asymptotical profiles of the model. The existence and uniqueness of the principal eigenvalue of problem \eqref{3-1-1} is proved in \autoref{theorem 3-1} by using the theory of resolvent positive operators with their perturbations. With the help of this eigenvalue, \autoref{theorem 4-3-a} shows that model \eqref{Zhou-Lin} is either vanishing or spreading. To provide the criteria for determining spreading and vanishing, the principal eigenvalues of the eigenvalue problem \eqref{3-1-1} with $[r,s]$ replaced by $[-h_{0}, h_{0}]$ and $(-\infty, +\infty)$ are next defined in $\lambda_{1}(h_{0})$ and $\lambda_{1}(\infty)$, respectively. We obtain the following criteria:
\begin{itemize}
\item{If $\lambda_{1}(\infty)>0$, then the diseases are vanishing, see \autoref{lemma 4-2-2};}
\item{If $\lambda_{1}(\infty)=0$, then the diseases are vanishing, see \autoref{lemma 4-2-3};}
\item{If $\lambda_{1}(\infty)< 0$ and $\lambda_{1}(h_{0})\leq 0$, then the diseases are spreading, see \autoref{theorem 4-3-2};}
\item{If $\lambda_{1}(\infty)< 0$ and $\lambda_{1}(h_{0})>0$, then
   \begin{itemize}
     \item[$\bullet$]{there exist two positive constants $\overline{\mu}$ and $\underline{\mu}$ such that vanishing occurs when $\mu_{1}+\mu_{2}\leq \underline{\mu}$ and spreading happens when $\mu_{1}+\mu_{2}>\overline{\mu}$, see \autoref{theorem 4-3-3},}
      \item[$\bullet$]{or if $\mu_{2}=\rho\mu_{1}$ where $\rho$ is a positive constant, then a positive constant $\mu^{*}$ can be found such that vanishing occurs when $\mu_{1}\leq \mu^{*}$ and spreading happens when $\mu_{1}>\mu^{*}$, see \autoref{theorem 4-3-3}.}
    \end{itemize}}
\end{itemize}

Our theoretical results suggest that the implementation of the impulse intervention is beneficial in combating the diseases transmitted by the faecal-oral route. Numerical simulation is used to further understand this finding. Additionally, the initial infected region also has a significant impact on the evolution of the diseases. Specifically, the smaller the initial infected area, the more advantageous for human beings to fight against the diseases. The effect of the nonlocal diffusion of individuals on the evolution of the diseases depends on the choice of the diffusion kernel functions.

When $\lambda_{1}(\infty)=0$, the conclusion that the diseases are vanishing is obtained under the assumption that the length of the infected region is finite. We strongly believe that this assumption can be removed. This will be the direction of our future research.

\section*{\bf Availability of data and material}
The authors confirm that the data supporting the findings of this study are available within the manuscript.

\section*{\bf Competing interests}
Authors have no conflict of interest.




\begin{thebibliography}{99}
\bibitem{Ail-Nelson-Lopze}
M. Ali, A.R. Nelson, A.L. Lopez, et al., Updated global burden of cholera in endemic countries, PLoS Negl Trop Dis, 2015, 9(6):e0003832.

\bibitem{Mogasale-Maskery}
V. Mogasale, B. Maskery, R.L. Ochiai, et al., Burden of typhoid fever in low-income and middle-income countries: a systematic, literature-based update with risk-factor adjustment,
Lancet Glob Health, 2014, 2(10):570-580.

\bibitem{Hay-collaborators}
S.I. Hay, A.E. Schumacher, H.H. Kyu, et al., Global age-sex-specific mortality, life expectancy, and population estimates in 204 countries and territories and 811 subnational locations, 1950-2021, and the impact of the COVID-19 pandemic: a comprehensive demographic analysis for the Global Burden of Disease Study 2021, Lancet, 2024, 403(10440):1989-2056.

\bibitem{Chen-Tang-Teng}
Q.L. Chen, S.Y. Tang, Z.D. Teng, et al., Long-time dynamics and semi-wave of a delayed nonlocal epidemic model with free boundaries, Proc. Roy. Soc. Edinburgh Sect. A, 2023, 1-47.

\bibitem{Capasso-Maddalena-JMB}
V. Capasso, L. Maddalena, Convergence to equilibrium states for a reaction-diffusion system modeling the spatial spread of a class of bacterial and viral diseases, J. Math. Biol., 1981, 13:173-184.

\bibitem{Heller-Mota}
L. Heller, C.R. Mota, D.B. Greco, COVID-19 faecal-oral transmission: Are we asking the right questions? Sci. Total Environ, 2020, 729:138919.

\bibitem{Ghosh-Kumar-Santiana}
S. Ghosh, M. Kumar, M. Santiana, et al., Enteric viruses replicate in salivary glands and infect through saliva, Nature, 2022, 607:345-350.

\bibitem{Yang-Gong-Sun}
J.Y. Yang, M.J. Gong, G.Q. Sun, Asymptotical profiles of an age-structured foot-and-mouth disease with nonlocal diffusion on a spatially heterogeneous environment, J. Differential Equations, 2023, 377:71-112.

\bibitem{Capasso-Paveri-Fontana}
V. Capasso, S.L. Paveri-Fontana, A mathematical model for the 1973 cholera epidemic in the European Mediterranean region, Rev. Epidemiol. Sante. Publique., 1979, 27:121-132.

\bibitem{Capasso-JMAA}
V. Capasso, Asymptotic stability for an integro-differential reaction-diffusion system, J. Math. Anal. Appl., 1984, 103:575-588.

\bibitem{Thieme-Zhao-JDE}
H.R. Thieme, X.Q. Zhao, Asymptotic speeds of spread and traveling waves for integral equations and delayed reaction-diffusion models, J. Differential Equations, 2003, 195:430-470.

\bibitem{Wu-Hsu-TAMS}
S.L. Wu, C.H. Hsu, Existence of entire solutions for delayed monostable epidemic models, Trans. Am. Math. Soc., 2016, 368:6033-6062.

\bibitem{Hsu-Yang-Nonlinearity}
C.H. Hsu, T.S. Yang, Existence, uniqueness, monotonicity and asymptotic behaviour of travelling waves for epidemic models, Nonlinearity, 2013, 26(1):121-139.

\bibitem{Du-Lin}
Y.H. Du, Z.G. Lin, Spreading-vanishing dichotomy in the diffusive logistic model with a free boundary, SIAM J. Math. Anal., 2010, 42(1):377-405.

\bibitem{Ahn-Baek-Lin}
I. Ahn, S. Baek, Z.G. Lin, The spreading fronts of an infective environment in a man-environment-man epidemic model, Appl. Math. Modelling, 2016, 40(15-16):7082-7101.

\bibitem{Wang-Du-3}
R. Wang, Y.H. Du, Long-time dynamics of a diffusive epidemic model with free boundaries, Discrete Contin. Dyn. Syst. Ser. B, 2021, 26(4):2201-2238.

\bibitem{Zhao-Ruan}
G.Y. Zhao, S.G. Ruan, Spatiotemporal dynamics in epidemic models with L\'{e}vy flights: A fractional diffusion approach, J Math Pures Appl, 2023, 173:243-277.

\bibitem{Murray}
J.D. Murray, Mathematical Biology II: Spatial Models and Biomedical Applications, third edition, Springer-Verlag, New York, 2003.

\bibitem{Cao-Du-Li-Li}
J.F. Cao, Y.H. Du, F. Li, et al., The dynamics of a Fisher-KPP nonlocal diffusion model with free boundaries, J. Funct. Anal., 2019, 277:2772-2814.

\bibitem{Wang-Du-1}
R. Wang, Y.H. Du, Long-time dynamics of a nonlocal epidemic model with free boundaries: Spreading-vanishing dichotomy, J. Differential Equations, 2022, 327:322-381.

\bibitem{Zhao-Zhang-Li}
M. Zhao, Y. Zhang, W.T. Li, et al., The dynamics of a degenerate epidemic model with nonlocal diffusion and free boundaries, J. Differential Equations, 2020, 269:3347-3386.

\bibitem{Wang-Du-2}
R. Wang, Y.H. Du, Long-time dynamics of an epidemic model with nonlocal diffusion and free boundaries: Spreading speed, Discrete Contin. Dynam. Systems, 2023, 43(1):121-161.

\bibitem{Du-Ni}
Y.H. Du, W.J. Ni, Spreading speed for monostable cooperative systems with nonlocal diffusion and free boundaries, part 1: semi-wave and a threshold condition,
J. Differential Equations, 2022, 308:369-420.

\bibitem{Zhao-Li-Du}
M. Zhao, W.T. Li, Y.H. Du, The effect of nonlocal reaction in an epidemic model with nonlocal diffusion and free boundaries, Commun. Pure Appl. Anal., 2020, 19(9):4599-4620.

\bibitem{Du-Li-Ni}
Y.H. Du, W.T. Li, W.J. Ni, et al., Finite or infinite spreading speed of an epidemic model with free boundary and double nonlocal effects, J. Dynam. Differential Equations, 2024, 36:1015-1063.

\bibitem{Chang-Du}
T.Y. Chang, Y.H. Du, Long-time dynamics of an epidemic model with nonlocal diffusion and free boundaries, Electron Res Arch, 2021, 30(1):289-313.

\bibitem{Lewis-Li}
M.A. Lewis, B.T. Li, Spreading speed, traveling waves, and minimal domain size in impulsive reaction-diffusion models, Bull. Math. Biol., 2012, 74:2383-2402.

\bibitem{Fazly-Lewis-Wang}
M. Fazly, M. Lewis, H. Wang, On impulsive reaction-diffusion models in higher dimensions, SIAM J. Appl. Math., 2017, 77(1):224-246.

\bibitem{Zhou-Lin-Santos}
Q. Zhou, Z.G. Lin, C.A. Santos, On an impulsive faecal-oral model in a periodically evolving environment, Chaos Solitons Fractals, 2025, 191:115825.

\bibitem{Liang-Yan-Xiang}
J.H. Liang, Q. Yan, C.C. Xiang, et al., A reaction-diffusion population growth equation with multiple pulse perturbations, Commun. Nonlinear Sci Numer. Simul., 2019, 74:122-137.

\bibitem{Fazly-Lewis-Wang-2}
M. Fazly, M. Lewi, H. Wang, Analysis of propagation for impulsive reaction-diffusion models, SIAM J. Appl. Math., 2020, 80(1):521-542.

\bibitem{Wu-Zhao-1}
R.W. Wu, X.Q. Zhao, Spatial invasion of a birth pulse population with nonlocal dispersal, SIAM J. Appl. Math., 2019, 79(3):1075-1097.

\bibitem{Zhang-Yi-Chen}
Y.R. Zhang, T.S. Yi, Y.M. Chen, Spreading dynamics of an impulsive reaction-diffusion model with shifting environments, J. Differential Equations, 2024, 381:1-19.

\bibitem{Wu-Zhao-2}
R.W. Wu, X.Q. Zhao, The evolution dynamics of an impulsive hybrid population model with spatial heterogeneity, Commun. Nonlinear Sci. Numer. Simul., 2022, 107:106181.

\bibitem{Liang-Zhang-Zhao}
X. Liang, L. Zhang, X.Q. Zhao, The principal eigenvalue for degenerate periodic reaction-diffusion systems, SIAM J. Math. Anal., 2017, 49:3603-3636.

\bibitem{Bao-Shen}
X. Bao, W. Shen, Criteria for the existence of principal eigenvalues of time periodic cooperative linear systems with nonlocal dispersal, Proc. Amer. Math. Soc., 2017, 145:2881-2894.

\bibitem{Feng-Li-Ruan-Xin}
Y.X. Feng, W.T. Li, S.G. Ruan, et al., Principal spectral theory of time-periodic nonlocal dispersal cooperative systems and applications, SIAM J. Math. Anal., 2024, 56:4040-4083.

\bibitem{Thieme}
H.R. Thieme, Spectral bound and reproduction number for infinite-dimensional population structure and time heterogeneity, SIAM J. Appl. Math., 2009, 70:188-211.

\bibitem{Thieme-1}
H.R. Thieme, Remarks on resolvent positive operators and their perturbation, Discrete Contin. Dynam. Systems, 1998, 4:73-90.

\bibitem{Brezis}
H. Brezis, Functional Analysis, Sobolev Spaces and Partial Differential Equations, Springer, New York, 2011.

\bibitem{Meyre-Nieberg}
P. Meyre-Nieberg, Banach Lattices, Springer-Verlag, New York, 1991.



\bibitem{Zhou-Lin-Pedersen}
Q. Zhou, Z.G. Lin, M. Pedersen, On an impulsive faecal-oral model in a moving infected environment, preprint, arXiv:2410.15798, 2024.

\bibitem{Chatelin}
F. Chatelin, The spectral approximation of linear operators with applications to the computation of eigenelements of differential and integral operators, SIAM Rev., 1981, 23(4):495-522.

\bibitem{Kuniya-Wang}
T. Kuniya, J.L. Wang, Global dynamics of an SIR epidemic model with nonlocal diffusion, Nonlinear Anal. Real World Appl., 2018, 43:262-282.

\bibitem{Wang-Yang}
X.Y. Wang, J.Y. Yang, Dynamics of a nonlocal dispersal foot-and-mouth disease model in a spatially heterogeneous environment, Acta Math. Sci., 2021, 41:552-572.


\end{thebibliography}
\end{document}